\documentclass[preprint,12pt]{elsarticle}

\usepackage{amssymb}
\usepackage{amsmath}
\usepackage{bm}
\usepackage{color}
\usepackage{caption}
\usepackage{exscale}
\usepackage{relsize}
\usepackage{graphicx}
\usepackage{subcaption}
\usepackage{wrapfig}
\usepackage{multicol}
\usepackage{verbatim}
\newtheorem{thm}{Theorem}
\newtheorem{rmk}{Remark}
\newtheorem{lem}{Lemma}

\usepackage{mathrsfs}

\journal{Journal of Computational Physics}

\begin{document}

\begin{frontmatter}



\title{A fourth-order multi-scale computational method and its convergence analysis for composite
Kirchhoff plates with microscopic periodic configurations}


\author[label1]{Hao Dong\corref{cor1}}\ead{donghao@mail.nwpu.edu.cn}
\cortext[cor1]{Corresponding author.}
\author[label2]{Liqun Cao}
\address[label1]{School of Mathematics and Statistics, Xidian University, Xi'an 710071, PR China}
%
\address[label2]{LSEC, ICMSEC, Academy of Mathematics and Systems Science, Chinese Academy of Sciences, Beijing 100190, PR China}
\begin{abstract}
The Kirchhoff plate model plays a vital role in modeling, computing and analyzing the mechanical behaviors of thin plate structures. This study propose a novel fourth-order multi-scale (FOMS) computational method for high-accuracy and efficient simulation of composite Kirchhoff plates with highly periodic heterogeneities.
At first, two-scale asymptotic expansion theory is employed to establish the high-accuracy fourth-order multi-scale computation model with novel fourth-order correctors for composite Kirchhoff plates, which are governed by fourth-order partial differential equation (PDE) with periodically oscillatory and highly discontinuous coefficients. Then, the locally point-wise error analysis is derived to theoretically illustrate the local balance preserving of fourth-order multi-scale model enabling high-accuracy multi-scale computation. Furthermore, a global error estimation with an explicit order for fourth-order multi-scale solutions is first demonstrated under appropriate assumptions. In contrast to the second- and third-order multi-scale solutions, only the fourth-order one is capable of providing an explicit error order estimate. Additionally, an efficient numerical algorithm is developed to conduct high-accuracy simulation for heterogeneous plate structures. Extensive numerical examples are provided to confirm the theoretical results for the computational convergence and accuracy of the proposed method. This work offers a higher-order (fourth-order) multi-scale computational framework that enables robust simulation and high-accuracy analysis to composite Kirchhoff plates.
\end{abstract}


\end{frontmatter}



\section{Introduction}
The tremendous progress in composite materials has been a key driver of the rapid development of modern science and technology. Engineering structures made of composite materials are widely used in engineering fields such as aerospace, aviation, marine, civil engineering, and architecture, etc. Particularly noteworthy are composite thin plate structures, whose mechanical properties-lightweight, high strength, and high reliability-have made them a focus of increased attention and application in the past few decades \cite{R1}. However, highly heterogeneous micro-structures of composite thin plates lead to tremendous cost to the high-accuracy modeling and simulation by use of classical numerical methods. Therefore, the development of high-accuracy and efficient numerical method is of paramount importance for effectively simulating the mechanical behaviors of composite thin plates.

In recent decades, scientists and engineers have developed a diverse array of multi-scale methods to overcome the challenging issues inherent in the multi-scale spatial nature of composite structures, such as asymptotic homogenization method (AHM) \cite{R2,R3}, multi-scale finite element method (MsFEM) \cite{R4,R5}, heterogeneous multi-scale method (HMM) \cite{R6}, variational multi-scale method (VMS) \cite{R7}, multi-scale eigenelement method (MEM) \cite{R8}, localized orthogonal decomposition method (LOD) \cite{R9} and their improved versions \cite{R10,R11}, among which the AHM is a quintessential approach with rigorous mathematical foundation and can balance the accuracy and efficiency of computational performance. Nevertheless, Cui and his research team revealed, via both theoretical analysis and engineering computation, the inadequate accuracy of classical AHM for highly heterogeneous composites. In the past thirty years, Cui and his research team systematically established a class of higher-order multi-scale approaches for precisely and efficiently simulating the thermal, mechanical and multiphysics behaviors of composite structures, as shown in references \cite{R12,R13,R14,R15,R16}. However, most of multi-scale approaches was primarily confined to simulate and analyze the two-dimensional (2D) or three-dimensional (3D) multi-scale problems of composite structures, which are governed by second-order partial differential equation in spatial scale. These multi-scale methods can not directly apply to simulate the bending behaviors of composite thin plates with highly periodic heterogeneities, which are modeled by fourth-order partial differential equation named as Kirchhoff plate model \cite{R1}. The Kirchhoff plate is a classical plate model in which only deflection is the dependent variable, which is suitable for analyzing and simulating the thin plates, and widely employed in various engineering fields. Hence, the development of effective multi-scale computational methods for composite thin plates holds significant scientific importance and practical engineering value.

To the best of our knowledge, some progress has been achieved regarding the problems of composite thin plates, both in mathematical theory and numerical computation. From a mathematics theory perspective, Pastukhova established G-convergence theory for higher-order elliptic operators without considering boundary condition \cite{R17}. Pastukhova employed the shift method to prove the operator error estimates for fourth-order elliptic equations in the unbounded domain without considering boundary condition \cite{R18}. Moreover, researcher derived the operator error estimates of the Dirichlet and Neumann problems for higher-order elliptic equations with periodic coefficients in the bounded domain \cite{R19,R20}. In addition, Niu et al. established the quantitative error estimates in homogenization and almost-periodic of higher-order elliptic systems \cite{R21,R22,R23}. Nevertheless, the aforementioned theoretical findings are not established on the basis of the plate model in practical engineering applications, and without effective numerical algorithms. In recent years, Xing et al. developed a two-scale asymptotic homogenization method for composite Kirchhoff plates with periodic micro-structures and obtained the two-scale convergence results for the two-scale asymptotic expansion solution with second-order corrector \cite{R24,R25}. However, the proposed homogenization method fails to capture the highly oscillatory behavior in the composite Kirchhoff plates with high fidelity. And also, the proposed two-scale convergence results are deficient in that they do not offer a quantitative error estimate. From a computational method perspective, some analytical approaches have been developed to obtain the effective mechanical parameters of composite plates \cite{R1,R26,R27,R28}. However, these analytical approaches are not suitable for composite plates with complicated microscopic configurations. Additionally, researchers developed a novel numerical implementation of asymptotic homogenization method to predict effective properties of periodic plates without any complicated mathematical derivation \cite{R29}. In reference \cite{R30}, a consistent multi-scale formulation is presented for the bending analysis of heterogeneous thin plate structures containing three dimensional reinforcements with in-plane periodicity. Thus a unified numerical implementation for thin plate analysis can be conveniently realized using the triangular elements with discretization flexibility and the quadratic Hermite triangular element with improved accuracy. Moreover, based on the theoretical framework of the extended multi-scale finite element method, an efficient multi-scale finite element method is developed for small-deflection analysis of thin composite plates with aperiodic microstructure characteristics \cite{R31}. Aliabadi et al. presented a novel approach to concurrent multi-scale analysis, which employed a plate model to formulate structures at both scales, and homogenisation was performed using the FFT-based approach, offering higher efficiency compared to conventional methods \cite{R32}. However, the aforementioned proposed approaches in references \cite{R24,R29,R30,R31,R32} are limited in their inability to implement high-accuracy simulation of composite plates with high-contrast material properties, and they are not equipped with a quantitative and explicit error estimate.

In this work, an innovative fourth-order multi-scale method is presented for high-accuracy and efficient simulation of fourth-order composite Kirchhoff plate problems. The pivotal contributions of this work are summarized as three pillars. (1) The fourth-order multi-scale computational model with novel fourth-order correctors is presented, which can preserve the local balance of transverse loads enabling high-resolution capturing of locally high-frequency oscillatory behavior in composite Kirchhoff plates. (2) The explicit error estimate for fourth-order multi-scale solution for fourth-order multi-scale PDE is derived for the first time, that provides a quantitative convergence rate for fourth-order multi-scale solution. (3) The efficient numerical algorithm is developed for effective multi-scale simulation of composite plate bending problems, thereby providing a robust numerical framework for real-world engineering applications.

The organization of this study is as follows. In Section 2, the establishment of fourth-order multi-scale computational model for composite Kirchhoff plates is introduced in detail. Section 3 conduct both local and global error analyses for multi-scale asymptotic solutions in the point-wise and energy senses. In particular, the explicit error order for fourth-order multi-scale problems is obtained for the first time. In Section 4, a high-accuracy and efficient numerical algorithm is developed to effectively simulate the bending problems of Kirchhoff plates with microscopic periodic configurations. In Section 5, numerical experiments are carried out to verify the computational performance of the proposed FOMS method including FOMS model and corresponding numerical algorithm. Eventually, concluding remarks and potential directions are drawn in Section 6. Throughout this paper, Einstein summation convention is utilized to streamline repetitive indices.
\section{Fourth-order multi-scale computational model for composite Kirchhoff plates}
\subsection{The mathematical setting of composite Kirchhoff plates}
Based on the classical plate theory in \cite{R1,R37,R40}, the following fourth-order partial differential equation with rapidly oscillating and highly discontinuous coefficients is presented for composite Kirchhoff plates with spatial periodic heterogeneity.
\begin{equation}
	\left\{  \begin{aligned}
		&   \frac{\partial^2 }{{\partial {x_i}\partial {x_j}}}\bigg[ {D_{ijkl}^\epsilon (\bm{x})\frac{{{\partial ^2 {\omega^\epsilon }(\bm{x})}}}{{\partial {x_k}\partial {x_l}}}} \bigg] {\rm{ = }}{q}(\bm{x}), \;\; \text{in}\;\; \Omega,  \\
		& {\omega^\epsilon }(\bm{x})= g_1(\bm{x}), \;\; \text{on}\;\; \partial \Omega,\\
		&\frac{{\partial {\omega^\epsilon }(\bm{x})}}{{\partial {x_j}}}{\hat n_j} = g_2(\bm{x}), \;\; \text{on}\;\; \partial \Omega,
	\end{aligned} \right.
\end{equation}
where $\Omega$ denotes a bounded convex domain in $\mathbb{R}^2$ with a Lipschitz continuous boundary $\partial\Omega$. $D_{ijkl}^\epsilon(\bm{x})$ is the fourth-order bending stiffness tensor. ${\omega^\epsilon }(\bm{x})$ denotes the undetermined transverse displacement on the reference middle surface of composite Kirchhoff plates, where the superscript $\epsilon$ represents the characteristic length of periodic unit cell (PUC) in composite Kirchhoff plates. ${q}(\bm{x})$ denotes the transverse load imposed on composite Kirchhoff plates. $g_1(\bm{x})$ and $g_2(\bm{x})$ denote the prescribed transverse displacement and the prescribed in-plane rotation on domain boundary $\Omega$. Besides, without loss of generality, $\bm{n}=(\hat n_1,\hat n_2)$ denotes the outward normal vector on the outer boundary $\partial \Omega$.

In accordance with multi-scale homogenization theory, we can define microscopic variable $\bm{y}=\bm{x}/\epsilon=(x_1/\epsilon,x_2/\epsilon)=(y_1,y_2)$ for periodic unit cell $\mathbf{Y}=[0,1]^2$. Consequently, material parameter $D_{ijkl}^\epsilon(\bm{x})$ can be formulated as new form $D_{ijkl}(\bm{y})$, which necessitate this material parameter is 1-periodic function in variable $\bm{y}$. To implement the subsequent multi-scale modeling and theoretical analysis, this section starts with the following assumptions for the above fourth-order governing PDE (1).
\begin{enumerate}
	\item[(A$_1$)]
	 The bending stiffness coefficient $D_{ijkl}^\epsilon (\bm{x}) \in {L^\infty }(\Omega )$.
	\item[(A$_2$)]
	Function $D_{ijkl}^\epsilon (\bm{x})$ is symmetric and uniform ellipticity, and there exist two constants $0 < \alpha  \le \beta $ such that,
	\begin{displaymath}
		\begin{aligned}
			&D_{ijkl}^\epsilon  = D_{ijlk}^\epsilon  = D_{klij}^\epsilon,\;\alpha {\eta _{ij}}{\eta _{ij}} \le D_{ijkl}^\epsilon (\bm{x}){\eta _{ij}}{\eta _{kl}} \le \beta {\eta _{ij}}{\eta _{ij}},
		\end{aligned}
	\end{displaymath}
	for arbitrary symmetric matrix $\{ {\eta _{ij}}\}  \in {\mathbb{R}^{2\times 2}}$.
	\item[(A$_3$)]
	Function ${q}(\bm{x}) \in {H^{-2}}(\Omega)$, function ${g}_1(\bm{x}) \in {H^{3/2}}(\partial\Omega)$, function ${g}_2(\bm{x}) \in {H^{1/2}}(\partial\Omega)$.
\end{enumerate}
\subsection{Higher-order multi-scale analysis for composite Kirchhoff plates}
Firstly, taking into account the microscopic periodicity and scale separation of composite plates, the following multi-scale chain rules can be determined for subsequent multi-scale modeling.
\begin{equation}
	\frac{\partial \Phi^\epsilon (\bm{x})}{{\partial {x_i}}} = \frac{\partial \Phi(\bm{x},\bm{y})}{{\partial {x_i}}} + \frac{1}{\epsilon }\frac{\partial \Phi(\bm{x},\bm{y})}{{\partial {y_i}}},
\end{equation}
\begin{equation}
	\frac{{{\partial ^2\Phi^\epsilon (\bm{x})}}}{{\partial {x_i}\partial {x_j}}} = \frac{{{\partial ^2\Phi(\bm{x},\bm{y})}}}{{\partial {x_i}\partial {x_j}}} + {\frac{1}{\epsilon }}\frac{{{\partial ^2\Phi(\bm{x},\bm{y})}}}{{\partial {x_i}\partial {y_j}}} + {\frac{1}{\epsilon }}\frac{{{\partial ^2\Phi(\bm{x},\bm{y})}}}{{\partial {y_i}\partial {x_j}}} + {\frac{1}{\epsilon^2 }}\frac{{{\partial ^2\Phi(\bm{x},\bm{y})}}}{{\partial {y_i}\partial {y_j}}}.
\end{equation}
To clarify and simplify subsequent multi-scale analysis, a series of differential operators is defined according to the aforementioned chain rules as follows.
\begin{equation}
	\left\{  \begin{aligned}
		& {B_0}\omega^\epsilon (\bm{x}) = \frac{\partial^2 }{{\partial {y_i}\partial {y_j}}}\bigg[ {{D_{ijkl}}(\bm{y})\frac{{{\partial ^2 \omega^\epsilon (\bm{x})}}}{{\partial {y_k}\partial {y_l}}}} \bigg],  \\
		&{B_1}\omega^\epsilon (\bm{x}) = \!2\frac{\partial^2 }{{\partial {x_i}\partial {y_j}}}\bigg[ {{D_{ijkl}}(\bm{y})\frac{{{\partial ^2 \omega^\epsilon (\bm{x})}}}{{\partial {y_k}\partial {y_l}}}} \bigg]\!+\!2\frac{\partial^2 }{{\partial {y_i}\partial {y_j}}}\bigg[ {{D_{ijkl}}(\bm{y})\frac{{{\partial ^2 \omega^\epsilon (\bm{x})}}}{{\partial {x_k}\partial {y_l}}}} \bigg],\\
		&{B_2}\omega^\epsilon (\bm{x})= \frac{\partial^2 }{{\partial {x_i}\partial {x_j}}}\bigg[ {{D_{ijkl}}(\bm{y})\frac{{{\partial ^2 \omega^\epsilon (\bm{x})}}}{{\partial {y_k}\partial {y_l}}}} \bigg]+\frac{\partial^2 }{{\partial {y_i}\partial {y_j}}}\bigg[ {{D_{ijkl}}(\bm{y})\frac{{{\partial ^2 \omega^\epsilon (\bm{x})}}}{{\partial {x_k}\partial {x_l}}}} \bigg]\\
&\quad\quad\quad\;\;+4\frac{\partial^2 }{{\partial {x_i}\partial {y_j}}}\bigg[ {{D_{ijkl}}(\bm{y})\frac{{{\partial ^2 \omega^\epsilon (\bm{x})}}}{{\partial {x_k}\partial {y_l}}}} \bigg],\\
		&{B_3}\omega^\epsilon (\bm{x})=\!2\frac{\partial^2 }{{\partial {y_i}\partial {x_j}}}\bigg[ {{D_{ijkl}}(\bm{y})\frac{{{\partial ^2 \omega^\epsilon (\bm{x})}}}{{\partial {x_k}\partial {x_l}}}} \bigg]\!+\!2\frac{\partial^2 }{{\partial {x_i}\partial {x_j}}}\bigg[ {{D_{ijkl}}(\bm{y})\frac{{{\partial ^2 \omega^\epsilon (\bm{x})}}}{{\partial {y_k}\partial {x_l}}}} \bigg], \\
		&{B_4}\omega^\epsilon (\bm{x}) = \frac{\partial^2 }{{\partial {x_i}\partial {x_j}}}\bigg[ {{D_{ijkl}}(\bm{y})\frac{{{\partial ^2 \omega^\epsilon (\bm{x})}}}{{\partial {x_k}\partial {x_l}}}} \bigg].
	\end{aligned} \right.
\end{equation}

Next, assume that the transverse displacement solution $\omega^\epsilon (\bm{x})$ of the fourth-order PDE (1) has the following asymptotic expansion form.
\begin{equation}
		\left\{  \begin{aligned}
		& \omega^\epsilon (\bm{x}) = \omega^{(0)}(\bm{x},\bm{y}) + \epsilon \omega^{(1)}(\bm{x},\bm{y}) + {\epsilon ^2}\omega^{(2)}(\bm{x},\bm{y})  \\
		&\qquad \;\;+ {\epsilon ^3}\omega^{(3)}(\bm{x},\bm{y}) + {\epsilon ^4}\omega^{(4)}(\bm{x},\bm{y}) + {\rm O}({\epsilon^5}).
	\end{aligned} \right.
\end{equation}
And also, substituting multi-scale asymptotic form (5) into multi-scale PDE (1), and utilizing the chain rules and differential operators provided in (2)-(4) gives the following series of equations by grouping the power-like terms of small periodic parameter $\epsilon$.
\begin{equation}
	O( {{\epsilon ^{ - 4}}} ):{B_0}\omega^{(0)}(\bm{x},\bm{y}) = 0.
\end{equation}
\begin{equation}
	O( {{\epsilon ^{ - 3}}}):{B_0}\omega^{(1)}(\bm{x},\bm{y}) + {B_1}\omega^{(0)}(\bm{x},\bm{y}) = 0.
\end{equation}
\begin{equation}
	O( {{\epsilon ^{ - 2}}}): {B_0}\omega^{(2)}(\bm{x},\bm{y}) + {B_1}\omega^{(1)}(\bm{x},\bm{y}) + {B_2}\omega^{(0)}(\bm{x},\bm{y}) = 0.
\end{equation}
\begin{equation}
	O( {{\epsilon ^{ - 1}}}): {B_0}\omega^{(3)}(\bm{x},\bm{y}) + {B_1}\omega^{(2)}(\bm{x},\bm{y}) + {B_2}\omega^{(1)}(\bm{x},\bm{y})+ {B_3}\omega^{(0)}(\bm{x},\bm{y}) = 0.
\end{equation}
\begin{equation}
	\begin{aligned}
		O( {{\epsilon ^0}}): &{B_0}\omega^{(4)}(\bm{x},\bm{y}) + {B_1}\omega^{(3)}(\bm{x},\bm{y}) + {B_2}\omega^{(2)}(\bm{x},\bm{y})\\
&+ {B_3}\omega^{(1)}(\bm{x},\bm{y})+ {B_4}\omega^{(0)}(\bm{x},\bm{y}) = q(\bm{x}).
	\end{aligned}
\end{equation}

According to the $O( {{\epsilon ^{ - 4}}} )$-order equation (6), we find that $\omega^{(0)}$ is independent of the microscopic variable $\bm{y}$, i.e.
\begin{equation}
	\omega^{(0)}(\bm{x},\bm{y}) = \omega^{(0)}(\bm{x}).
\end{equation}
Thus $\omega^{(0)}(\bm{x})$ is called as macroscopic homogenized solution.

From the definition of $\omega^{(0)}(\bm{x})$, we can reduce the $O( {{\epsilon ^{ - 3}}} )$-order equation (7) to the following equation.
\begin{equation}
	{B_0}\omega^{(1)}(\bm{x},\bm{y}) = 0.
\end{equation}
After that, based on equation (12) and the periodicity of function $\omega^{(1)}(\bm{x},\bm{y})$ in microscopic variable $\bm{y}$, it is clear that we have the following.
\begin{equation}
	\omega^{(1)}(\bm{x},\bm{y}) = 0,
\end{equation}
Combining the conclusions of (11) and (13), we can simplify the $O( {{\epsilon ^{ - 2}}} )$-order equation (8) as the following equation.
\begin{equation}
	\frac{\partial^2 }{{\partial {y_i}\partial {y_j}}}\bigg[ {{D_{ijkl}}(\bm{y})\frac{{{\partial ^2 \omega^{(2)} (\bm{x},\bm{y})}}}{{\partial {y_k}\partial {y_l}}}} \bigg] = -\frac{\partial^2 D_{ij\alpha_1\alpha_2}(\bm{y})}{{\partial {y_i}\partial {y_j}}}\frac{{{\partial ^2 \omega^{(0)} (\bm{x})}}}{{\partial {x_{\alpha_1}}\partial {x_{\alpha_2}}}}.
\end{equation}
In virtue of the governing equation (14) concerning $\omega^{(2)} (\bm{x},\bm{y})$, we can define the following particular separation-of-variable form for $\omega^{(2)} (\bm{x},\bm{y})$
\begin{equation}
	\omega^{(2)}(\bm{x},\bm{y}) = \mathcal{N}_2^{\alpha_1\alpha_2}(\bm{y}) \frac{{{\partial ^2 \omega^{(0)} (\bm{x})}}}{{\partial {x_{\alpha_1}}\partial {x_{\alpha_2}}}},
\end{equation}
where $\mathcal{N}_2^{\alpha_1\alpha_2}(\bm{y})$ denotes the second-order cell functions. Moreover, substituting (15) into (14), the following equations are derived for solving $\mathcal{N}_2^{\alpha_1\alpha_2}(\bm{y})$ by attaching clamped boundary condition.
\begin{equation}
	\left\{  \begin{aligned}
		&\frac{\partial^2 }{{\partial {y_i}\partial {y_j}}}\bigg[ {{D_{ijkl}}(\bm{y})\frac{{{\partial ^2 \mathcal{N}_2^{\alpha_1\alpha_2}(\bm{y})}}}{{\partial {y_k}\partial {y_l}}}} \bigg] = -\frac{{\partial ^2 {D_{ij{\alpha _1}{\alpha _2}}}(\bm{y})}}{{\partial {y_i}\partial {y_j}}},\;\;\bm{y} \in \mathbf{Y}, \\
		&\mathcal{N}_2^{\alpha_1\alpha_2}(\bm{y})= 0,\;\;\qquad\bm{y} \in \partial \mathbf{Y},\\
		&\frac{{\partial \mathcal{N}_2^{\alpha_1\alpha_2}(\bm{y})}}{{\partial {y_j}}}{\hat n_j}=0,\;\;\bm{y} \in \partial \mathbf{Y}.
	\end{aligned} \right.
\end{equation}

Furthermore, together with the facts (11), (13) and (15), the $O( {{\epsilon ^{ - 1}}} )$-order equation (9) can be naturally reduced to the following equation.
\begin{equation}
\begin{aligned}
&\frac{\partial^2 }{{\partial {y_i}\partial {y_j}}}\bigg[ {{D_{ijkl}}(\bm{y})\frac{{{\partial ^2 \omega^{(3)} (\bm{x},\bm{y})}}}{{\partial {y_k}\partial {y_l}}}} \bigg] = -2\frac{\partial }{{\partial {y_j}}}\bigg[ {{D_{\alpha_1jkl}}(\bm{y})\frac{{{\partial ^2 \mathcal{N}_2^{\alpha_2\alpha_3}(\bm{y})}}}{{\partial {y_k}\partial {y_l}}}} \bigg]{\frac{{{\partial ^3 \omega^{(0)}(\bm{x})}}}{{\partial {x_{\alpha_1}}\partial {x_{\alpha_2}}\partial {x_{\alpha_3}}}}}\\
&-2\frac{\partial^2 }{{\partial {y_i}\partial {y_j}}}\bigg[ {{D_{ij\alpha_1l}}(\bm{y})\frac{{{\partial \mathcal{N}_2^{\alpha_2\alpha_3}(\bm{y})}}}{{\partial {y_l}}}} \bigg]{\frac{{{\partial ^3 \omega^{(0)}(\bm{x})}}}{{\partial {x_{\alpha_1}}\partial {x_{\alpha_2}}\partial {x_{\alpha_3}}}}}-2\frac{\partial {D_{i\alpha_1\alpha_2\alpha_3}}(\bm{y})}{{\partial {y_i}}}{\frac{{{\partial ^3 \omega^{(0)}(\bm{x})}}}{{\partial {x_{\alpha_1}}\partial {x_{\alpha_2}}\partial {x_{\alpha_3}}}}}.
\end{aligned}
\end{equation}
According to (17), the same technique to define (15) suggests $\omega^{(3)}(\bm{x},\bm{y})$ can be expressed in a separation-of-variable form as below.
\begin{equation}
	\omega^{(3)}(\bm{x},\bm{y}) = \mathcal{N}_3^{\alpha_1\alpha_2\alpha_3}(\bm{y}) \frac{{{\partial ^3 \omega^{(0)} (\bm{x})}}}{{\partial {x_{\alpha_1}}\partial {x_{\alpha_2}}\partial {x_{\alpha_3}}}},
\end{equation}
where $\mathcal{N}_3^{\alpha_1\alpha_2\alpha_3}(\bm{y})$ denotes the third-order cell functions. Again, substituting (18) into (17), we can obtain the governing equations for solving $\mathcal{N}_3^{\alpha_1\alpha_2\alpha_3}(\bm{y})$ by applying clamped boundary condition.
\begin{equation}
	\left\{  \begin{aligned}
		&\frac{\partial^2 }{{\partial {y_i}\partial {y_j}}}\bigg[ {{D_{ijkl}}(\bm{y})\frac{{{\partial ^2 \mathcal{N}_3^{\alpha_1\alpha_2\alpha_3}(\bm{y})}}}{{\partial {y_k}\partial {y_l}}}} \bigg] = -2\frac{\partial }{{\partial {y_j}}}\bigg[ {{D_{\alpha_1jkl}}(\bm{y})\frac{{{\partial ^2 \mathcal{N}_2^{\alpha_2\alpha_3}(\bm{y})}}}{{\partial {y_k}\partial {y_l}}}} \bigg]\\
&-2\frac{\partial^2 }{{\partial {y_i}\partial {y_j}}}\bigg[ {{D_{ij\alpha_1l}}(\bm{y})\frac{{{\partial \mathcal{N}_2^{\alpha_2\alpha_3}(\bm{y})}}}{{\partial {y_l}}}} \bigg]-2\frac{\partial {D_{i\alpha_1\alpha_2\alpha_3}}(\bm{y})}{{\partial {y_i}}},\;\;\bm{y} \in \mathbf{Y}, \\
		&\mathcal{N}_3^{\alpha_1\alpha_2\alpha_3}(\bm{y})= 0,\;\;\qquad\bm{y} \in \partial \mathbf{Y},\\
		&\frac{{\partial \mathcal{N}_3^{\alpha_1\alpha_2\alpha_3}(\bm{y})}}{{\partial {y_j}}}{\hat n_j}=0,\;\;\bm{y} \in \partial \mathbf{Y}.
	\end{aligned} \right.
\end{equation}

Subsequently, performing an integral on both sides of the $O( {{\epsilon ^{0}}} )$-order equation (10) over PUC $\mathbf{Y}$, exploiting the Gauss theorem on equation (10) and replacing herein $\omega^{(2)}$ by its formula (15) simultaneously, these procedures lead to derive the macroscopic homogenized equation for composite Kirchhoff plates as follows.
\begin{equation}
	\left\{  \begin{aligned}
		& - \frac{\partial^2 }{\partial {x_i}\partial {x_j}}\bigg[{\widehat D}_{ijkl}\frac{{{\partial ^2{\omega^{(0)}}(\bm{x})}}}{\partial {x_k}\partial {x_l}}\bigg] =q(\bm{x}) ,\;\; \text{in}\;\; \Omega ,\\
		&{\omega^{(0)}}(\bm{x}) = g_1(\bm{x}),\;\; \text{on}\;\; \partial \Omega,\\
		&\frac{{\partial \omega^{(0)}(\bm{x})}}{{\partial {x_j}}}{\hat n_j} = g_2(\bm{x}), \;\; \text{on}\;\; \partial \Omega,
	\end{aligned} \right.
\end{equation}
where the homogenized bending stiffness parameters ${\widehat D}_{ijkl}$ can be evaluated by the following formula.
\begin{equation}
	{\widehat D_{ijkl}} = \frac{1}{{\left|\mathbf{Y}\right|}}\int_\mathbf{Y} {\bigg[{D_{ijkl}}(\bm{y})+ {{D_{ij\alpha_1\alpha_2}}(\bm{y})\frac{{{\partial ^2 \mathcal{N}_2^{kl}(\bm{y})}}}{{\partial {y_{\alpha_1}}\partial {y_{\alpha_2}}}}}\bigg]} d\mathbf{Y}.
\end{equation}

Now, we initiate to establish the crucial fourth-order corrector $\omega^{(4)} (\bm{x},\bm{y})$. By inserting the particular expressions (11), (13), (15) and (18) into equation (10), and replacing $q(\bm{x})$ in equation (10) with its equivalent definition in macroscopic homogenized equation (20), we can reduce the $O( {{\epsilon ^{0}}} )$-order equation (10) as the following equation.
\begin{equation}
\begin{aligned}
&\frac{\partial^2 }{{\partial {y_i}\partial {y_j}}}\bigg[ {{D_{ijkl}}(\bm{y})\frac{{{\partial ^2 \omega^{(4)} (\bm{x},\bm{y})}}}{{\partial {y_k}\partial {y_l}}}} \bigg] = \widehat D_{\alpha_1\alpha_2\alpha_3\alpha_4}{\frac{{{\partial ^4 \omega^{(0)}(\bm{x})}}}{{\partial {x_{\alpha_1}}\partial {x_{\alpha_2}}\partial {x_{\alpha_3}}\partial {x_{\alpha_4}}}}}\\
&-D_{\alpha_1\alpha_2\alpha_3\alpha_4}(\bm{y}){\frac{{{\partial ^4 \omega^{(0)}(\bm{x})}}}{{\partial {x_{\alpha_1}}\partial {x_{\alpha_2}}\partial {x_{\alpha_3}}\partial {x_{\alpha_4}}}}}-2\frac{\partial }{{\partial {y_j}}}\bigg[ {{D_{\alpha_1jkl}}(\bm{y})\frac{{{\partial ^2 \mathcal{N}_3^{\alpha_2\alpha_3\alpha_4}(\bm{y})}}}{{\partial {y_k}\partial {y_l}}}} \bigg]{\frac{{{\partial ^4 \omega^{(0)}(\bm{x})}}}{{\partial {x_{\alpha_1}}\partial {x_{\alpha_2}}\partial {x_{\alpha_3}}\partial {x_{\alpha_4}}}}}\\
&-2\frac{\partial^2 }{{\partial {y_i}\partial {y_j}}}\bigg[ {{D_{ij\alpha_1l}}(\bm{y})\frac{{{\partial \mathcal{N}_3^{\alpha_2\alpha_3\alpha_4}(\bm{y})}}}{{\partial {y_l}}}} \bigg]{\frac{{{\partial ^4 \omega^{(0)}(\bm{x})}}}{{\partial {x_{\alpha_1}}\partial {x_{\alpha_2}}\partial {x_{\alpha_3}}\partial {x_{\alpha_4}}}}}\\
&-{{D_{\alpha_1\alpha_2kl}}(\bm{y})\frac{{{\partial^2 \mathcal{N}_2^{\alpha_3\alpha_4}(\bm{y})}}}{{\partial {y_k}\partial {y_l}}}}{\frac{{{\partial ^4 \omega^{(0)}(\bm{x})}}}{{\partial {x_{\alpha_1}}\partial {x_{\alpha_2}}\partial {x_{\alpha_3}}\partial {x_{\alpha_4}}}}}-\frac{\partial^2 }{{\partial {y_i}\partial {y_j}}}\bigg[{D_{ij\alpha_1\alpha_2}}(\bm{y})\mathcal{N}_2^{\alpha_3\alpha_4}(\bm{y})\bigg]{\frac{{{\partial ^4 \omega^{(0)}(\bm{x})}}}{{\partial {x_{\alpha_1}}\partial {x_{\alpha_2}}\partial {x_{\alpha_3}}\partial {x_{\alpha_4}}}}}\\
&-4\frac{\partial}{{\partial {y_j}}}\bigg[ {{D_{\alpha_1j\alpha_2l}}(\bm{y})\frac{{{\partial \mathcal{N}_2^{\alpha_3\alpha_4}(\bm{y})}}}{{\partial {y_l}}}} \bigg]{\frac{{{\partial ^4 \omega^{(0)}(\bm{x})}}}{{\partial {x_{\alpha_1}}\partial {x_{\alpha_2}}\partial {x_{\alpha_3}}\partial {x_{\alpha_4}}}}}.
\end{aligned}
\end{equation}
As observed in equation (22), we can establish the separation-of-variable expression for $omega^{(4)}(\bm{x},\bm{y})$ as follows.
\begin{equation}
	\omega^{(4)}(\bm{x},\bm{y}) = \mathcal{N}_4^{\alpha_1\alpha_2\alpha_3\alpha_4}(\bm{y}) \frac{{{\partial ^4 \omega^{(0)} (\bm{x})}}}{{\partial {x_{\alpha_1}}\partial {x_{\alpha_2}}\partial {x_{\alpha_3}}\partial {x_{\alpha_4}}}},
\end{equation}
where $\mathcal{N}_4^{\alpha_1\alpha_2\alpha_3\alpha_4}(\bm{y})$ denotes the fourth-order cell functions. Moreover, substituting (23) into (22), one can deduce the following governing equation for computing $\mathcal{N}_4^{\alpha_1\alpha_2\alpha_3\alpha_4}(\bm{y})$ by imposing clamped boundary condition.
\begin{equation}
	\left\{  \begin{aligned}
		&\frac{\partial^2 }{{\partial {y_i}\partial {y_j}}}\bigg[ {{D_{ijkl}}(\bm{y})\frac{{{\partial ^2 \mathcal{N}_4^{\alpha_1\alpha_2\alpha_3\alpha_4}(\bm{y})}}}{{\partial {y_k}\partial {y_l}}}} \bigg] = \widehat D_{\alpha_1\alpha_2\alpha_3\alpha_4}-D_{\alpha_1\alpha_2\alpha_3\alpha_4}(\bm{y})\\
&-2\frac{\partial }{{\partial {y_j}}}\bigg[ {{D_{\alpha_1jkl}}(\bm{y})\frac{{{\partial ^2 \mathcal{N}_3^{\alpha_2\alpha_3\alpha_4}(\bm{y})}}}{{\partial {y_k}\partial {y_l}}}} \bigg]-2\frac{\partial^2 }{{\partial {y_i}\partial {y_j}}}\bigg[ {{D_{ij\alpha_1l}}(\bm{y})\frac{{{\partial \mathcal{N}_3^{\alpha_2\alpha_3\alpha_4}(\bm{y})}}}{{\partial {y_l}}}} \bigg]\\
&-{{D_{\alpha_1\alpha_2kl}}(\bm{y})\frac{{{\partial^2 \mathcal{N}_2^{\alpha_3\alpha_4}(\bm{y})}}}{{\partial {y_k}\partial {y_l}}}}-\frac{\partial^2 }{{\partial {y_i}\partial {y_j}}}\bigg[{D_{ij\alpha_1\alpha_2}}(\bm{y})\mathcal{N}_2^{\alpha_3\alpha_4}(\bm{y})\bigg]\\
&-4\frac{\partial}{{\partial {y_j}}}\bigg[ {{D_{\alpha_1j\alpha_2l}}(\bm{y})\frac{{{\partial \mathcal{N}_2^{\alpha_3\alpha_4}(\bm{y})}}}{{\partial {y_l}}}} \bigg],\;\;\bm{y} \in \mathbf{Y}, \\
		&\mathcal{N}_4^{\alpha_1\alpha_2\alpha_3\alpha_4}(\bm{y})= 0,\;\;\qquad\bm{y} \in \partial \mathbf{Y},\\
		&\frac{{\partial \mathcal{N}_4^{\alpha_1\alpha_2\alpha_3\alpha_4}(\bm{y})}}{{\partial {y_j}}}{\hat n_j}=0,\;\;\bm{y} \in \partial \mathbf{Y}.
	\end{aligned} \right.
\end{equation}
\begin{rmk}
On the basis of Lax-Milgram theorem and the assumption (A$_2$), the existence and uniqueness of solutions for auxiliary cell problems (16), (19) and (24) can be proved in \cite{R33,R34}.
\end{rmk}
\begin{rmk}
On the basis of Lax-Milgram theorem and the assumptions (A$_2$) and (A$_3$), the existence and uniqueness of the macroscopic homogenized solution to the fourth-order homogenized PDE (20) can be proved.
\end{rmk}

Summing up, the above-mentioned multi-scale analysis can be concluded with the following theorem.
\begin{thm}
	The fourth-order multi-scale asymptotic solution of the plate bending problems of composite Kirchhoff plates is given as follows
	\begin{equation}
		\begin{aligned}
			&\omega^{(4\varepsilon)} (\bm{x}) = \omega^{(0)}(\bm{x}) + {\epsilon ^2}\mathcal{N}_2^{\alpha_1\alpha_2}(\bm{y}) \frac{{{\partial ^2 \omega^{(0)} (\bm{x})}}}{{\partial {x_{\alpha_1}}\partial {x_{\alpha_2}}}}+ {\epsilon ^3}\mathcal{N}_3^{\alpha_1\alpha_2\alpha_3}(\bm{y}) \frac{{{\partial ^3 \omega^{(0)} (\bm{x})}}}{{\partial {x_{\alpha_1}}\partial {x_{\alpha_2}}\partial {x_{\alpha_3}}}}\\
			&\;\qquad \quad +{\epsilon ^4}\mathcal{N}_4^{\alpha_1\alpha_2\alpha_3\alpha_4}(\bm{y}) \frac{{{\partial ^4 \omega^{(0)} (\bm{x})}}}{{\partial {x_{\alpha_1}}\partial {x_{\alpha_2}}\partial {x_{\alpha_3}}\partial {x_{\alpha_4}}}}.
		\end{aligned}
	\end{equation}
\end{thm}

According to the Kirchhoff plate theory, the displacement field of the composite plate for bending problems can be determined as follows.
\begin{equation}
		\begin{aligned}
u_1^\epsilon\approx u_1^{(4\epsilon)}=-x_3\frac{\partial \omega^{(4\epsilon)} (\bm{x})}{\partial x_1}
		\end{aligned}
	\end{equation}
\begin{equation}
		\begin{aligned}
u_2^\epsilon\approx u_2^{(4\epsilon)}=-x_3\frac{\partial \omega^{(4\epsilon)} (\bm{x})}{\partial x_2}
		\end{aligned}
	\end{equation}	
\begin{equation}
		\begin{aligned}
u_3^\epsilon\approx \omega^{(4\epsilon)} (\bm{x})
		\end{aligned}
	\end{equation}		
\section{Theoretical error analysis of multi-scale asymptotic solutions}
This section presents the theoretical error analysis of multi-scale asymptotic solutions for bending problems of composite Kirchhoff plates in the local and global senses in detail.
\subsection{Local error analysis in the point-wise sense}
First define the second-order multi-scale solution and third-order multi-scale solution for composite Kirchhoff plates as below.
\begin{equation}
	\omega^{(2\epsilon )} = \omega^{(0)} + \epsilon \omega^{(1)} + {\epsilon ^2}\omega^{(2)} = \omega^{(0)} + {\epsilon ^2}\omega^{(2)}.
\end{equation}
\begin{equation}
	\omega^{(3\epsilon )} = \omega^{(0)} + \epsilon \omega^{(1)} + {\epsilon ^2}\omega^{(2)} + {\epsilon ^3}\omega^{(3)}= \omega^{(0)} + {\epsilon ^2}\omega^{(2)}+ {\epsilon ^3}\omega^{(3)}.
\end{equation}
The fourth-order multi-scale solution $\omega^{(4\epsilon )}$ is defined as formula (25).

Furthermore, define the following residual functions for different kinds of multi-scale solutions.
\begin{equation}
	Z_{\Delta}^{(2\epsilon )}(\bm{x}) = \omega^\epsilon(\bm{x})  - \omega^{(2\epsilon )}(\bm{x}),\;\;\bm{x} \in \Omega.
\end{equation}
\begin{equation}
	Z_{\Delta}^{(3\epsilon )}(\bm{x}) = \omega^\epsilon(\bm{x})  - \omega^{(3\epsilon )}(\bm{x}),\;\;\bm{x} \in \Omega.
\end{equation}
\begin{equation}
	Z_{\Delta}^{(4\epsilon )}(\bm{x}) = \omega^\epsilon(\bm{x})  - \omega^{(4\epsilon )}(\bm{x}),\;\;\bm{x} \in \Omega.
\end{equation}

After that, substituting the above residual functions into the multi-scale fourth-order PDE (1) respectively, the residual equations for the second-order multi-scale solution, third-order multi-scale solution and fourth-order multi-scale solution are derived as follows.
\begin{equation}
		\left\{  \begin{aligned}
			&   \frac{\partial^2 }{{\partial {x_i}\partial {x_j}}}\bigg[ {D_{ijkl}^\epsilon (\bm{x})\frac{{{\partial ^2 Z_{\Delta}^{(2\epsilon )}(\bm{x})}}}{{\partial {x_k}\partial {x_l}}}} \bigg]\\
&\quad\quad\quad\quad\quad\quad\quad = -\epsilon^{-1}({B_1}\omega^{(2)}+ {B_2}\omega^{(1)}+ {B_3}\omega^{(0)})\\
&\quad\quad\quad\quad\quad\quad\quad-\epsilon^0({B_2}\omega^{(2)}+ {B_3}\omega^{(1)}+ {B_4}\omega^{(0)}-q)\\
&\quad\quad\quad\quad\quad\quad\quad-\epsilon({B_3}\omega^{(2)}+ {B_4}\omega^{(1)})-\epsilon^2B_4\omega^{(2)}:=\epsilon^{-1}E_{2}(\bm{x},\bm{y}), \;\; \text{in}\;\; \Omega,  \\
		& Z_{\Delta}^{(2\epsilon )}(\bm{x})= - {\epsilon ^2}\omega^{(2)}(\bm{x},\bm{y}):={\epsilon ^2}\chi_{2}(\bm{x}), \;\; \text{on}\;\; \partial \Omega ,\\
		&\frac{{\partial Z_{\Delta}^{(2\epsilon )}(\bm{x})}}{{\partial {x_j}}}{\hat n_j}= \frac{{\partial \omega^{\epsilon}(\bm{x})}}{{\partial {x_j}}}{\hat n_j}-\frac{{\partial \omega^{(2\epsilon )}(\bm{x})}}{{\partial {x_j}}}{\hat n_j} := \epsilon I_{2j}(\bm{x}){\hat n_j}, \;\; \text{on}\;\; \partial \Omega.
		\end{aligned} \right.
\end{equation}
\begin{equation}
		\left\{  \begin{aligned}
			&   \frac{\partial^2 }{{\partial {x_i}\partial {x_j}}}\bigg[ {D_{ijkl}^\epsilon (\bm{x})\frac{{{\partial ^2 Z_{\Delta}^{(3\epsilon )}(\bm{x})}}}{{\partial {x_k}\partial {x_l}}}} \bigg]\\
&\quad\quad\quad\quad\quad= -\epsilon^0({B_1}\omega^{(3)}+{B_2}\omega^{(2)}+ {B_3}\omega^{(1)}+ {B_4}\omega^{(0)}-q)\\
&\quad\quad\quad\quad\quad-\epsilon({B_2}\omega^{(3)}+{B_3}\omega^{(2)}+ {B_4}\omega^{(1)})-\epsilon^2(B_3\omega^{(3)}+B_4\omega^{(2)})\\
&\quad\quad\quad\quad\quad-\epsilon^3B_4\omega^{(3)}:=\epsilon^{0}E_{3}(\bm{x},\bm{y}), \;\; \text{in}\;\; \Omega,  \\
		& Z_{\Delta}^{(3\epsilon )}(\bm{x})= -{\epsilon ^2}\omega^{(2)}(\bm{x},\bm{y})- {\epsilon ^3}\omega^{(3)}(\bm{x},\bm{y}):={\epsilon ^2}\chi_{3}(\bm{x}), \;\; \text{on}\;\; \partial \Omega ,\\
		&\frac{{\partial Z_{\Delta}^{(3\epsilon )}(\bm{x})}}{{\partial {x_j}}}{\hat n_j}= \frac{{\partial \omega^{\epsilon}(\bm{x})}}{{\partial {x_j}}}{\hat n_j}-\frac{{\partial \omega^{(3\epsilon )}(\bm{x})}}{{\partial {x_j}}}{\hat n_j} :=\epsilon I_{3j}(\bm{x}){\hat n_j}, \;\; \text{on}\;\; \partial \Omega.
		\end{aligned} \right.
\end{equation}
\begin{equation}
		\left\{  \begin{aligned}
			&   \frac{\partial^2 }{{\partial {x_i}\partial {x_j}}}\bigg[ {D_{ijkl}^\epsilon (\bm{x})\frac{{{\partial ^2 Z_{\Delta}^{(4\epsilon )}(\bm{x})}}}{{\partial {x_k}\partial {x_l}}}} \bigg]\\
&\quad\quad\quad\quad\quad= -\epsilon({B_1}\omega^{(4)}+{B_2}\omega^{(3)}+{B_3}\omega^{(2)}+ {B_4}\omega^{(1)})\\
&\quad\quad\quad\quad\quad-\epsilon^2(B_2\omega^{(4)}+B_3\omega^{(3)}+B_4\omega^{(2)})\\
&\quad\quad\quad\quad\quad-\epsilon^3(B_3\omega^{(4)}+B_4\omega^{(3)})-\epsilon^4B_4\omega^{(4)}:=\epsilon E_{4}(\bm{x},\bm{y}), \;\; \text{in}\;\; \Omega,  \\
		& Z_{\Delta}^{(4\epsilon )}(\bm{x})=  -{\epsilon ^2}\omega^{(2)}(\bm{x},\bm{y})- {\epsilon ^3}\omega^{(3)}(\bm{x},\bm{y}) - {\epsilon ^4}\omega^{(4)}(\bm{x},\bm{y}):={\epsilon ^2}\chi_{4}(\bm{x}), \;\; \text{on}\;\; \partial \Omega ,\\
		&\frac{{\partial Z_{\Delta}^{(4\epsilon )}(\bm{x})}}{{\partial {x_j}}}{\hat n_j} = \frac{{\partial \omega^{\epsilon}(\bm{x})}}{{\partial {x_j}}}{\hat n_j}-\frac{{\partial \omega^{(4\epsilon )}(\bm{x})}}{{\partial {x_j}}}{\hat n_j}: =\epsilon I_{4j}(\bm{x}){\hat n_j}, \;\; \text{on}\;\; \partial \Omega.
		\end{aligned} \right.
\end{equation}

Analyzing the residual equations (34) and (35), we can clearly find that second-order multi-scale solution and third-order multi-scale solution fail to preserve local mechanical balance due to the error terms $O(\epsilon^{-1})$ and $O(\epsilon^{0})$ in their right sides with the quantitative relationship $\epsilon^{-1}\ll \epsilon^{0}\ll \epsilon$. In contrast, the error term of fourth-order multi-scale solution in the residual equation (36) shall preserve the local mechanical balance and provide the high-accuracy performance when $\epsilon$ is a small constant. Additionally, when the microstructural parameter $\epsilon  \to 0$, the right-hand error term of the residual equation (36) also tends to zero, indicating that the the proposed fourth-order multi-scale solution $\omega^{(4\epsilon )}$ converges to the exact solution $\omega^{\epsilon}$ in a local point-wise sense. This capability is the primary impetus for current study to develop the FOMS solution manifesting high-fidelity computation performance for composite Kirchhoff plates.

\subsection{Global error estimate in the energy sense}
To the best of our knowledge, the following error estimate has been established for fourth-order Kirchhoff plate model in \cite{R33}.
\begin{equation}
{\left\| \omega^{\epsilon}(\bm{x})\right\|_{{{{H^2}(\Omega )}}}}\!\leq\!C_{\Omega}{\left\|q(\bm{x}) \right\|_{{H^{-2}}(\Omega )}} + C_{\Omega}{\left\|g_1(\bm{x}) \right\|_{{H^{3/2}}(\Omega )}} + C_{\Omega}{\left\|g_2(\bm{x}) \right\|_{{H^{1/2}}(\Omega )}}.
\end{equation}
However, this error estimate can not be employed to obtain the explicit error order for fourth-order multi-scale solution $\omega^{(4\epsilon )}$ with residual equation (36). The reason is that the appearance of $\epsilon ^{-2}$-order terms lead to an order reduction phenomenon when using trace theorem to ${\left\|{\epsilon ^2}\chi_{4}(\bm{x})\right\|_{{H^{3/2}}(\Omega )}}$, which makes it impossible to obtain an explicit error estimate.

In order to an explicit error estimate for fourth-order multi-scale solution $\omega^{(4\epsilon )}$, we need to present some assumptions as follows.
\begin{enumerate}
	\item[(B$_1$)]
	Suppose $\Omega  \subset \mathbb{R}^2$ is a bounded convex domain and a union of integral unit cells, i.e. $\bar{\Omega}=\cup_{\mathbf{z}\in I_{\epsilon}}\epsilon(\mathbf{z}+\bar{\mathbf{Y}})$, where the index set $I_{\epsilon}=\{\mathbf{z}=(z_1,z_2)\in Z^2,\epsilon(\mathbf{z}+\bar{\mathbf{Y}})\subset \bar{\Omega}\}$. Besides, let $E_\mathbf{z}=\epsilon(\mathbf{z}+\mathbf{Y})$ be the translational unit cell and $\partial E_\mathbf{z}$ be its boundary with $E_\mathbf{z}=\ell_1\bigcup\ell_2\bigcup\ell_3\bigcup\ell_4$, as shown in Fig.~1.
	\item[(B$_2$)]
	The stiffness function ${D_{ijkl}}(\bm{y})$ is piecewise constants on the microscopic PUC $\mathbf{Y}$.
	\item[(B$_3$)]
	Define ${\bigtriangleup_1}$ and ${\bigtriangleup_2}$ to be middle hyperplanes on the microscopic PUC $\mathbf{Y}$, and assume that the microscopic PUC $\mathbf{Y}$ is geometrically symmetric with respect to ${\bigtriangleup_1}$ and ${\bigtriangleup_2}$, as exhibited in Fig.~1.
\end{enumerate}
\begin{figure}[!htb]
	\centering
	\begin{minipage}[c]{0.5\textwidth}
		\centering
		\includegraphics[width=0.8\linewidth,totalheight=2.1in]{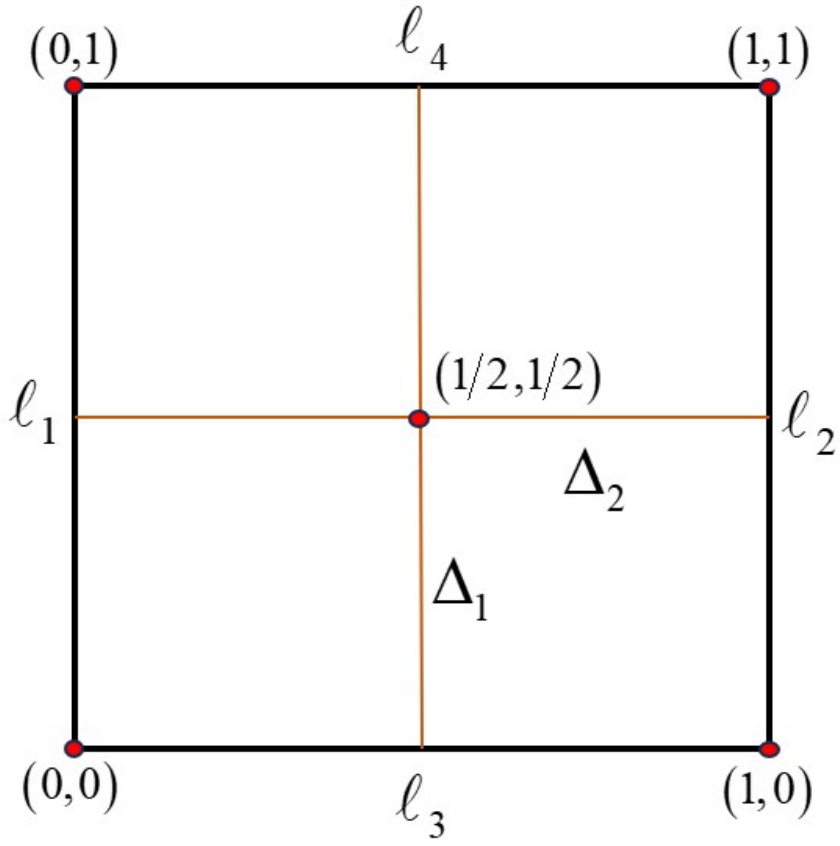} \\
	\end{minipage}
	\caption{The symmetry and boundaries of microscopic PUC $\mathbf{Y}$.}\label{f1}
\end{figure}

\begin{lem}
Define the differential operator $\displaystyle{\sigma _{\mathbf{Y}}}({\Phi} ) = {\hat n_i}\frac{\partial}{\partial y_j}\big[{D_{ijkl}}(\bm{y})\frac{{\partial^2 {\Phi}}}{{\partial {y_k}\partial {y_l}}}\big]$, and under the assumptions (A$_1$)-(A$_2$) and (B$_2$)-(B$_3$), it can be proven that ${\sigma _{\mathbf{Y}}}(\mathcal{N}_2^{\alpha_1\alpha_2})$ is continuous on the boundary $E_\mathbf{z}$.
\end{lem}
$\mathbf{Proof:}$
Inspired by the references \cite{R15,R16,R35} and the definition (16) of second-order cell function $\mathcal{N}_2^{\alpha_1\alpha_2}$, we define the following two functions.
\begin{equation}
\Lambda_1^{\alpha_1\alpha_2}(\bm{y})=\frac{\partial }{{\partial {y_j}}}\bigg[ {{D_{1jkl}}(\bm{y})\frac{{{\partial ^2 \mathcal{N}_2^{\alpha_1\alpha_2}(\bm{y})}}}{{\partial {y_k}\partial {y_l}}}} \bigg]+\frac{{\partial {D_{1j{\alpha _1}{\alpha _2}}}(\bm{y})}}{{\partial {y_j}}},
\end{equation}
\begin{equation}
\Lambda_2^{\alpha_1\alpha_2}(\bm{y})=\frac{\partial }{{\partial {y_j}}}\bigg[ {{D_{2jkl}}(\bm{y})\frac{{{\partial ^2 \mathcal{N}_2^{\alpha_1\alpha_2}(\bm{y})}}}{{\partial {y_k}\partial {y_l}}}} \bigg]+\frac{{\partial {D_{2j{\alpha _1}{\alpha _2}}}(\bm{y})}}{{\partial {y_j}}}.
\end{equation}
Then, we can derive the following equality according to the definition (16).
\begin{equation}
\frac{\partial \Lambda_1^{\alpha_1\alpha_2}(\bm{y})}{\partial y_1}+\frac{\partial \Lambda_2^{\alpha_1\alpha_2}(\bm{y})}{\partial y_2}=0.
\end{equation}
Furthermore, we set the test function $\displaystyle v(\bm{y})=(e^{i2\pi m_1y_1}-1)(e^{i2\pi m_2y_2}-1)$, $m_1\neq 0$, $m_2\neq 0$. Next, employing the assumptions (B$_2$)-(B$_3$), and the same approach in references \cite{R15,R16,R35}, and defining $v_1(\bm{y})=e^{i2\pi m_2y_2}-1$ and $v_2(\bm{y})=e^{i2\pi m_1y_1}-1$, it follows that
\begin{equation}
\int_{\mathbf{Y}}\Big(\frac{\partial \Lambda_1^{\alpha_1\alpha_2}(\bm{y})}{\partial y_1}+\frac{\partial \Lambda_2^{\alpha_1\alpha_2}(\bm{y})}{\partial y_2}\Big)v_1(\bm{y})d{\mathbf{Y}}=\int_{\ell_1}[\Lambda_1^{\alpha_1\alpha_2}(\bm{y})]v_1(\bm{y})ds=0,
\end{equation}
\begin{equation}
\int_{\mathbf{Y}}\Big(\frac{\partial \Lambda_1^{\alpha_1\alpha_2}(\bm{y})}{\partial y_1}+\frac{\partial \Lambda_2^{\alpha_1\alpha_2}(\bm{y})}{\partial y_2}\Big)v_2(\bm{y})d{\mathbf{Y}}=\int_{\ell_3}[\Lambda_2^{\alpha_1\alpha_2}(\bm{y})]v_2(\bm{y})ds=0,
\end{equation}
where $[\Lambda_1^{\alpha_1\alpha_2}(\bm{y})]$ and $[\Lambda_2^{\alpha_1\alpha_2}(\bm{y})]$ denotes the jump of the functions $\Lambda_1^{\alpha_1\alpha_2}(\bm{y})$ and $\Lambda_2^{\alpha_1\alpha_2}(\bm{y})$, respectively. Moreover, by virtue of the completeness of the function families $\displaystyle\{e^{i2\pi m_1y_1}\}_{m_1=-\infty}^{+\infty}$ and $\displaystyle\{e^{i2\pi m_2y_2}\}_{m_2=-\infty}^{+\infty}$, it can be derived that $[\Lambda_1^{\alpha_1\alpha_2}(\bm{y})]|_{\ell_1}=0$ and $[\Lambda_2^{\alpha_1\alpha_2}(\bm{y})]|_{\ell_3}=0$.

Moreover, the following two equalities can be easily obtained
\begin{equation}
\int_{\mathbf{Y}}\Big(\frac{\partial \Lambda_1^{\alpha_1\alpha_2}(\bm{y})}{\partial y_1}+\frac{\partial \Lambda_2^{\alpha_1\alpha_2}(\bm{y})}{\partial y_2}\Big)d{\mathbf{Y}}=\int_{\ell_1}[\Lambda_1^{\alpha_1\alpha_2}(\bm{y})]ds+\int_{\ell_3}[\Lambda_2^{\alpha_1\alpha_2}(\bm{y})]ds=0.
\end{equation}
\begin{equation}
\int_{\mathbf{Y}}\Big(\frac{\partial \Lambda_1^{\alpha_1\alpha_2}(\bm{y})}{\partial y_1}+\frac{\partial \Lambda_2^{\alpha_1\alpha_2}(\bm{y})}{\partial y_2}\Big)d{\mathbf{Y}}=\int_{\ell_1}[{\sigma _{\mathbf{Y}}}(\mathcal{N}_2^{\alpha_1\alpha_2})]ds+\int_{\ell_3}[{\sigma _{\mathbf{Y}}}(\mathcal{N}_2^{\alpha_1\alpha_2})]ds=0.
\end{equation}
Hence, we can deduce that ${\sigma _{\mathbf{Y}}}(\mathcal{N}_2^{\alpha_1\alpha_2})]ds+\int_{\ell_3}[{\sigma _{\mathbf{Y}}}(\mathcal{N}_2^{\alpha_1\alpha_2})$ is continuous on $E_\mathbf{z}$ on the basis of $[\Lambda_1^{\alpha_1\alpha_2}(\bm{y})]|_{\ell_1}=0$ and $[\Lambda_2^{\alpha_1\alpha_2}(\bm{y})]|_{\ell_3}=0$.
\begin{lem}
Under the assumptions (A$_1$)-(A$_2$) and (B$_2$)-(B$_3$), it can be proven that ${\sigma _{\mathbf{Y}}}(\mathcal{N}_3^{\alpha_1\alpha_2\alpha_3})$ and ${\sigma _{\mathbf{Y}}}(\mathcal{N}_4^{\alpha_1\alpha_2\alpha_3\alpha_4})$ are continuous on the boundary $E_\mathbf{z}$.
\end{lem}
$\mathbf{Proof:}$
Inspired by the same approach in references \cite{R15,R16,R35}, the continuous property of ${\sigma _{\mathbf{Y}}}(\mathcal{N}_3^{\alpha_1\alpha_2\alpha_3})$ and ${\sigma _{\mathbf{Y}}}(\mathcal{N}_4^{\alpha_1\alpha_2\alpha_3\alpha_4})$ can be proved.
\begin{lem}
Based on the lemmas 1 and 2, it can be proven that the first-order and second-order derivatives of all microscopic cell functions $\mathcal{N}_2^{\alpha_1\alpha_2}$, $\mathcal{N}_3^{\alpha_1\alpha_2\alpha_3}$ and $\mathcal{N}_4^{\alpha_1\alpha_2\alpha_3\alpha_4}$ are continuous on the boundary $E_\mathbf{z}$.
\end{lem}
$\mathbf{Proof:}$
From the lemmas 1 and 2, we have acknowledged that the third-order derivatives of all microscopic cell functions are continuous on the boundary $E_\mathbf{z}$. Hence, it is no doubt that first-order and second-order derivatives of all microscopic cell functions are continuous on the boundary $E_\mathbf{z}$.

Based on the aforementioned assumptions and lemmas, we now present the ultimate result of global error estimate for the FOMS solution $\omega^{(4\epsilon)}(\bm{x})$ of the fourth-order Kirchhoff plate equation (1) as the subsequent theorem.
\begin{thm}
Consider $\omega^\epsilon(\bm{x})$ as the weak solution of multi-scale fourth-order equation (1), while $\omega^{(0)}(\bm{x})$ represents the weak solution of corresponding homogenized equation (20). Let ${\omega^{(4\epsilon)}}(\bm{x})$ denotes the FOMS solution given by formula (25). Assuming the above hypotheses (A$_1$)-(A$_3$) and (B$_1$)-(B$_3$) hold, we derive the subsequent global error estimation if $\omega^{(0)}(\bm{x})\in H^8(\Omega)$.
\begin{equation}
{\begin{aligned}
\big\|{\omega^{\epsilon}}(\bm{x})-{\omega^{(4\epsilon)}}(\bm{x})\big\|_{H_0^2(\Omega)}\leq C_\Omega\epsilon.
\end{aligned}}
\end{equation}
\end{thm}
Here $C_\Omega$ is a positive constant dependent of $\Omega$, but irrespective of $\epsilon$.\\
$\mathbf{Proof:}$
The residual equation and associated boundary conditions in (36) are utilized to obtain the explicit error estimate (40). Firstly, multiplying both sides of the residual equation (36) by the test function $Z_{\Delta}^{(4\epsilon )}(\bm{x})$ and using integration on $\Omega$, the following equation is obtained
\begin{equation}
	\begin{aligned}
\int_\Omega \frac{\partial^2 }{{\partial {x_i}\partial {x_j}}}\bigg[ {D_{ijkl}^\epsilon (\bm{x})\frac{{{\partial ^2 Z_{\Delta}^{(4\epsilon )}(\bm{x})}}}{{\partial {x_k}\partial {x_l}}}} \bigg]Z_{\Delta}^{(4\epsilon )}(\bm{x})d\Omega\!=\!\int_\Omega \epsilon {E_{4}(\bm{x},\bm{y})} Z_{\Delta}^{(4\epsilon )}(\bm{x})d\Omega.
	\end{aligned}
\end{equation}
Immediately after applying Green's formula twice to the left term of equation (39), then we have
\begin{equation}
	\begin{aligned}
&\int_\Omega {D_{ijkl}^\epsilon (\bm{x})\frac{{{\partial ^2 Z_{\Delta}^{(4\epsilon )}(\bm{x})}}}{{\partial {x_k}\partial {x_l}}}}\frac{\partial^2 Z_{\Delta}^{(4\epsilon )}(\bm{x})}{{\partial {x_i}\partial {x_j}}}d\Omega\\
&+\int_{\partial\Omega} \frac{\partial }{{\partial {x_j}}}\bigg[ {D_{ijkl}^\epsilon (\bm{x})\frac{{{\partial ^2 Z_{\Delta}^{(4\epsilon )}(\bm{x})}}}{{\partial {x_k}\partial {x_l}}}} \bigg]Z_{\Delta}^{(4\epsilon )}(\bm{x}){\hat n_i}ds\\
&-\int_{\partial\Omega} {D_{ijkl}^\epsilon (\bm{x})\frac{{{\partial ^2 Z_{\Delta}^{(4\epsilon )}(\bm{x})}}}{{\partial {x_k}\partial {x_l}}}}\frac{\partial Z_{\Delta}^{(4\epsilon )}(\bm{x})}{{\partial {x_i}}}{\hat n_j}ds\\
&+\sum_{\mathbf{z}\in I_{\epsilon}}\int_{\partial E_\mathbf{z}} \frac{\partial }{{\partial {x_j}}}\bigg[ {D_{ijkl}^\epsilon (\bm{x})\frac{{{\partial ^2 Z_{\Delta}^{(4\epsilon )}(\bm{x})}}}{{\partial {x_k}\partial {x_l}}}} \bigg]Z_{\Delta}^{(4\epsilon )}(\bm{x}){\hat n_i}d{\Gamma _{\bm{y}}}\\
&-\sum_{\mathbf{z}\in I_{\epsilon}}\int_{\partial E_\mathbf{z}} {D_{ijkl}^\epsilon (\bm{x})\frac{{{\partial ^2 Z_{\Delta}^{(4\epsilon )}(\bm{x})}}}{{\partial {x_k}\partial {x_l}}}}\frac{\partial Z_{\Delta}^{(4\epsilon )}(\bm{x})}{{\partial {x_i}}}{\hat n_j}d{\Gamma _{\bm{y}}}\\
&=\int_\Omega \epsilon {E_{4}(\bm{x},\bm{y})} Z_{\Delta}^{(4\epsilon )}(\bm{x})d\Omega,
	\end{aligned}
\end{equation}
where the integral terms on $\partial E_\mathbf{z}$ arise from employing Green's formula on interface $\partial E_\mathbf{z}$ between two adjacent PUCs.

On recalling the assumption B$_1$ and the boundary conditions of microscopic cell functions, thus we can derive
\begin{equation}
Z_{\Delta}^{(4\epsilon )}(\bm{x})={\epsilon ^2}\chi_{4}(\bm{x})=0, \;\; \text{on}\;\; \partial \Omega,
\end{equation}
\begin{equation}
\frac{{\partial Z_{\Delta}^{(4\epsilon )}(\bm{x})}}{{\partial {x_j}}}{\hat n_j}=\epsilon I_{4j}(\bm{x}){\hat n_j}=0, \;\; \text{on}\;\; \partial \Omega.
\end{equation}
Considering the assumption B$_2$, and the lemmas 1-3, thus we shall obtain
\begin{small}
\begin{equation}
	\begin{aligned}
		&\sum_{\mathbf{z}\in I_{\epsilon}}\int_{\partial E_\mathbf{z}} \frac{\partial }{{\partial {x_j}}}\bigg[ {D_{ijkl}^\epsilon (\bm{x})\frac{{{\partial ^2 Z_{\Delta}^{(4\epsilon )}(\bm{x})}}}{{\partial {x_k}\partial {x_l}}}} \bigg]Z_{\Delta}^{(4\epsilon )}(\bm{x}){\hat n_i}d{\Gamma _{\bm{y}}}\\
        & = \sum_{\mathbf{z}\in I_{\epsilon}}\int_{\partial E_\mathbf{z}}\!\frac{\partial }{{\partial {x_j}}}\bigg[ {D_{ijkl}^\epsilon (\bm{x})\frac{{{\partial ^2 \omega^{\epsilon}(\bm{x})}}}{{\partial {x_k}\partial {x_l}}}} \bigg]Z_{\Delta}^{(4\epsilon )}(\bm{x}){\hat n_i}d{\Gamma _{\bm{y}}}\!-\!\sum_{\mathbf{z}\in I_{\epsilon}}\int_{\partial E_\mathbf{z}}\!\frac{\partial }{{\partial {x_j}}}\bigg[ {D_{ijkl}^\epsilon (\bm{x})\frac{{{\partial ^2 \omega^{(4\epsilon)}(\bm{x})}}}{{\partial {x_k}\partial {x_l}}}} \bigg]Z_{\Delta}^{(4\epsilon )}(\bm{x}){\hat n_i}d{\Gamma _{\bm{y}}}\\
		& =0-\epsilon^0\sum_{\mathbf{z}\in I_{\epsilon}}\int_{\partial E_\mathbf{z}} \frac{\partial }{{\partial {x_j}}}\bigg[ {D_{ijkl}(\bm{y})\frac{{{\partial ^2 }}}{{\partial {x_k}\partial {x_l}}}}\Big(\omega^{(0)} + {\epsilon ^2}\mathcal{N}_2^{\alpha_1\alpha_2}\frac{{{\partial ^2 \omega^{(0)} }}}{{\partial {x_{\alpha_1}}\partial {x_{\alpha_2}}}}\\
&+ {\epsilon ^3}\mathcal{N}_3^{\alpha_1\alpha_2\alpha_3}\frac{{{\partial ^3 \omega^{(0)}}}}{{\partial {x_{\alpha_1}}\partial {x_{\alpha_2}}\partial {x_{\alpha_3}}}}+{\epsilon ^4}\mathcal{N}_4^{\alpha_1\alpha_2\alpha_3\alpha_4} \frac{{{\partial ^4 \omega^{(0)}}}}{{\partial {x_{\alpha_1}}\partial {x_{\alpha_2}}\partial {x_{\alpha_3}}\partial {x_{\alpha_4}}}}\Big) \bigg]Z_{\Delta}^{(4\epsilon )}(\bm{x}){\hat n_i}d{\Gamma _{\bm{y}}}\\
        &-\epsilon^{-1}\sum_{\mathbf{z}\in I_{\epsilon}}\int_{\partial E_\mathbf{z}} \frac{\partial }{{\partial {x_j}}}\bigg[ {D_{ijkl}(\bm{y})\frac{{{\partial ^2 }}}{{\partial {x_k}\partial {y_l}}}}\Big(\omega^{(0)} + {\epsilon ^2}\mathcal{N}_2^{\alpha_1\alpha_2}\frac{{{\partial ^2 \omega^{(0)} }}}{{\partial {x_{\alpha_1}}\partial {x_{\alpha_2}}}}\\
&+ {\epsilon ^3}\mathcal{N}_3^{\alpha_1\alpha_2\alpha_3}\frac{{{\partial ^3 \omega^{(0)}}}}{{\partial {x_{\alpha_1}}\partial {x_{\alpha_2}}\partial {x_{\alpha_3}}}}+{\epsilon ^4}\mathcal{N}_4^{\alpha_1\alpha_2\alpha_3\alpha_4} \frac{{{\partial ^4 \omega^{(0)}}}}{{\partial {x_{\alpha_1}}\partial {x_{\alpha_2}}\partial {x_{\alpha_3}}\partial {x_{\alpha_4}}}}\Big) \bigg]Z_{\Delta}^{(4\epsilon )}(\bm{x}){\hat n_i}d{\Gamma _{\bm{y}}}\\
        &-\epsilon^{-1}\sum_{\mathbf{z}\in I_{\epsilon}}\int_{\partial E_\mathbf{z}} \frac{\partial }{{\partial {x_j}}}\bigg[ {D_{ijkl}(\bm{y})\frac{{{\partial ^2 }}}{{\partial {y_k}\partial {x_l}}}}\Big(\omega^{(0)} + {\epsilon ^2}\mathcal{N}_2^{\alpha_1\alpha_2}\frac{{{\partial ^2 \omega^{(0)} }}}{{\partial {x_{\alpha_1}}\partial {x_{\alpha_2}}}}\\
&+ {\epsilon ^3}\mathcal{N}_3^{\alpha_1\alpha_2\alpha_3}\frac{{{\partial ^3 \omega^{(0)}}}}{{\partial {x_{\alpha_1}}\partial {x_{\alpha_2}}\partial {x_{\alpha_3}}}}+{\epsilon ^4}\mathcal{N}_4^{\alpha_1\alpha_2\alpha_3\alpha_4} \frac{{{\partial ^4 \omega^{(0)}}}}{{\partial {x_{\alpha_1}}\partial {x_{\alpha_2}}\partial {x_{\alpha_3}}\partial {x_{\alpha_4}}}}\Big) \bigg]Z_{\Delta}^{(4\epsilon )}(\bm{x}){\hat n_i}d{\Gamma _{\bm{y}}}\\
        &-\epsilon^{-1}\sum_{\mathbf{z}\in I_{\epsilon}}\int_{\partial E_\mathbf{z}} \frac{\partial }{{\partial {y_j}}}\bigg[ {D_{ijkl}(\bm{y})\frac{{{\partial ^2 }}}{{\partial {x_k}\partial {x_l}}}}\Big(\omega^{(0)} + {\epsilon ^2}\mathcal{N}_2^{\alpha_1\alpha_2}\frac{{{\partial ^2 \omega^{(0)} }}}{{\partial {x_{\alpha_1}}\partial {x_{\alpha_2}}}}\\
&+ {\epsilon ^3}\mathcal{N}_3^{\alpha_1\alpha_2\alpha_3}\frac{{{\partial ^3 \omega^{(0)}}}}{{\partial {x_{\alpha_1}}\partial {x_{\alpha_2}}\partial {x_{\alpha_3}}}}+{\epsilon ^4}\mathcal{N}_4^{\alpha_1\alpha_2\alpha_3\alpha_4} \frac{{{\partial ^4 \omega^{(0)}}}}{{\partial {x_{\alpha_1}}\partial {x_{\alpha_2}}\partial {x_{\alpha_3}}\partial {x_{\alpha_4}}}}\Big) \bigg]Z_{\Delta}^{(4\epsilon )}(\bm{x}){\hat n_i}d{\Gamma _{\bm{y}}}\\
        &-\epsilon^{-2}\sum_{\mathbf{z}\in I_{\epsilon}}\int_{\partial E_\mathbf{z}} \frac{\partial }{{\partial {x_j}}}\bigg[ {D_{ijkl}(\bm{y})\frac{{{\partial ^2 }}}{{\partial {y_k}\partial {y_l}}}}\Big(\omega^{(0)} + {\epsilon ^2}\mathcal{N}_2^{\alpha_1\alpha_2}\frac{{{\partial ^2 \omega^{(0)} }}}{{\partial {x_{\alpha_1}}\partial {x_{\alpha_2}}}}\\
&+ {\epsilon ^3}\mathcal{N}_3^{\alpha_1\alpha_2\alpha_3}\frac{{{\partial ^3 \omega^{(0)}}}}{{\partial {x_{\alpha_1}}\partial {x_{\alpha_2}}\partial {x_{\alpha_3}}}}+{\epsilon ^4}\mathcal{N}_4^{\alpha_1\alpha_2\alpha_3\alpha_4} \frac{{{\partial ^4 \omega^{(0)}}}}{{\partial {x_{\alpha_1}}\partial {x_{\alpha_2}}\partial {x_{\alpha_3}}\partial {x_{\alpha_4}}}}\Big) \bigg]Z_{\Delta}^{(4\epsilon )}(\bm{x}){\hat n_i}d{\Gamma _{\bm{y}}}\\
        &-\epsilon^{-2}\sum_{\mathbf{z}\in I_{\epsilon}}\int_{\partial E_\mathbf{z}} \frac{\partial }{{\partial {y_j}}}\bigg[ {D_{ijkl}(\bm{y})\frac{{{\partial ^2 }}}{{\partial {x_k}\partial {y_l}}}}\Big(\omega^{(0)} + {\epsilon ^2}\mathcal{N}_2^{\alpha_1\alpha_2}\frac{{{\partial ^2 \omega^{(0)} }}}{{\partial {x_{\alpha_1}}\partial {x_{\alpha_2}}}}\\
&+ {\epsilon ^3}\mathcal{N}_3^{\alpha_1\alpha_2\alpha_3}\frac{{{\partial ^3 \omega^{(0)}}}}{{\partial {x_{\alpha_1}}\partial {x_{\alpha_2}}\partial {x_{\alpha_3}}}}+{\epsilon ^4}\mathcal{N}_4^{\alpha_1\alpha_2\alpha_3\alpha_4} \frac{{{\partial ^4 \omega^{(0)}}}}{{\partial {x_{\alpha_1}}\partial {x_{\alpha_2}}\partial {x_{\alpha_3}}\partial {x_{\alpha_4}}}}\Big) \bigg]Z_{\Delta}^{(4\epsilon )}(\bm{x}){\hat n_i}d{\Gamma _{\bm{y}}}\\
        &-\epsilon^{-2}\sum_{\mathbf{z}\in I_{\epsilon}}\int_{\partial E_\mathbf{z}} \frac{\partial }{{\partial {y_j}}}\bigg[ {D_{ijkl}(\bm{y})\frac{{{\partial ^2 }}}{{\partial {y_k}\partial {x_l}}}}\Big(\omega^{(0)} + {\epsilon ^2}\mathcal{N}_2^{\alpha_1\alpha_2}\frac{{{\partial ^2 \omega^{(0)} }}}{{\partial {x_{\alpha_1}}\partial {x_{\alpha_2}}}}\\
&+ {\epsilon ^3}\mathcal{N}_3^{\alpha_1\alpha_2\alpha_3}\frac{{{\partial ^3 \omega^{(0)}}}}{{\partial {x_{\alpha_1}}\partial {x_{\alpha_2}}\partial {x_{\alpha_3}}}}+{\epsilon ^4}\mathcal{N}_4^{\alpha_1\alpha_2\alpha_3\alpha_4} \frac{{{\partial ^4 \omega^{(0)}}}}{{\partial {x_{\alpha_1}}\partial {x_{\alpha_2}}\partial {x_{\alpha_3}}\partial {x_{\alpha_4}}}}\Big) \bigg]Z_{\Delta}^{(4\epsilon )}(\bm{x}){\hat n_i}d{\Gamma _{\bm{y}}}\\
        &-\epsilon^{-3}\sum_{\mathbf{z}\in I_{\epsilon}}\int_{\partial E_\mathbf{z}} \frac{\partial }{{\partial {y_j}}}\bigg[ {D_{ijkl}(\bm{y})\frac{{{\partial ^2 }}}{{\partial {y_k}\partial {y_l}}}}\Big(\omega^{(0)} + {\epsilon ^2}\mathcal{N}_2^{\alpha_1\alpha_2}\frac{{{\partial ^2 \omega^{(0)} }}}{{\partial {x_{\alpha_1}}\partial {x_{\alpha_2}}}}\\
&+ {\epsilon ^3}\mathcal{N}_3^{\alpha_1\alpha_2\alpha_3}\frac{{{\partial ^3 \omega^{(0)}}}}{{\partial {x_{\alpha_1}}\partial {x_{\alpha_2}}\partial {x_{\alpha_3}}}}+{\epsilon ^4}\mathcal{N}_4^{\alpha_1\alpha_2\alpha_3\alpha_4} \frac{{{\partial ^4 \omega^{(0)}}}}{{\partial {x_{\alpha_1}}\partial {x_{\alpha_2}}\partial {x_{\alpha_3}}\partial {x_{\alpha_4}}}}\Big) \bigg]Z_{\Delta}^{(4\epsilon )}(\bm{x}){\hat n_i}d{\Gamma _{\bm{y}}}=0.
	\end{aligned}
\end{equation}
\end{small}
Furthermore, by means of the assumption B$_2$, and the lemmas 1-3, it is not difficult to check that
\begin{equation}
	\begin{aligned}
\sum_{\mathbf{z}\in I_{\epsilon}}\int_{\partial E_\mathbf{z}} {D_{ijkl}^\epsilon (\bm{x})\frac{{{\partial ^2 Z_{\Delta}^{(4\epsilon )}(\bm{x})}}}{{\partial {x_k}\partial {x_l}}}}\frac{\partial Z_{\Delta}^{(4\epsilon )}(\bm{x})}{{\partial {x_i}}}{\hat n_j}d{\Gamma _{\bm{y}}}= 0.
	\end{aligned}
\end{equation}

Substituting the boundary conditions (41)-(44) into equation (40), thus equation (40) can be simplified as
\begin{equation}
	\begin{aligned}
&\int_\Omega {D_{ijkl}^\epsilon (\bm{x})\frac{{{\partial ^2 Z_{\Delta}^{(4\epsilon )}(\bm{x})}}}{{\partial {x_k}\partial {x_l}}}}\frac{\partial^2 Z_{\Delta}^{(4\epsilon )}(\bm{x})}{{\partial {x_i}\partial {x_j}}}d\Omega=\int_\Omega \epsilon {E_{4}(\bm{x},\bm{y})} Z_{\Delta}^{(4\epsilon )}(\bm{x})d\Omega.
	\end{aligned}
\end{equation}
Furthermore, on the basis of assumptions A$_1$-A$_3$, and the priori estimate (37) of fourth-order equation, the following inequality holds
\begin{equation}
{\left\| Z_{\Delta}^{(4\epsilon )}(\bm{x})\right\|_{{{{H_0^2}(\Omega )}}}}\leq C_{\Omega}{\left\|\epsilon {E_{4}(\bm{x},\bm{y})} \right\|_{{H^{-2}}(\Omega )}}.
\end{equation}
Finally, we obtain that the explicit error estimate in the integral sense for the higher-order multi-scale solution of the fourth-order multi-scale plate problem (1) as below.
\begin{equation}
\left\|\omega^\epsilon(\bm{x}) - \omega^{(4\epsilon )}(\bm{x})\right\|_{H_0^2(\Omega)}\le C_{\Omega}\epsilon.
\end{equation}
Hence the proof of Theorem 2 is complete.
\section{Multi-scale numerical algorithm}
Based on the theoretical results established in the previous chapters, we now develop the multi-scale numerical algorithm for fourth-order Kirchhoff plate problem in composite plates (1). The detailed multi-scale algorithm is illustrated as below.
\begin{enumerate}
\item[Step 1.]
Identify the geometric configuration of PUC $\mathbf{Y}=[0,1]^2$ and generate a family of triangular finite element meshes $\mathcal{T}_{h_1}=\{K\}$ for PUC $\mathbf{Y}$, where $h_1=$max$_K\{h_K\}$.
\item[Step 2.]
Define the Morley finite element space $\mathbb{P}_{h_1}^2(\mathbf{Y})$ or Hsieh-Clough-Tocher (HCT) finite element space $\mathbb{P}_{h_1}^{\text{HCT}}(\mathbf{Y})$ for solving auxiliary cell problems. Next, employ FEM to solve the second-order cell functions defined by (16) on $\mathbb{P}_{h_1}^2(\mathbf{Y})$ or $\mathbb{P}_{h_1}^{\text{HCT}}(\mathbf{Y})$. Note that classical periodic boundary condition of auxiliary cell problems is replaced by clamped boundary condition for practical numerical implementation. The detailed computational scheme with test function $\upsilon^{h_1}(\bm{y})\in\mathbb{P}_{h_1}^2(\mathbf{Y})/\mathbb{P}_{h_1}^{\text{HCT}}(\mathbf{Y})$ is establishes for solving second-order cell problems (16) as below.
\begin{equation}
\!\!\!\int_\mathbf{Y}\!{{D_{ijkl}}(\bm{y})\frac{{{\partial ^2 \mathcal{N}_2^{\alpha_1\alpha_2}(\bm{y})}}}{{\partial {y_k}\partial {y_l}}}}\!\frac{\partial^2 \upsilon^{h_1}(\bm{y})}{{\partial {y_i}\partial {y_j}}}d\mathbf{Y}\!=\! -\!\int_\mathbf{Y}\!{D_{ij{\alpha _1}{\alpha _2}}}(\bm{y})\!\frac{{\partial ^2 \upsilon^{h_1}(\bm{y})}}{{\partial {y_i}\partial {y_j}}}d\mathbf{Y}.
\end{equation}
\item[Step 3.]
Employ FEM to solve the third-order cell functions defined by (19) on $\mathbb{P}_{h_1}^2(\mathbf{Y})$ or $\mathbb{P}_{h_1}^{\text{HCT}}(\mathbf{Y})$.
\item[Step 4.]
Evaluate the homogenized bending stiffness parameters ${\widehat D}_{ijkl}$ at macro-scale based on the formula (21).
\item[Step 5.]
Identify the geometric configuration of the macroscopic region $\Omega$ and $\mathcal{T}_{h_0}=\{e\}$ be a family of triangular finite element meshes for the macroscopic region $\Omega$, where $h_0=$max$_e\{h_e\}$.
\item[Step 6.]
Define the Morley finite element space $\mathbb{P}_{h_0}^2(\Omega)$ or HCT finite element space $\mathbb{P}_{h_0}^{\text{HCT}}(\Omega)$ \cite{R39} for macroscopic homogenized Kirchhoff plate equation (20). The specific computational scheme with test function $\phi^{h_0}(\bm{x})\in\mathbb{P}_{h_0}^2(\Omega)/\mathbb{P}_{h_0}^{\text{HCT}}(\Omega)$ is establishes for macroscopic homogenized Kirchhoff plate equation (20).
\begin{equation}
\left\{ \begin{aligned}
		& \int_\Omega{\widehat D}_{ijkl}\frac{{{\partial ^2{\omega^{(0)}}(\bm{x})}}}{\partial {x_k}\partial {x_l}}\frac{\partial^2 \phi^{h_0}(\bm{x})}{\partial {x_i}\partial {x_j}}d\Omega=\int_\Omega q(\bm{x})\phi^{h_0}(\bm{x})d\Omega,\\
		&{\omega^{(0)}}(\bm{x}) = g_1(\bm{x})\;\;\text{on}\;\; \partial \Omega,\\
		&\frac{{\partial \omega^{(0)}(\bm{x})}}{{\partial {x_j}}}{\hat n_j} = g_2(\bm{x})\;\;\text{on}\;\; \partial \Omega.
\end{aligned}\right.
\end{equation}
\item[Step 7.]
Employ FEM to solve the fourth-order cell functions defined by (24) on $\mathbb{P}_{h_1}^2(\mathbf{Y})$ or $\mathbb{P}_{h_1}^{\text{HCT}}(\mathbf{Y})$.
\item[Step 8.]
The average approach on relative elements is utilized to compute the spatial derivatives of macroscopic homogenized solution ${\omega^{(0)}}(\bm{x})$ in \cite{R15,R16,R35,R36}. Furthermore, the fourth-order multi-scale solution ${\omega^{(4\epsilon)}}(\bm{x})$ for transverse displacement field of composite Kirchhoff plates is computed by the formula (25).
\end{enumerate}

\section{Numerical experiments and results}
This section presents several numerical examples to validate the computational accuracy and efficiency of the proposed FOMS methodology. The numerical experiments are conducted on a HP desktop workstation equipped with an Intel Core i9-13900H processor (2.60 GHz) and 32.0 GB of internal memory, and all numerical simulations are performed based on Freefem++ software \cite{R38}.

In addition, since it is scarcely impossible to obtain the exact solutions for the plate bending problems of composite Kirchhoff plates, we use the direct numerical simulation (DNS) solution $\omega_{\text{D}}^\epsilon(\bm{x})$ as reference solution for evaluating the computational performance of the proposed FOMS approach. Furthermore, we denote some error notations $e_0(\bm{x})=\omega^0(\bm{x})-\omega_{\text{DNS}}^\epsilon(\bm{x})$, $e_2(\bm{x})=\omega^{(2\epsilon)}(\bm{x})-\omega_{\text{DNS}}^\epsilon(\bm{x})$,  $e_3(\bm{x})=\omega^{(3\epsilon)}(\bm{x})-\omega_{\text{DNS}}^\epsilon(\bm{x})$ and  $e_4(\bm{x})=\omega^{(4\epsilon)}(\bm{x})-\omega_{\text{DNS}}^\epsilon(\bm{x})$ for subsequent error analysis.

\subsection{Validation of the computational performance of FOMS method}
In this example, a composite thin plate with clamped boundary condition is investigated, whose multi-scale structure $\Omega=(x_1,x_2)=[0,1]\times[0,1]$cm$^2$, microscopic unit cell $\mathbf{Y}=(y_1,y_2)=[0,1]\times[0,1]$ and macroscopic homogenization structure are depicted in Fig.~2. This composite plate with the characteristic periodic parameter $\epsilon=1/8$ consists of two component materials, where Young's modulus and Poisson's ratio of the matrix material (yellow color) are 50.0GPa and 0.20, respectively. For the inclusion material (red color), these values are 8.0MPa and 0.20. In addition, the transverse load imposed on this composite plate is set to 1500.0N/cm$^2$.
\begin{figure}[!htb]
\centering
\begin{minipage}[c]{0.32\textwidth}
  \centering
  \includegraphics[width=0.9\linewidth,totalheight=1.5in]{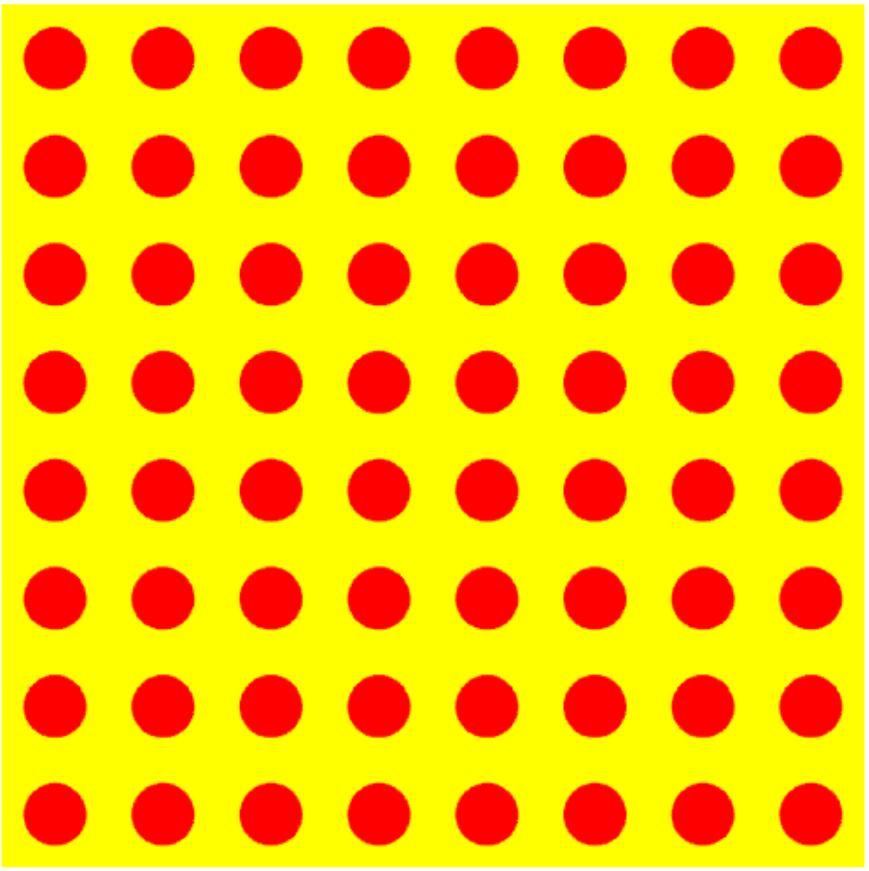} \\
  (a)
\end{minipage}
\begin{minipage}[c]{0.32\textwidth}
  \centering
  \includegraphics[width=0.9\linewidth,totalheight=1.5in]{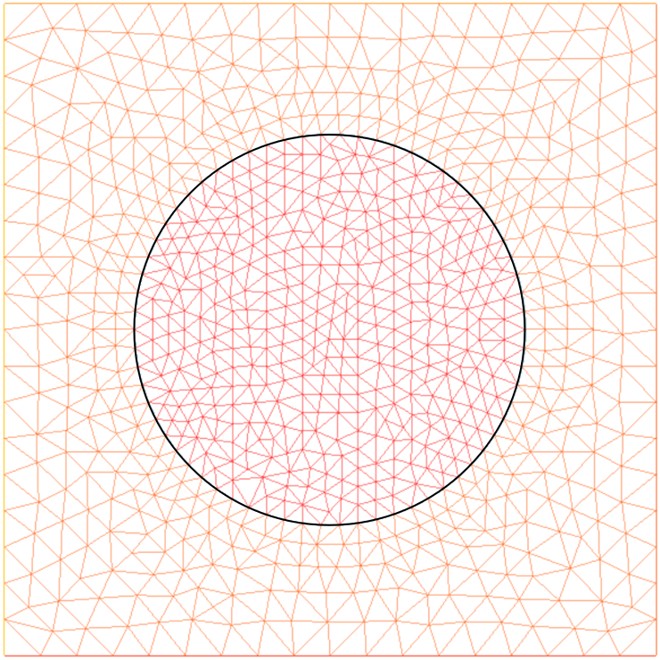} \\
  (b)
\end{minipage}
\begin{minipage}[c]{0.32\textwidth}
  \centering
  \includegraphics[width=0.9\linewidth,totalheight=1.5in]{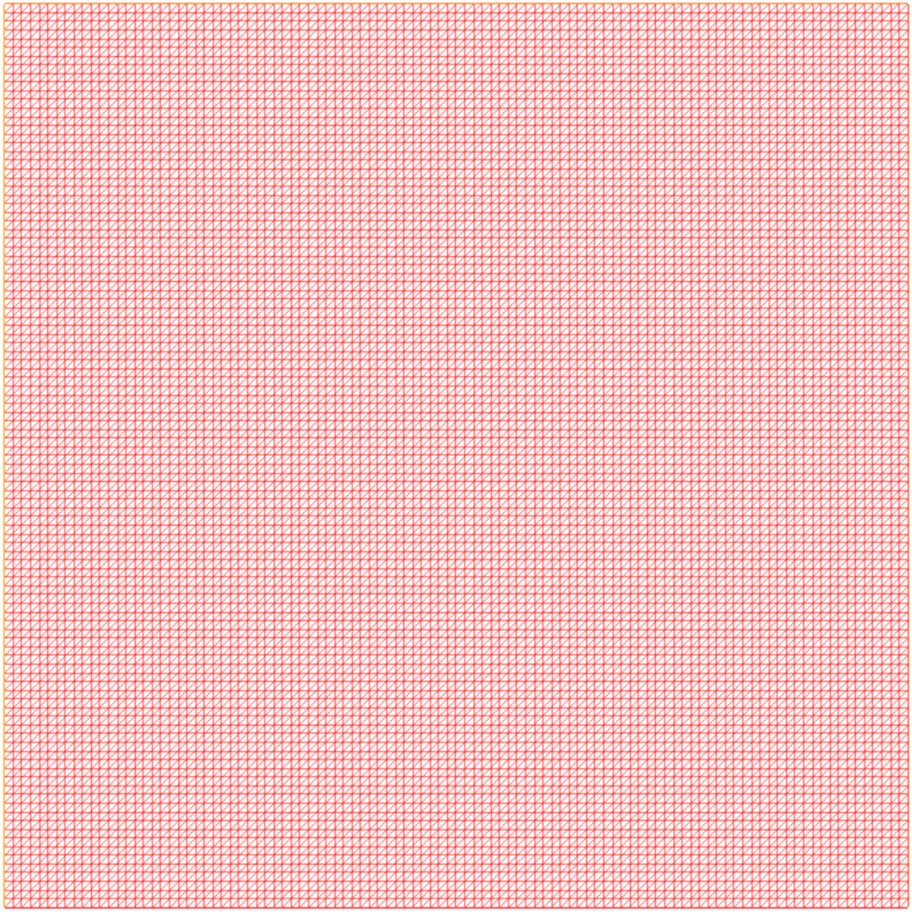} \\
  (c)
\end{minipage}
\caption{(a) The multi-scale structure of composite plate $\Omega$; (b) the microscopic unit cell $\mathbf{Y}$; (c) the macroscopic homogenization structure.}\label{f2}
\end{figure}

Next, the proposed FOMS computational method based on Morley finite element and HCT finite element is employed to simulate this composite Kirchhoff plate. Table 1 presents the detail comparison of computational resource cost of the FOMS method and reference direct numerical simulation.
\begin{table}[!htb]{\caption{Comparison of computational cost.}\label{t2}}
\centering
\begin{tabular}{cccc}
\hline
 & Multi-scale problem & Cell problem & Homogenized problem \\
\hline
Number of elements & 209920 & 1620 & 21632\\
Number of nodes    & 105361 & 847 & 11025\\
\hline
& Direct numerical simulation &  \multicolumn{2}{c}{Multi-scale algorithm} \\
\hline
Time/Morley & 10.419s  & \multicolumn{2}{c}{7.594s}\\
\hline
Time/HCT & 91.642s  & \multicolumn{2}{c}{20.552s}\\
\hline
\end{tabular}
\end{table}

Furthermore, Table 2 displays the detailed computational accuracy comparison of of the FOMS method and reference direct numerical simulation.
\begin{table}[h]{\caption{Comparison of computational accuracy.}\label{t2}}
\centering
\begin{tabular}{ccccccccc}
\hline
Morley & $\frac{||e_0||_{L^2}}{||\omega_{\text{D}}^\epsilon||_{L^2}}$ & $\frac{||e_2||_{L^2}}{||\omega_{\text{D}}^\epsilon||_{L^2}}$ & $\frac{||e_3||_{L^2}}{||\omega_{\text{D}}^\epsilon||_{L^2}}$ & $\frac{||e_4||_{L^2}}{||\omega_{\text{D}}^\epsilon||_{L^2}}$ & $\frac{|e_0|_{H^1}}{|\omega_{\text{D}}^\epsilon|_{H^1}}$ & $\frac{|e_2|_{H^1}}{|\omega_{\text{D}}^\epsilon|_{H^1}}$ & $\frac{|e_3|_{H^1}}{|\omega_{\text{D}}^\epsilon|_{H^1}}$ & $\frac{|e_4|_{H^1}}{|\omega_{\text{D}}^\epsilon|_{H^1}}$\\
\hline
Percentage \% & 5.001 & 4.854 & 4.858 & 1.464 & 56.760 & 56.375 & 56.346 & 16.527\\
\hline
HCT & $\frac{||e_0||_{L^2}}{||\omega_{\text{D}}^\epsilon||_{L^2}}$ & $\frac{||e_2||_{L^2}}{||\omega_{\text{D}}^\epsilon||_{L^2}}$ & $\frac{||e_3||_{L^2}}{||\omega_{\text{D}}^\epsilon||_{L^2}}$ & $\frac{||e_4||_{L^2}}{||\omega_{\text{D}}^\epsilon||_{L^2}}$ & $\frac{|e_0|_{H^1}}{|\omega_{\text{D}}^\epsilon|_{H^1}}$ & $\frac{|e_2|_{H^1}}{|\omega_{\text{D}}^\epsilon|_{H^1}}$ & $\frac{|e_3|_{H^1}}{|\omega_{\text{D}}^\epsilon|_{H^1}}$ & $\frac{|e_4|_{H^1}}{|\omega_{\text{D}}^\epsilon|_{H^1}}$\\
\hline
Percentage \% & 4.845 & 4.718 & 4.722 & 1.565 & 56.346 & 55.986 & 55.957 & 24.837\\
\hline
\end{tabular}
\end{table}

Figs.~3 and 4 are plotted to visualize the simulative results for solutions $\omega^{(0)}$, $\omega^{(2\epsilon)}$, $\omega^{(3\epsilon)}$, $\omega^{(4\epsilon)}$ and $\omega^\epsilon_{\text{D}}$ based on Morley finite element and HCT finite element, respectively.
\begin{figure}[!htb]
\centering
\begin{minipage}[c]{0.32\textwidth}
  \centering
  \includegraphics[width=45mm]{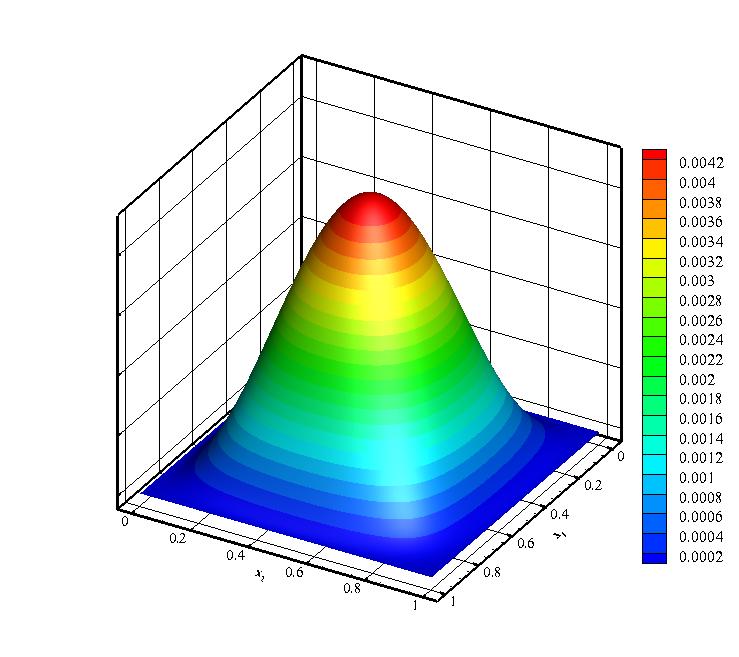}\\
  (a)
\end{minipage}
\begin{minipage}[c]{0.32\textwidth}
  \centering
  \includegraphics[width=45mm]{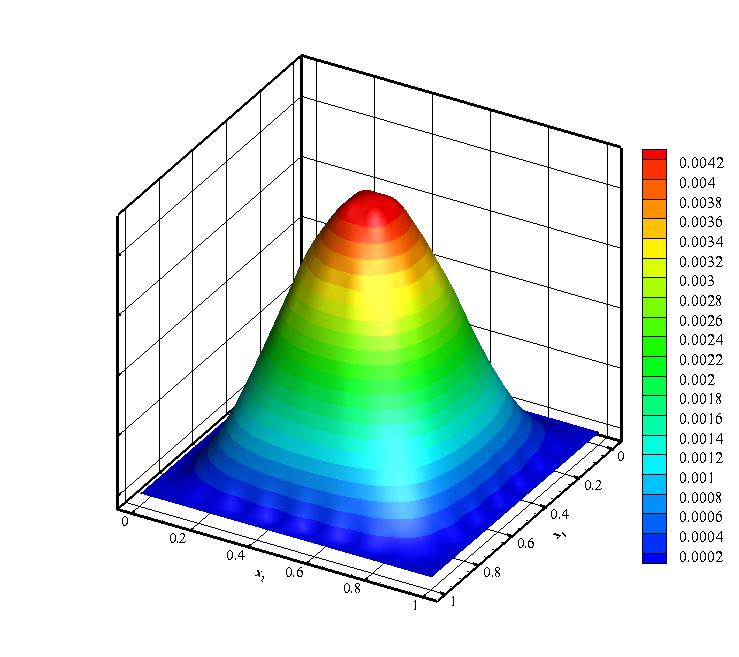}\\
  (b)
\end{minipage}
\begin{minipage}[c]{0.32\textwidth}
  \centering
  \includegraphics[width=45mm]{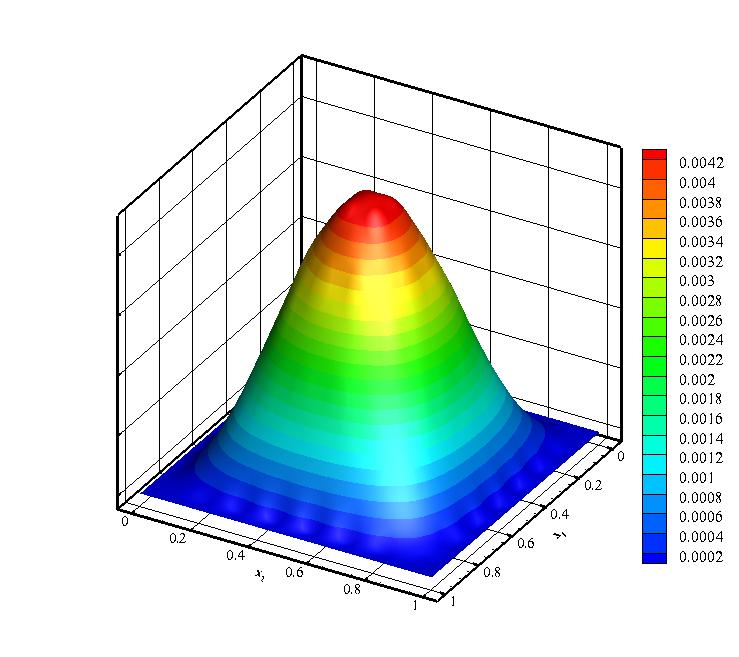}\\
  (c)
\end{minipage}
\begin{minipage}[c]{0.32\textwidth}
  \centering
  \includegraphics[width=45mm]{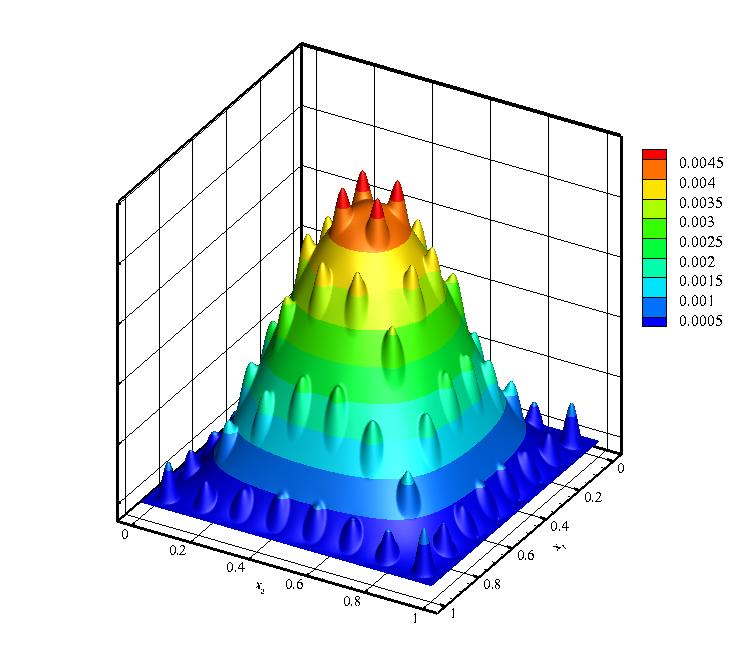}\\
  (d)
\end{minipage}
\begin{minipage}[c]{0.32\textwidth}
  \centering
  \includegraphics[width=45mm]{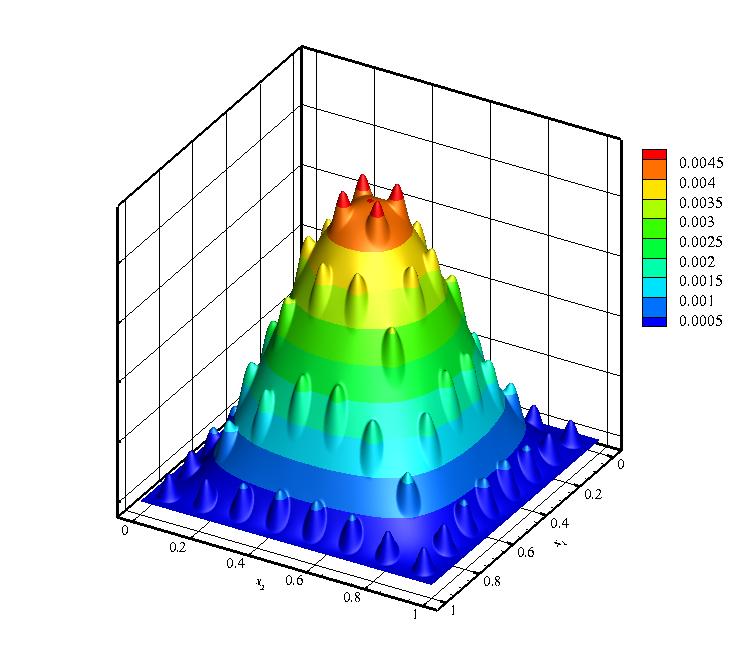}\\
  (e)
\end{minipage}
\caption{The transverse displacement of composite Kirchhoff plate computed by Morley finite element: (a) $\omega^{(0)}$; (b) $\omega^{(2\epsilon)}$; (c) $\omega^{(3\epsilon)}$; (d) $\omega^{(4\epsilon)}$; (e) $\omega^\epsilon_{\text{D}}$.}\label{f8}
\end{figure}
\begin{figure}[!htb]
\centering
\begin{minipage}[c]{0.32\textwidth}
  \centering
  \includegraphics[width=45mm]{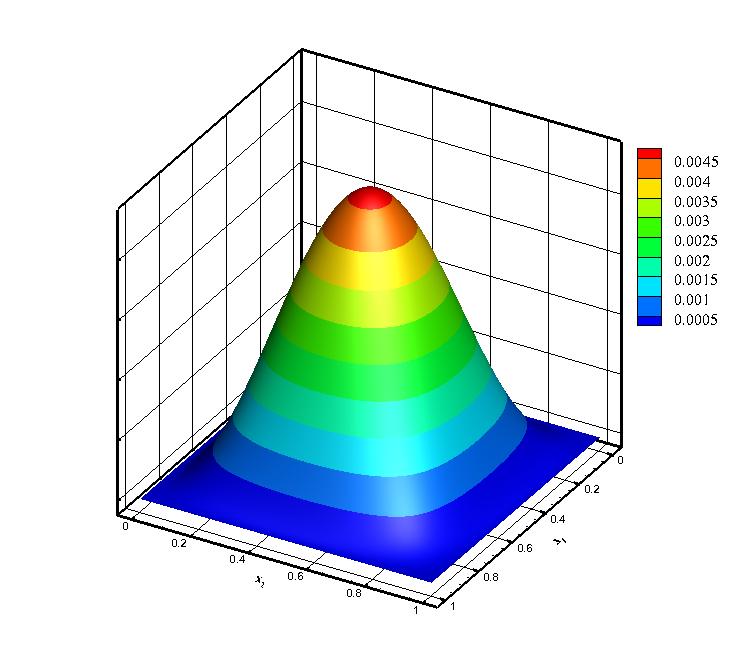}\\
  (a)
\end{minipage}
\begin{minipage}[c]{0.32\textwidth}
  \centering
  \includegraphics[width=45mm]{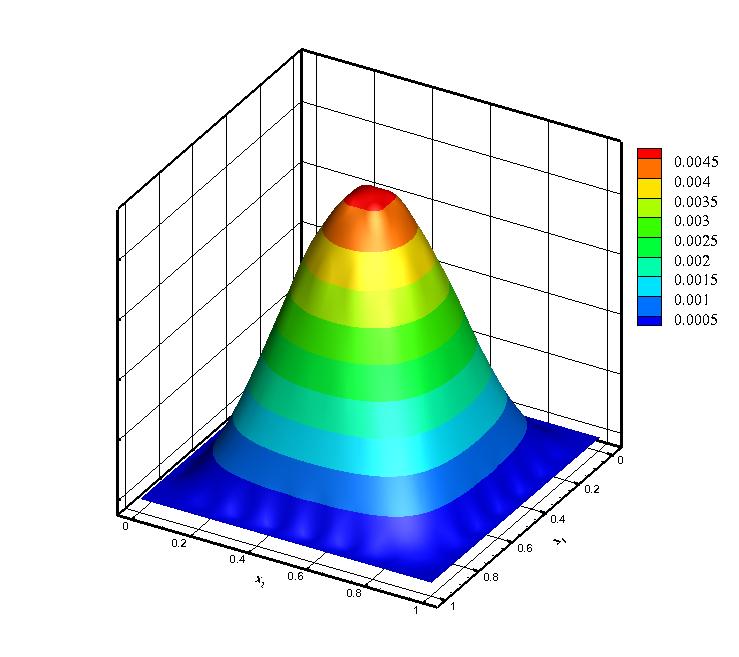}\\
  (b)
\end{minipage}
\begin{minipage}[c]{0.32\textwidth}
  \centering
  \includegraphics[width=45mm]{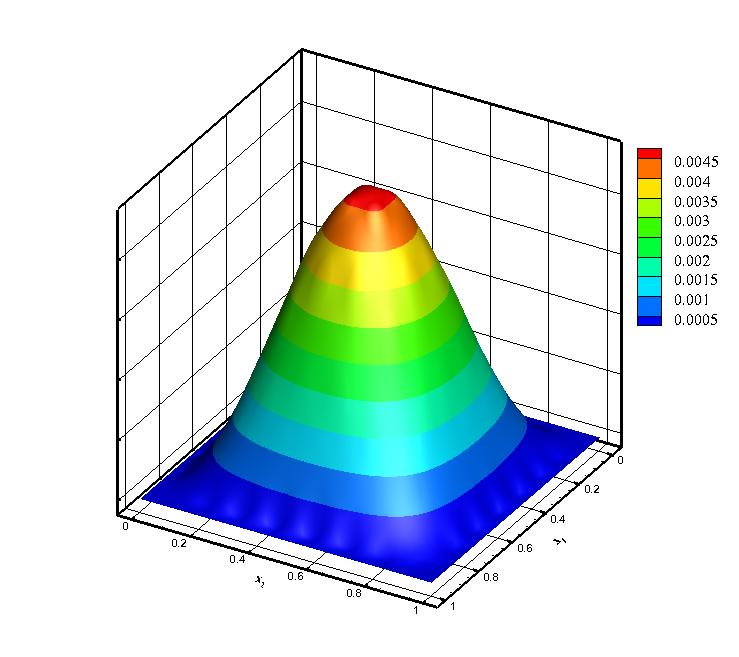}\\
  (c)
\end{minipage}
\begin{minipage}[c]{0.32\textwidth}
  \centering
  \includegraphics[width=45mm]{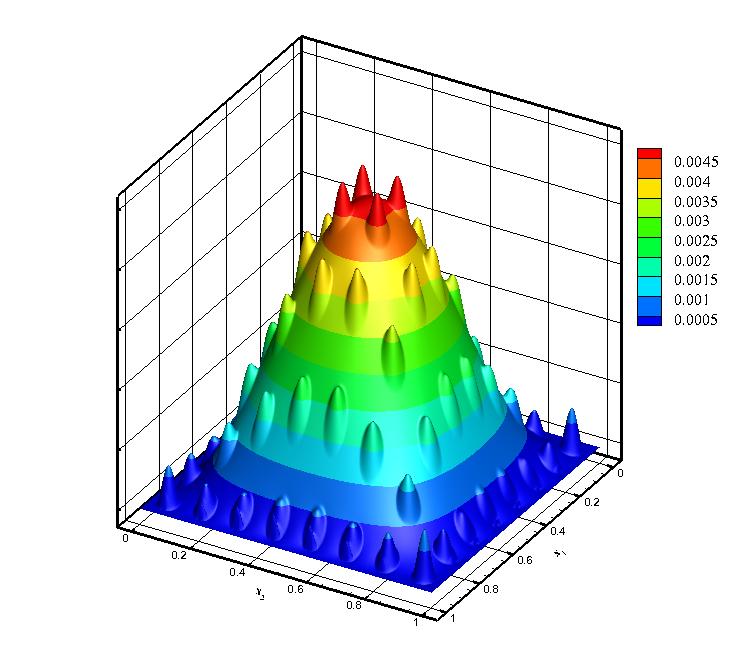}\\
  (d)
\end{minipage}
\begin{minipage}[c]{0.32\textwidth}
  \centering
  \includegraphics[width=45mm]{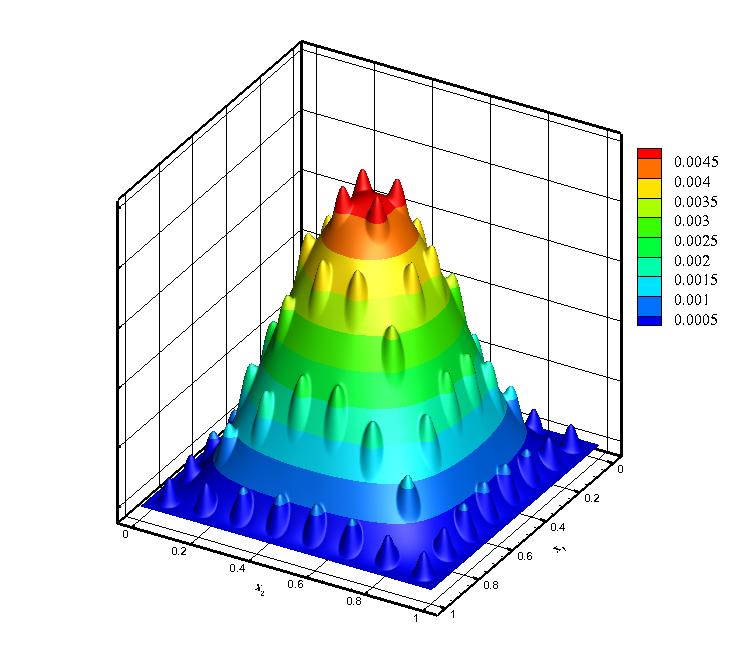}\\
  (e)
\end{minipage}
\caption{The transverse displacement of composite Kirchhoff plate computed by HCT finite element: (a) $\omega^{(0)}$; (b) $\omega^{(2\epsilon)}$; (c) $\omega^{(3\epsilon)}$; (d) $\omega^{(4\epsilon)}$; (e) $\omega^\epsilon_{\text{D}}$.}\label{f8}
\end{figure}

Beside, Figs.~5-8 display the simulative results for the numerical gradients of solutions $\omega^{(0)}$, $\omega^{(2\epsilon)}$, $\omega^{(3\epsilon)}$, $\omega^{(4\epsilon)}$ and $\omega^\epsilon_{\text{D}}$ based on Morley finite element and HCT finite element, respectively.
\begin{figure}[!htb]
\centering
\begin{minipage}[c]{0.32\textwidth}
  \centering
  \includegraphics[width=45mm]{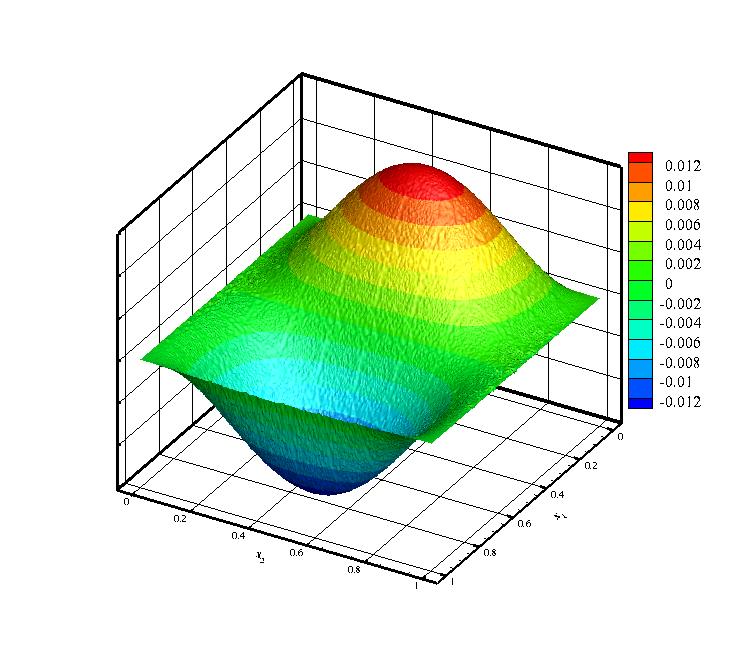}\\
  (a)
\end{minipage}
\begin{minipage}[c]{0.32\textwidth}
  \centering
  \includegraphics[width=45mm]{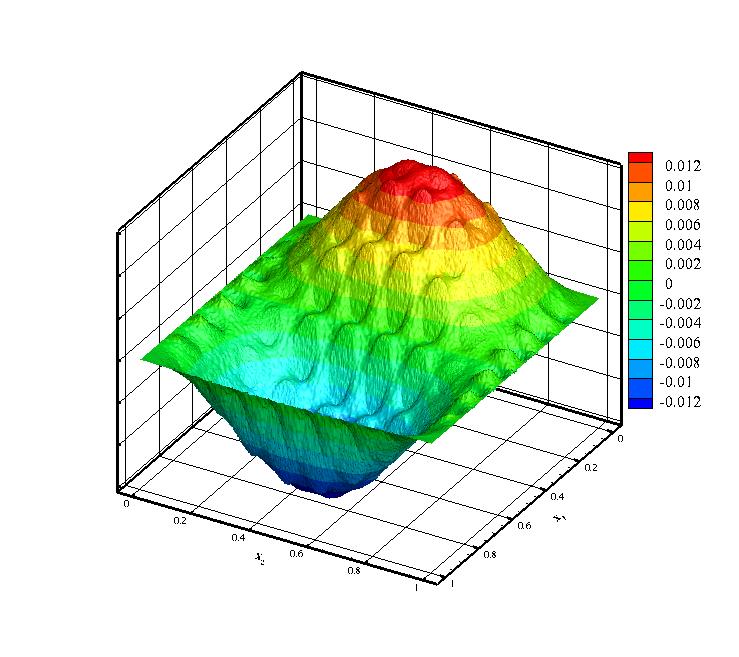}\\
  (b)
\end{minipage}
\begin{minipage}[c]{0.32\textwidth}
  \centering
  \includegraphics[width=45mm]{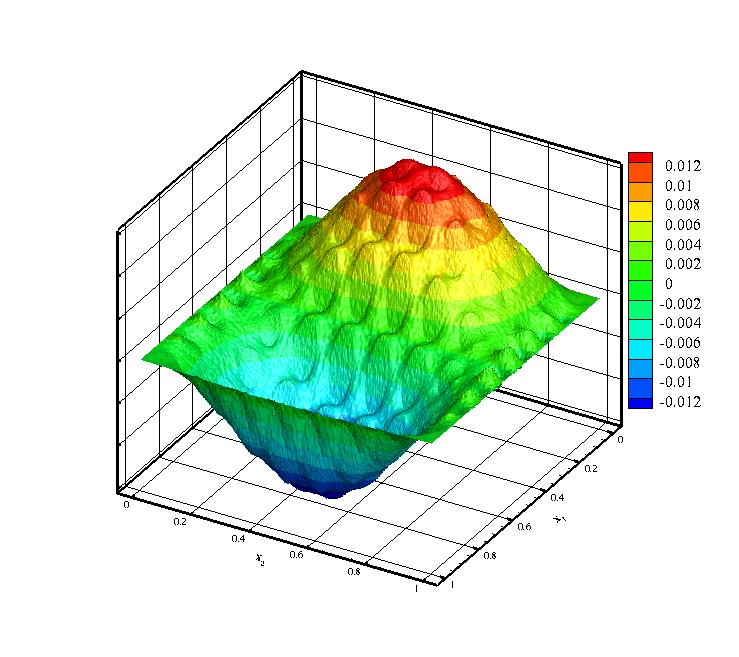}\\
  (c)
\end{minipage}
\begin{minipage}[c]{0.32\textwidth}
  \centering
  \includegraphics[width=45mm]{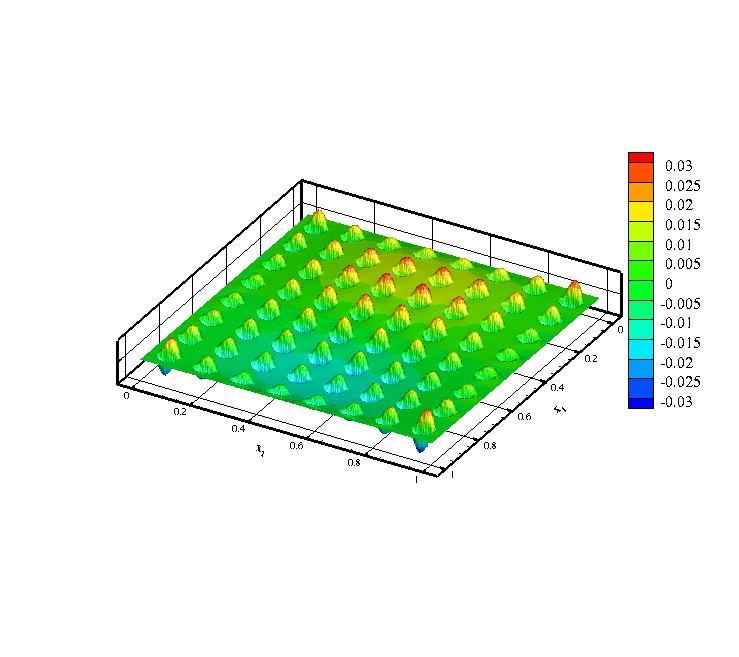}\\
  (d)
\end{minipage}
\begin{minipage}[c]{0.32\textwidth}
  \centering
  \includegraphics[width=45mm]{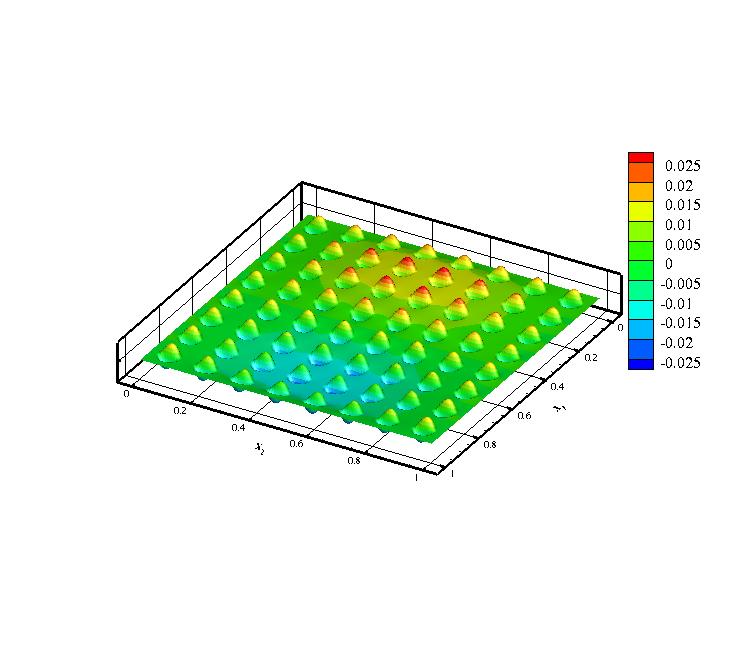}\\
  (e)
\end{minipage}
\caption{The $x_1$-direction gradient of transverse displacement of composite Kirchhoff plate computed by Morley finite element: (a) $\frac{\partial\omega^{(0)}}{\partial x_1}$; (b) $\frac{\partial\omega^{(2\epsilon)}}{\partial x_1}$; (c) $\frac{\partial\omega^{(3\epsilon)}}{\partial x_1}$; (d) $\frac{\partial\omega^{(4\epsilon)}}{\partial x_1}$; (e) $\frac{\partial\omega^\epsilon_{\text{D}}}{\partial x_1}$.}\label{f8}
\end{figure}
\begin{figure}[!htb]
\centering
\begin{minipage}[c]{0.32\textwidth}
  \centering
  \includegraphics[width=45mm]{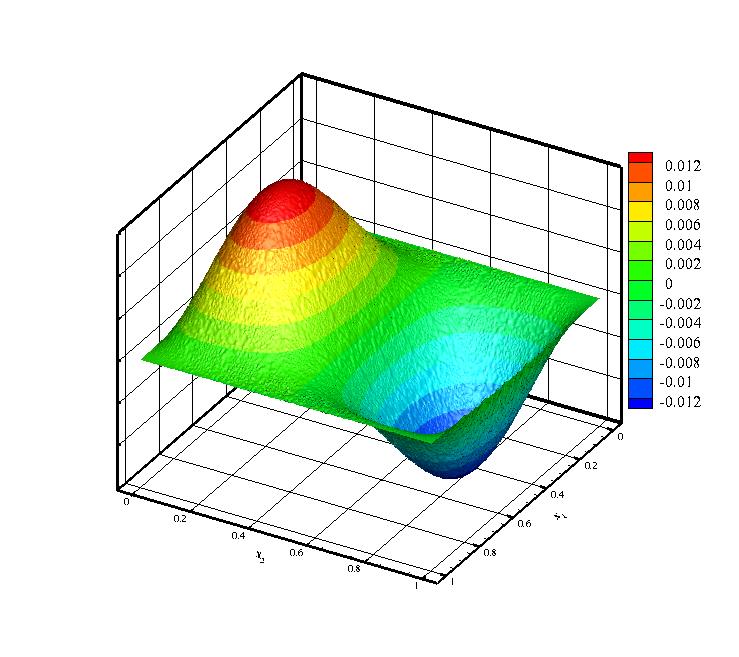}\\
  (a)
\end{minipage}
\begin{minipage}[c]{0.32\textwidth}
  \centering
  \includegraphics[width=45mm]{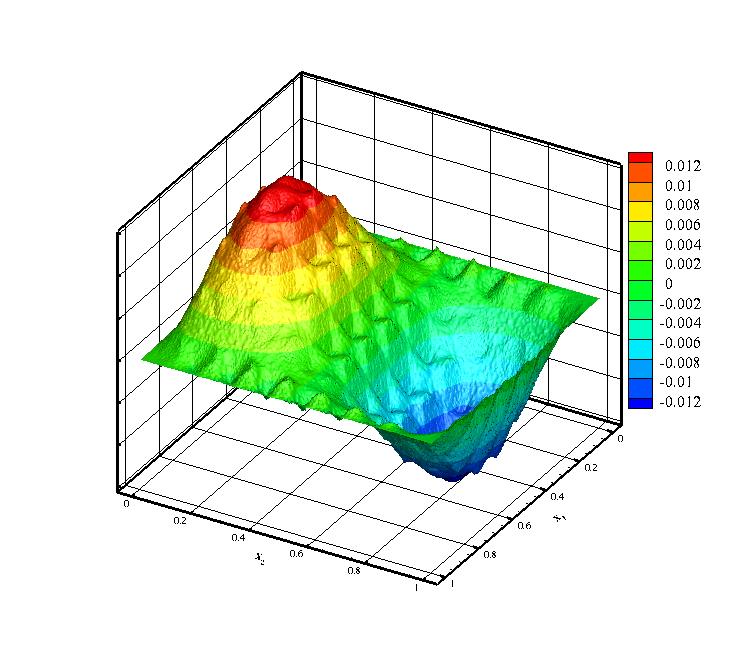}\\
  (b)
\end{minipage}
\begin{minipage}[c]{0.32\textwidth}
  \centering
  \includegraphics[width=45mm]{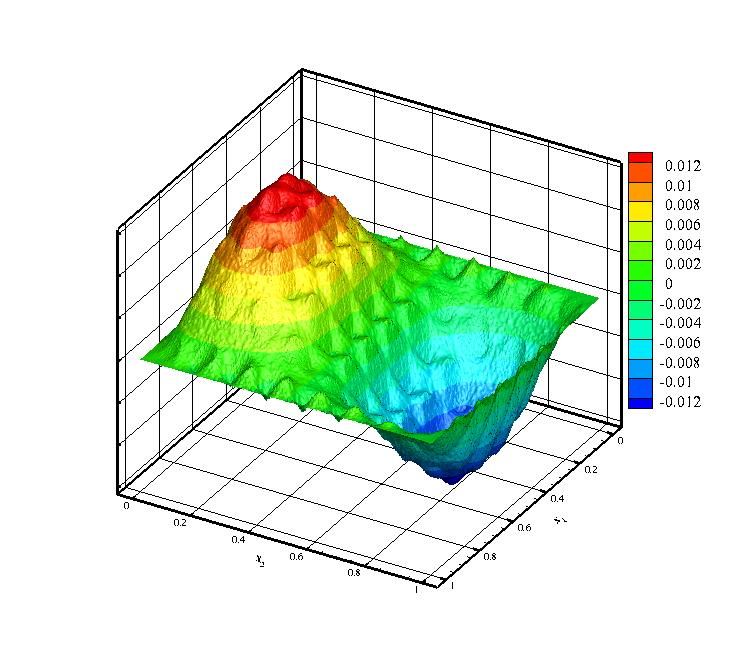}\\
  (c)
\end{minipage}
\begin{minipage}[c]{0.32\textwidth}
  \centering
  \includegraphics[width=45mm]{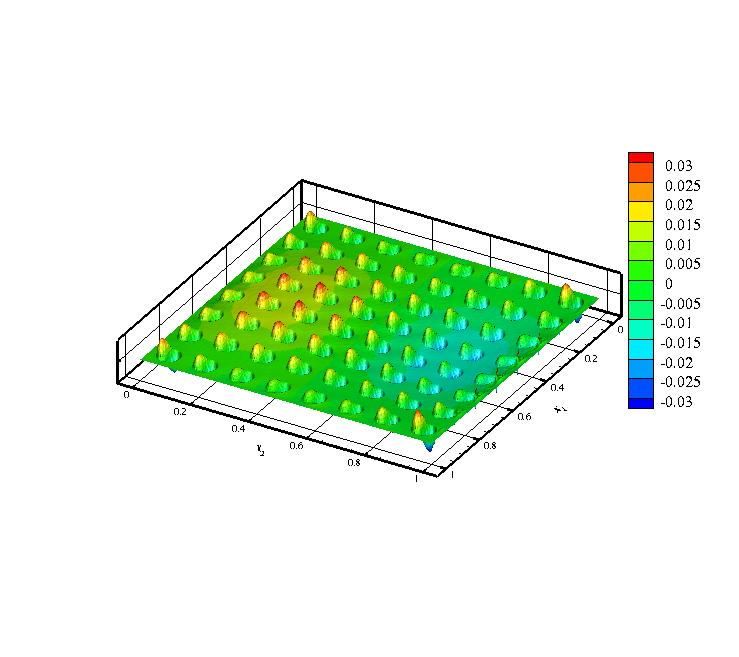}\\
  (d)
\end{minipage}
\begin{minipage}[c]{0.32\textwidth}
  \centering
  \includegraphics[width=45mm]{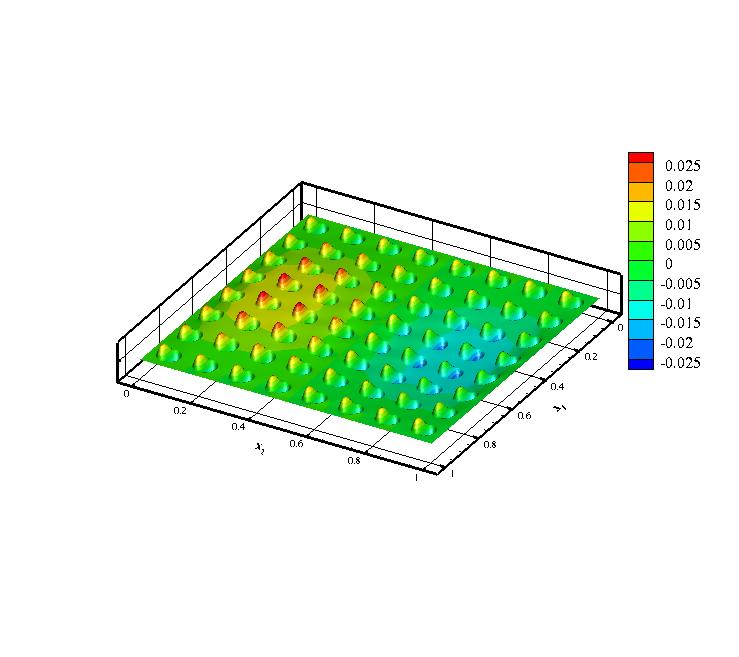}\\
  (e)
\end{minipage}
\caption{The $x_2$-direction gradient of transverse displacement of composite Kirchhoff plate computed by Morley finite element: (a) $\frac{\partial\omega^{(0)}}{\partial x_2}$; (b) $\frac{\partial\omega^{(2\epsilon)}}{\partial x_2}$; (c) $\frac{\partial\omega^{(3\epsilon)}}{\partial x_2}$; (d) $\frac{\partial\omega^{(4\epsilon)}}{\partial x_2}$; (e) $\frac{\partial\omega^\epsilon_{\text{D}}}{\partial x_2}$.}\label{f8}
\end{figure}
\begin{figure}[!htb]
\centering
\begin{minipage}[c]{0.32\textwidth}
  \centering
  \includegraphics[width=45mm]{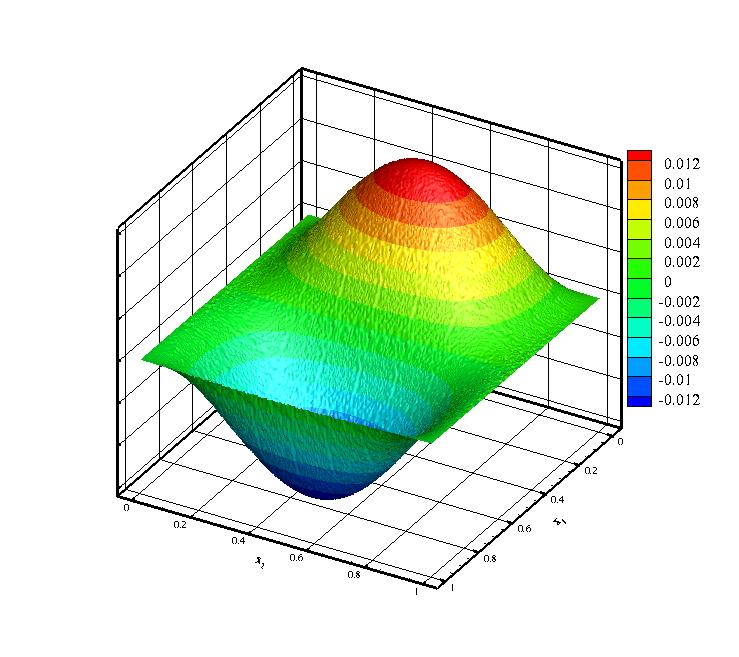}\\
  (a)
\end{minipage}
\begin{minipage}[c]{0.32\textwidth}
  \centering
  \includegraphics[width=45mm]{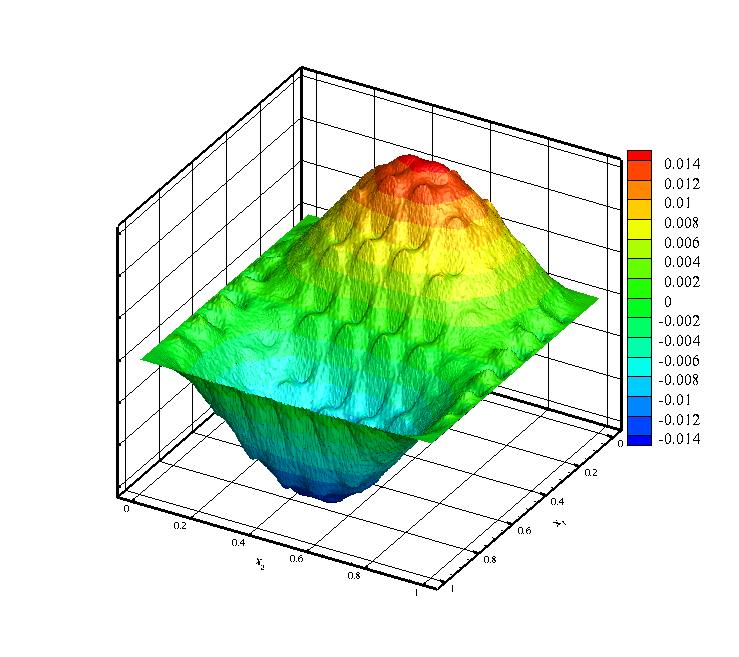}\\
  (b)
\end{minipage}
\begin{minipage}[c]{0.32\textwidth}
  \centering
  \includegraphics[width=45mm]{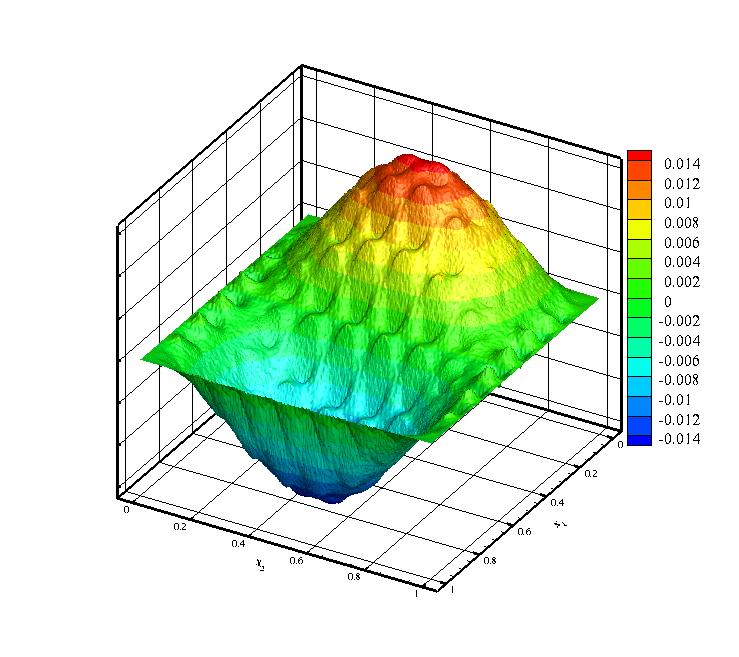}\\
  (c)
\end{minipage}
\begin{minipage}[c]{0.32\textwidth}
  \centering
  \includegraphics[width=45mm]{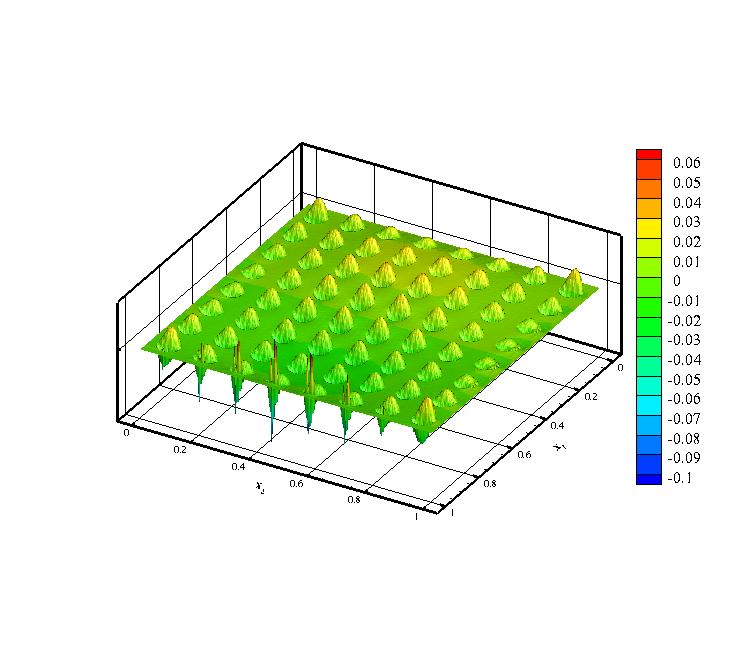}\\
  (d)
\end{minipage}
\begin{minipage}[c]{0.32\textwidth}
  \centering
  \includegraphics[width=45mm]{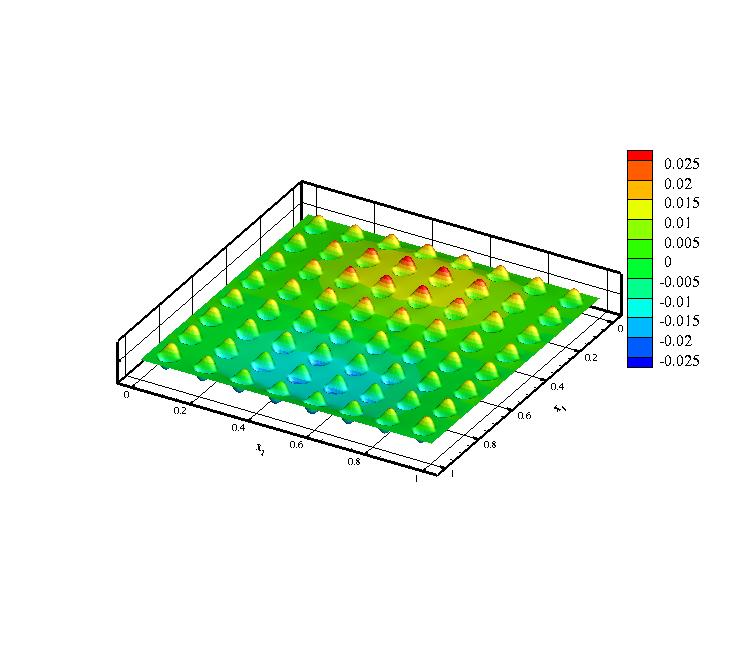}\\
  (e)
\end{minipage}
\caption{The $x_1$-direction gradient of transverse displacement of composite Kirchhoff plate computed by HCT finite element: (a) $\frac{\partial\omega^{(0)}}{\partial x_1}$; (b) $\frac{\partial\omega^{(2\epsilon)}}{\partial x_1}$; (c) $\frac{\partial\omega^{(3\epsilon)}}{\partial x_1}$; (d) $\frac{\partial\omega^{(4\epsilon)}}{\partial x_1}$; (e) $\frac{\partial\omega^\epsilon_{\text{D}}}{\partial x_1}$.}\label{f8}
\end{figure}
\begin{figure}[!htb]
\centering
\begin{minipage}[c]{0.32\textwidth}
  \centering
  \includegraphics[width=45mm]{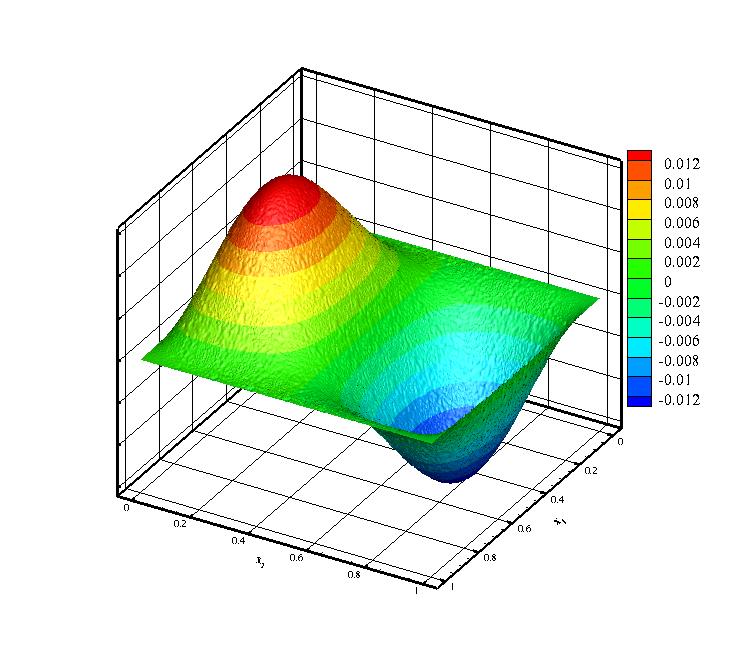}\\
  (a)
\end{minipage}
\begin{minipage}[c]{0.32\textwidth}
  \centering
  \includegraphics[width=45mm]{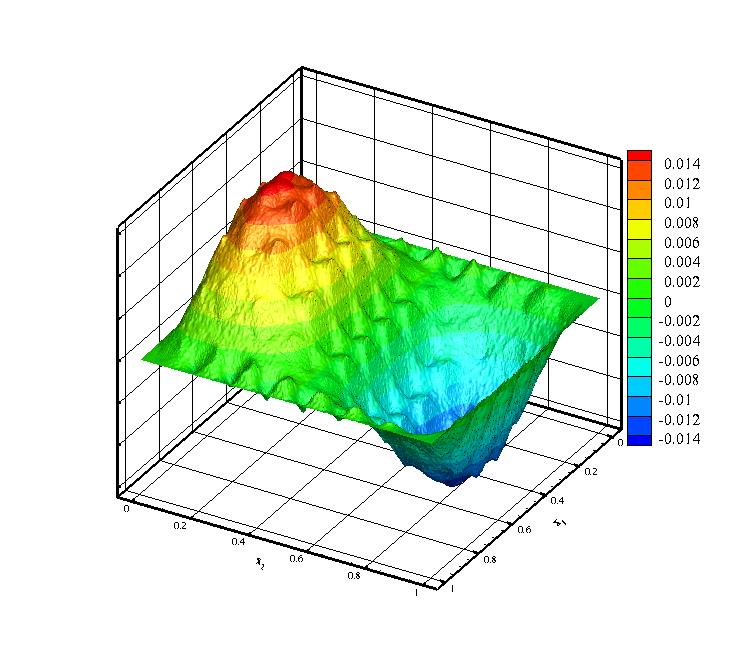}\\
  (b)
\end{minipage}
\begin{minipage}[c]{0.32\textwidth}
  \centering
  \includegraphics[width=45mm]{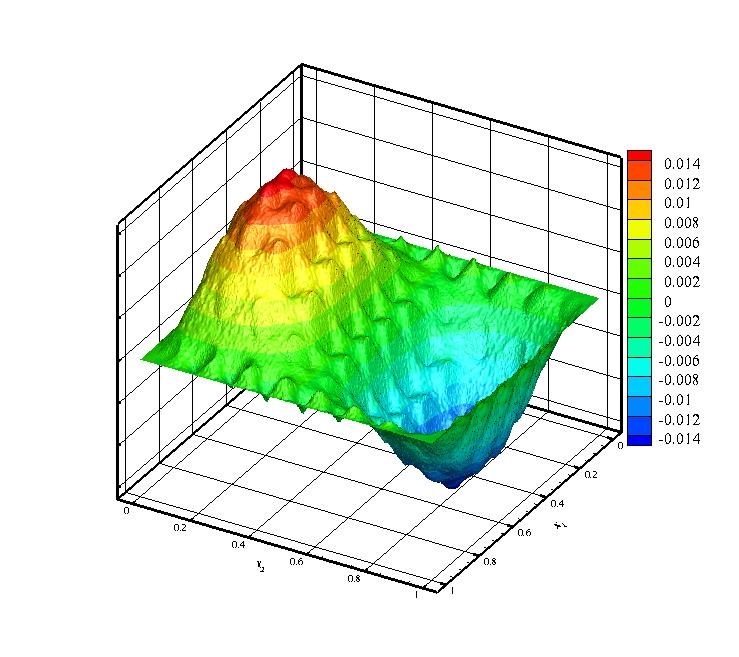}\\
  (c)
\end{minipage}
\begin{minipage}[c]{0.32\textwidth}
  \centering
  \includegraphics[width=45mm]{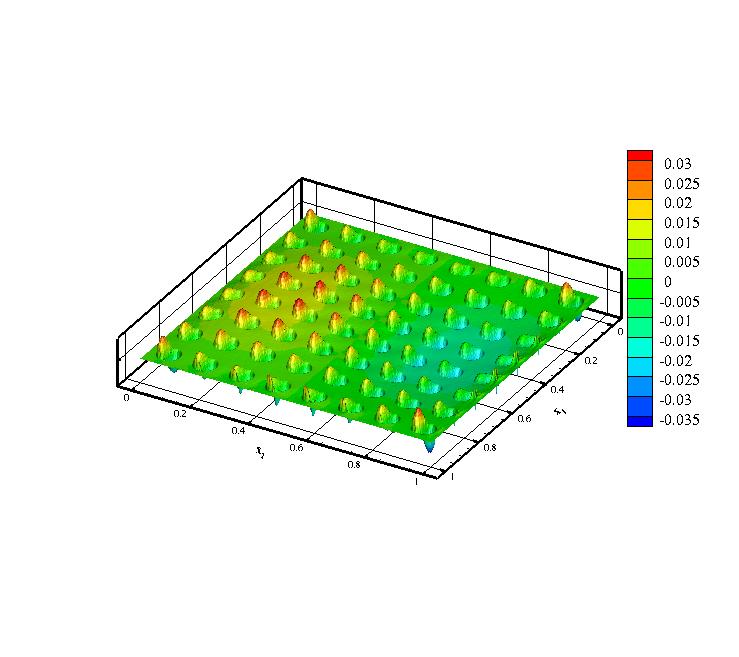}\\
  (d)
\end{minipage}
\begin{minipage}[c]{0.32\textwidth}
  \centering
  \includegraphics[width=45mm]{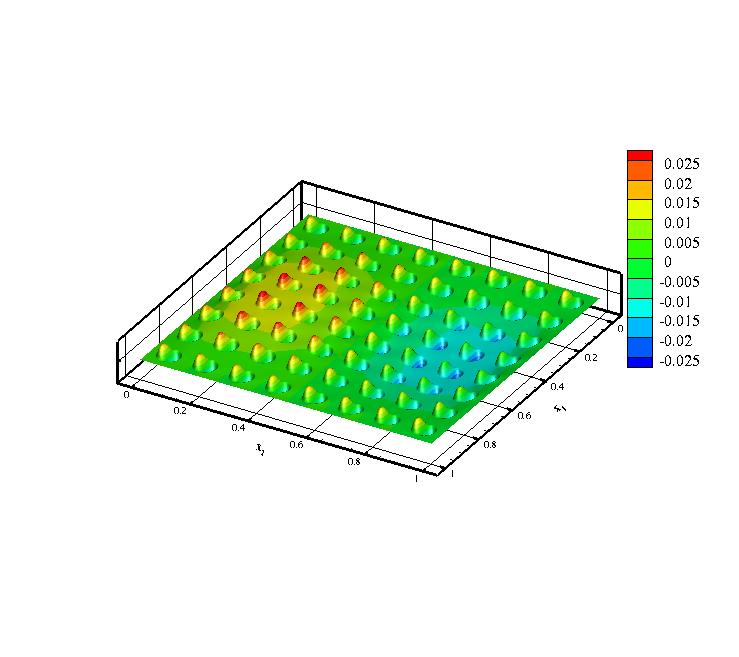}\\
  (e)
\end{minipage}
\caption{The $x_2$-direction gradient of transverse displacement of composite Kirchhoff plate computed by HCT finite element: (a) $\frac{\partial\omega^{(0)}}{\partial x_2}$; (b) $\frac{\partial\omega^{(2\epsilon)}}{\partial x_2}$; (c) $\frac{\partial\omega^{(3\epsilon)}}{\partial x_2}$; (d) $\frac{\partial\omega^{(4\epsilon)}}{\partial x_2}$; (e) $\frac{\partial\omega^\epsilon_{\text{D}}}{\partial x_2}$.}\label{f8}
\end{figure}

As plainly shown in Table 1, the FOMS approach proposed herein demonstrates a substantial reduction in computational resource, particularly in terms of CPU memory and time, when contrasted with high-resolution DNS. Specifically, the proposed FOMS approach achieves an approximate saving of 27.11\% and 77.57\% in computational time when employing Morley and HCT finite elements, respectively. Furthermore, the numerical accuracy analysis in Table 2 clearly illustrates that the FOMS solution achieves a significant enhancement in computational accuracy compared with macroscopic homogenized solution, second-order multi-scale solution and third-order multi-scale solution, especially for H$^1$ semi-norm. Moreover, we can conclude that the computational efficiency and accuracy of the FOMS approach based on Morley finite element is superior to those of FOMS approach based on HCT finite element. Similarly, the computational results in Figs.~3 and 4 also demonstrate only FOMS solution have the capacity to accurately capture the highly microscopic oscillation of composite thin plate, and the computational accuracy of macroscopic homogenized, second-order and third-order solutions is significantly inferior. Additionally,the numerical results in Figs.~5-8 demonstrate the proposed FOMS method can accurately capture the gradient behaviors of composite thin plate. In addition, the numerical error observed in the gradient results from the HCT element at the plate boundary, when using the same mesh, are likely due to inaccurate computation of higher-order derivatives for the macroscopic homogenized solution. Hence, the FOMS method based on Morley finite element is employed to the numerical simulation in the subsequent numerical examples. 
\subsection{Application to composite Kirchhoff plates with different kinds of inclusions}
In this example, other four kinds of composite Kirchhoff plates with clamped boundary condition are simulated, whose geometric configurations of microscopic PUCs are illustrated in Fig.~9. The investigated composite plates with the characteristic periodic parameter $\epsilon=1/8$ are composed of two constituent materials. And the Young's modulus and Poisson's ratio of the matrix material are 60.0GPa and 0.35, respectively. For the inclusion material, these values are 2.4MPa and 0.35. In addition, the transverse load imposed on this composite plate is set to 2500.0N/cm$^2$.
\begin{figure}[!htb]
\centering
\begin{minipage}[c]{0.24\textwidth}
  \centering
  \includegraphics[width=0.9\linewidth,totalheight=1.2in]{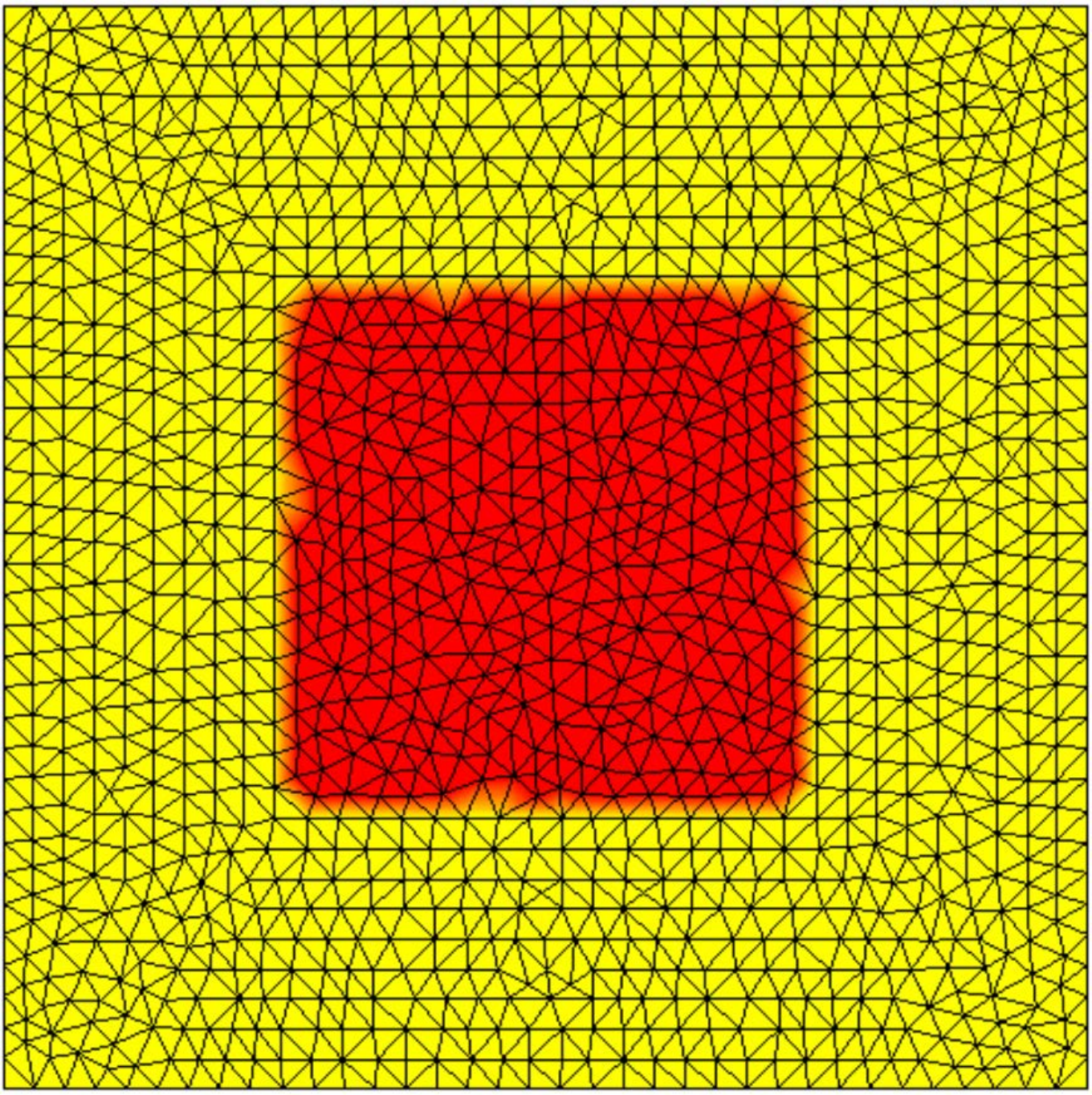} \\
  (a)
\end{minipage}
\begin{minipage}[c]{0.24\textwidth}
  \centering
  \includegraphics[width=0.9\linewidth,totalheight=1.2in]{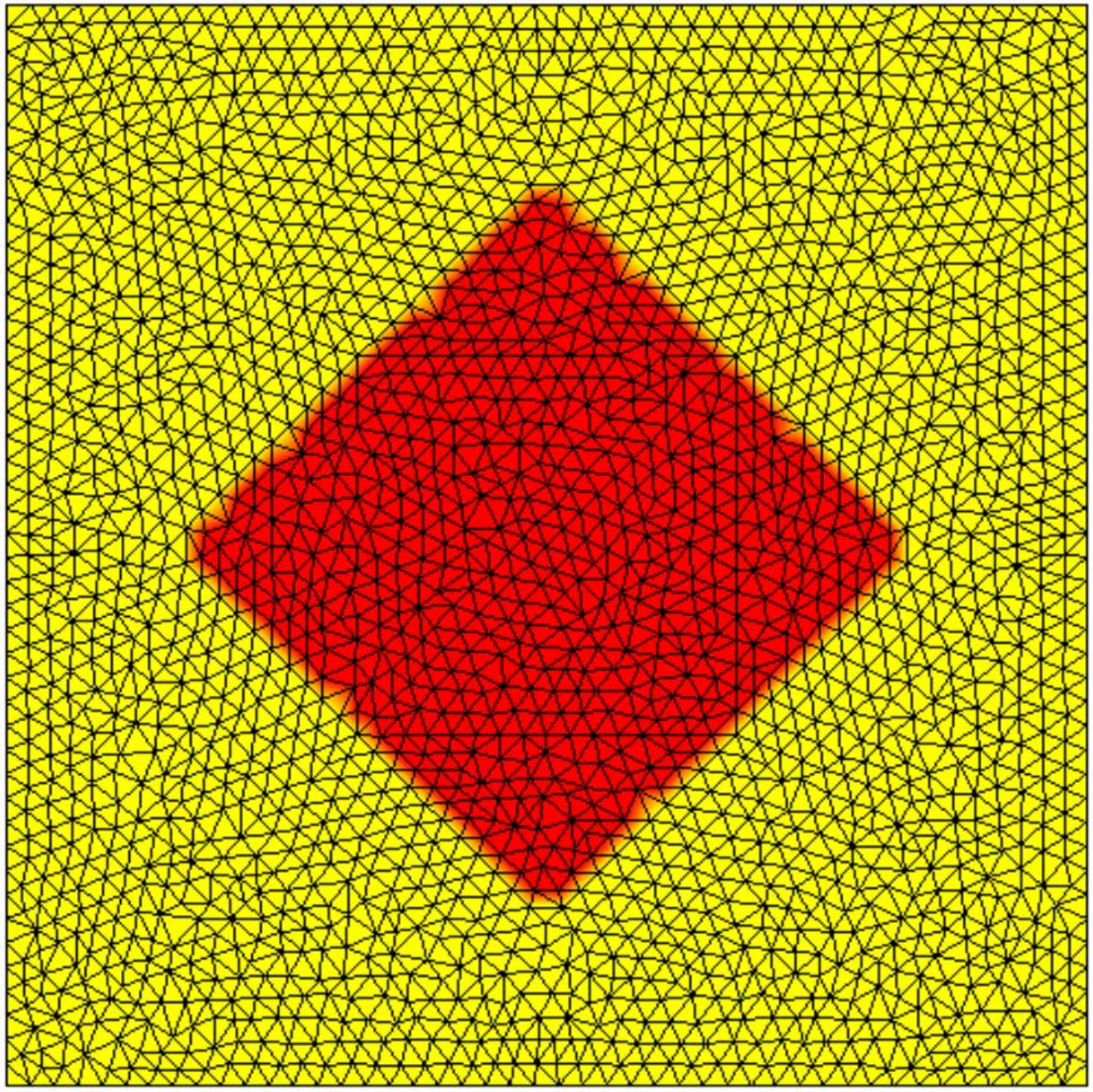} \\
  (b)
\end{minipage}
\begin{minipage}[c]{0.24\textwidth}
  \centering
  \includegraphics[width=0.9\linewidth,totalheight=1.2in]{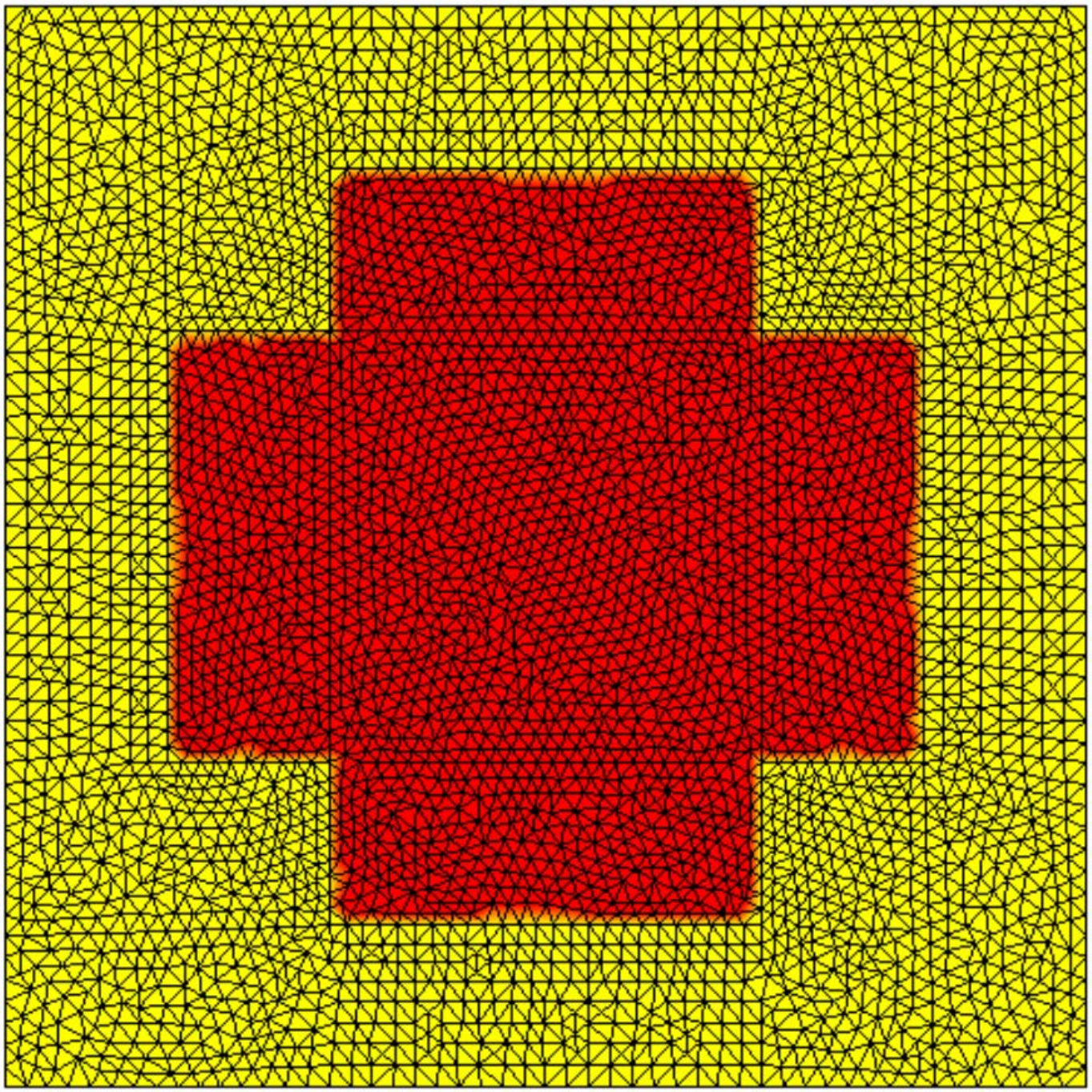} \\
  (c)
\end{minipage}
\begin{minipage}[c]{0.24\textwidth}
  \centering
  \includegraphics[width=0.9\linewidth,totalheight=1.2in]{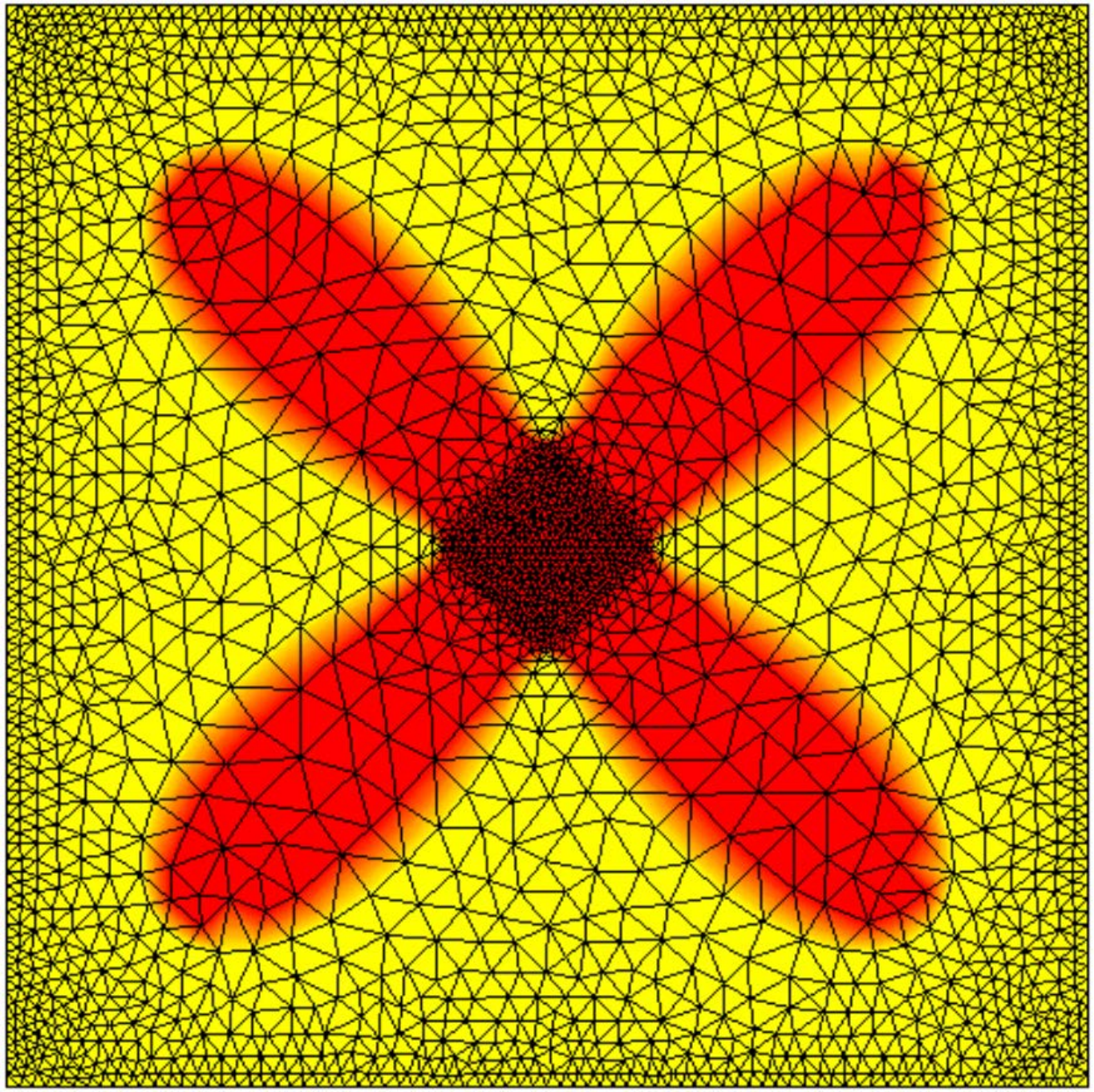} \\
  (d)
\end{minipage}
\caption{Four kinds of composite Kirchhoff plates: (a) type-I PUC $\mathbf{Y}$; (b) type-II PUC $\mathbf{Y}$; (c) type-III PUC $\mathbf{Y}$; (d) type-IV PUC $\mathbf{Y}$.}\label{f2}
\end{figure}

Based on the conclusion in subsection 5.1, the FOMS method based on Morley finite element is employed in this subsection. After numerical simulation, Table 3 shows the computational accuracy of the proposed FOMS method based on Morley finite element for composite Kirchhoff plates with four kinds of microscopic PUCs.
\begin{table}[h]{\caption{Comparison of computational accuracy.}\label{t2}}
\centering
\begin{tabular}{ccccccccc}
\hline
FOMS-Morley & $\frac{||e_0||_{L^2}}{||\omega_{\text{D}}^\epsilon||_{L^2}}$ & $\frac{||e_2||_{L^2}}{||\omega_{\text{D}}^\epsilon||_{L^2}}$ & $\frac{||e_3||_{L^2}}{||\omega_{\text{D}}^\epsilon||_{L^2}}$ & $\frac{||e_4||_{L^2}}{||\omega_{\text{D}}^\epsilon||_{L^2}}$ & $\frac{|e_0|_{H^1}}{|\omega_{\text{D}}^\epsilon|_{H^1}}$ & $\frac{|e_2|_{H^1}}{|\omega_{\text{D}}^\epsilon|_{H^1}}$ & $\frac{|e_3|_{H^1}}{|\omega_{\text{D}}^\epsilon|_{H^1}}$ & $\frac{|e_4|_{H^1}}{|\omega_{\text{D}}^\epsilon|_{H^1}}$\\
\hline
Type-I/Percentage \% & 9.937 & 9.893 & 9.897 & 2.733 & 87.058 & 86.856 & 86.829 & 21.630\\
\hline
Type-II/Percentage \% & 10.864 & 10.690 & 10.693 & 3.495 & 85.233 & 85.009 & 84.985 & 19.425\\
\hline
Type-III/Percentage \% & 19.570 & 19.442 & 19.445 & 5.693 & 93.222 & 93.090 & 93.069 & 22.923\\
\hline
Type-IV/Percentage \% & 1.495 & 1.031 & 1.015 & 0.345 & 29.887 & 26.199 & 25.710 & 8.586\\
\hline
\end{tabular}
\end{table}

Furthermore, Figs.~10-13 display the numerical results for solutions $\omega^{(0)}$, $\omega^{(2\epsilon)}$, $\omega^{(3\epsilon)}$, $\omega^{(4\epsilon)}$ and $\omega^\epsilon_{\text{D}}$ based on Morley finite element of composite Kirchhoff plates with four kinds of microscopic PUCs, respectively.
\begin{figure}[!htb]
\centering
\begin{minipage}[c]{0.32\textwidth}
  \centering
  \includegraphics[width=45mm]{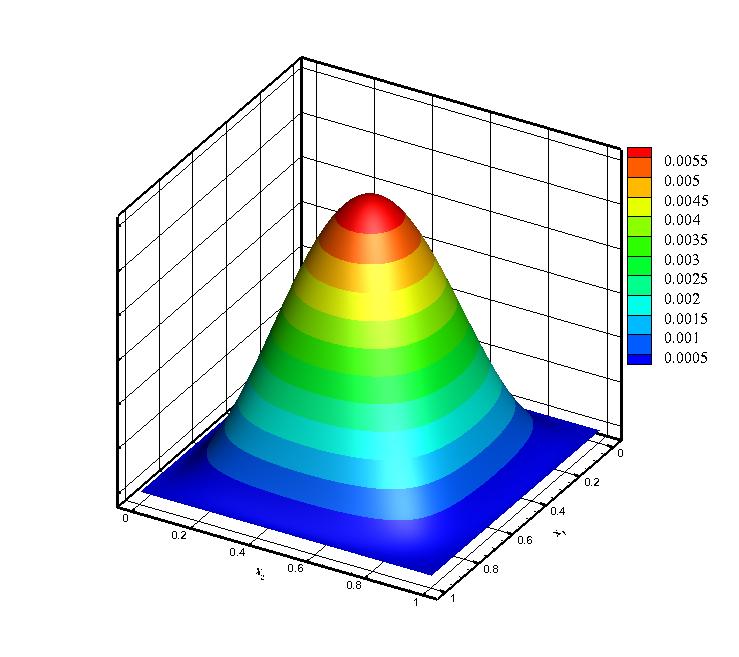}\\
  (a)
\end{minipage}
\begin{minipage}[c]{0.32\textwidth}
  \centering
  \includegraphics[width=45mm]{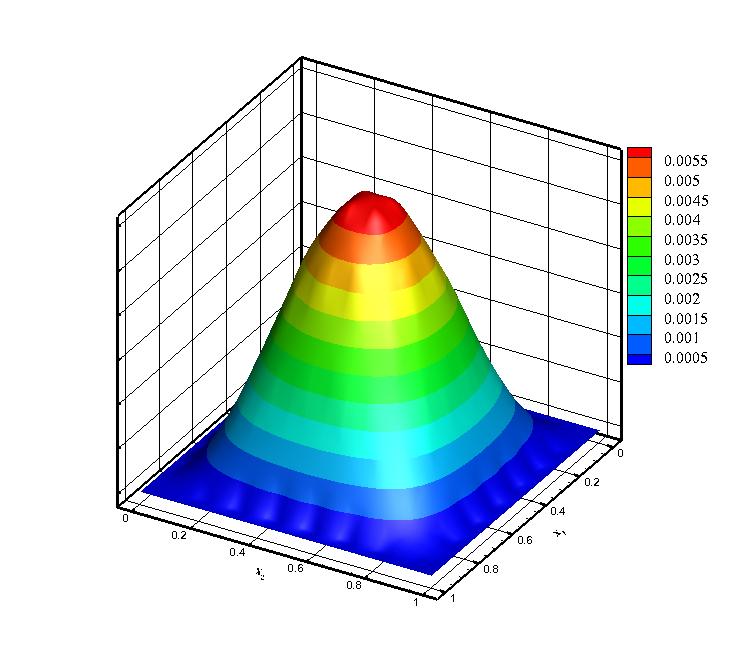}\\
  (b)
\end{minipage}
\begin{minipage}[c]{0.32\textwidth}
  \centering
  \includegraphics[width=45mm]{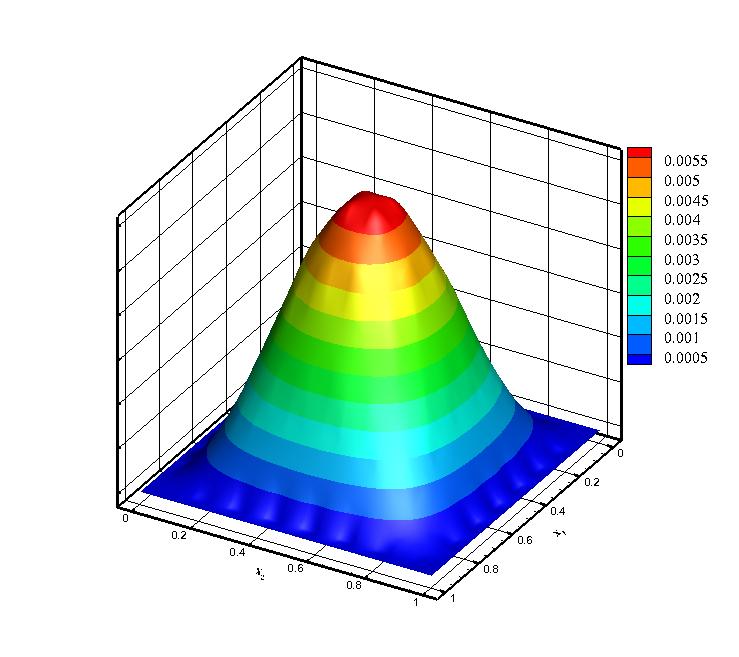}\\
  (c)
\end{minipage}
\begin{minipage}[c]{0.32\textwidth}
  \centering
  \includegraphics[width=45mm]{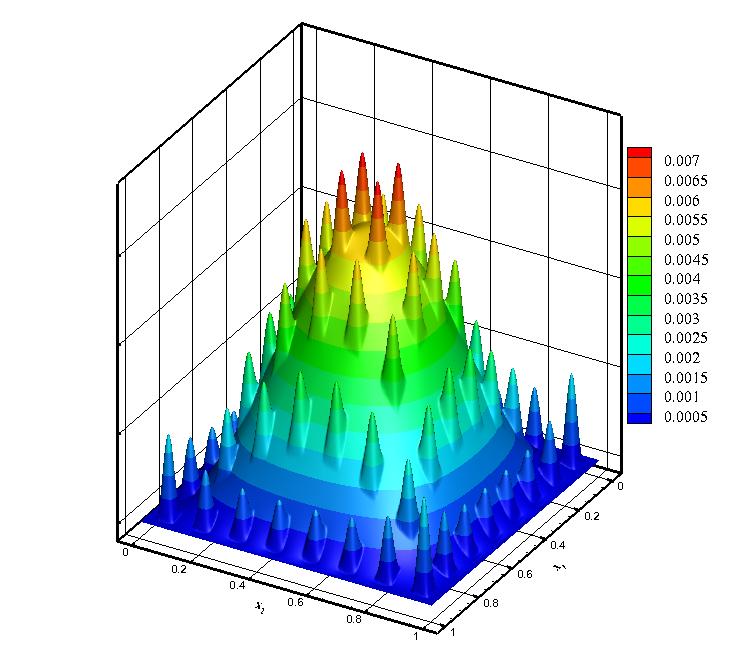}\\
  (d)
\end{minipage}
\begin{minipage}[c]{0.32\textwidth}
  \centering
  \includegraphics[width=45mm]{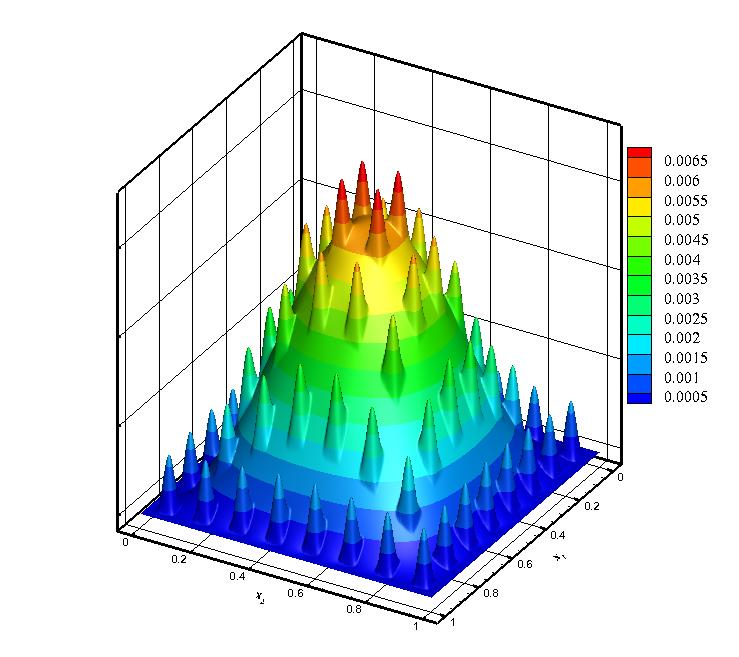}\\
  (e)
\end{minipage}
\caption{The transverse displacement of composite Kirchhoff plate with type-I PUC computed by Morley finite element: (a) $\omega^{(0)}$; (b) $\omega^{(2\epsilon)}$; (c) $\omega^{(3\epsilon)}$; (d) $\omega^{(4\epsilon)}$; (e) $\omega^\epsilon_{\text{D}}$.}\label{f8}
\end{figure}

\begin{figure}[!htb]
\centering
\begin{minipage}[c]{0.32\textwidth}
  \centering
  \includegraphics[width=45mm]{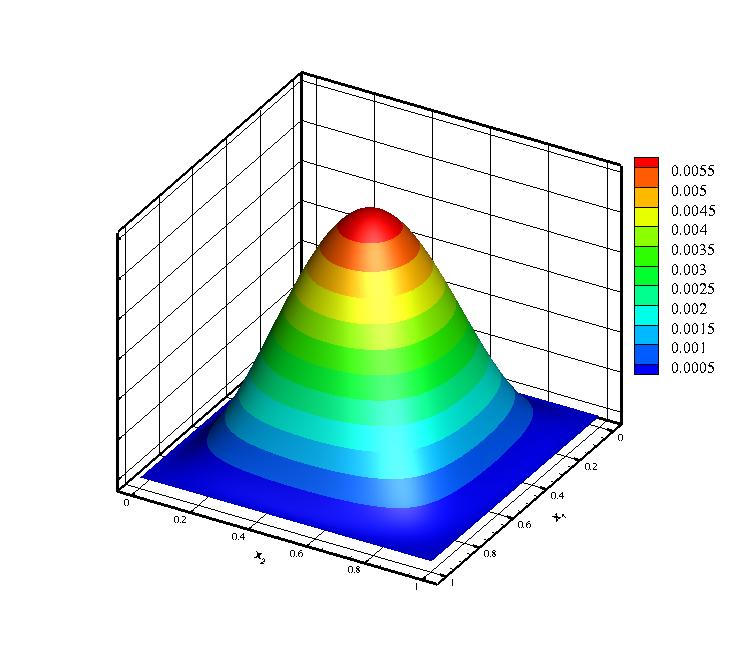}\\
  (a)
\end{minipage}
\begin{minipage}[c]{0.32\textwidth}
  \centering
  \includegraphics[width=45mm]{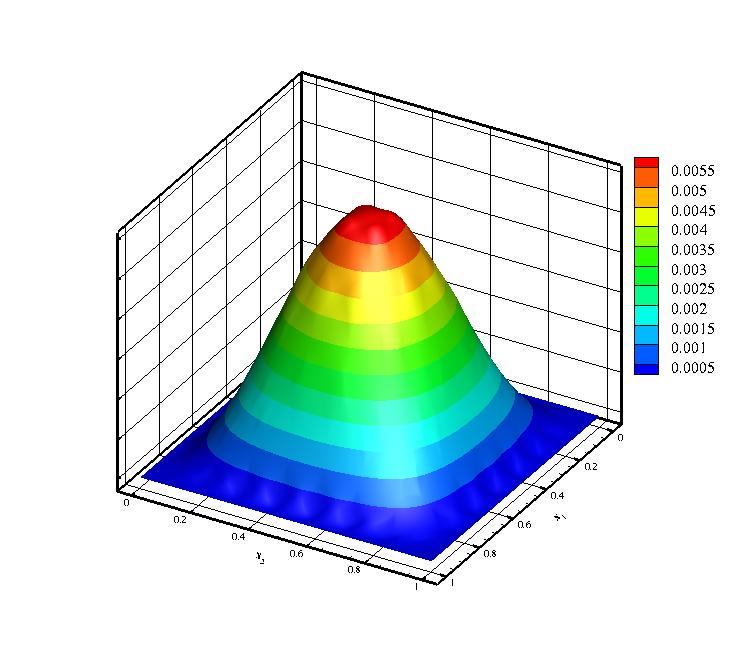}\\
  (b)
\end{minipage}
\begin{minipage}[c]{0.32\textwidth}
  \centering
  \includegraphics[width=45mm]{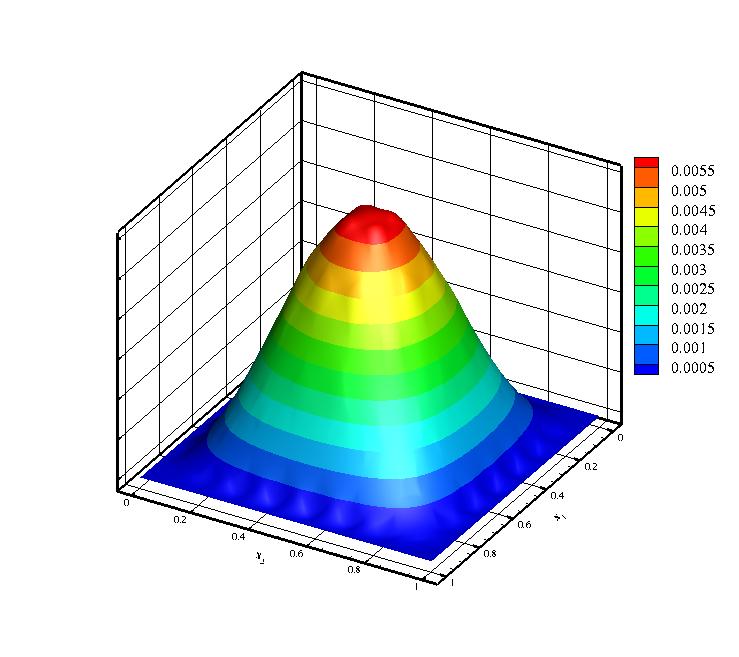}\\
  (c)
\end{minipage}
\begin{minipage}[c]{0.32\textwidth}
  \centering
  \includegraphics[width=45mm]{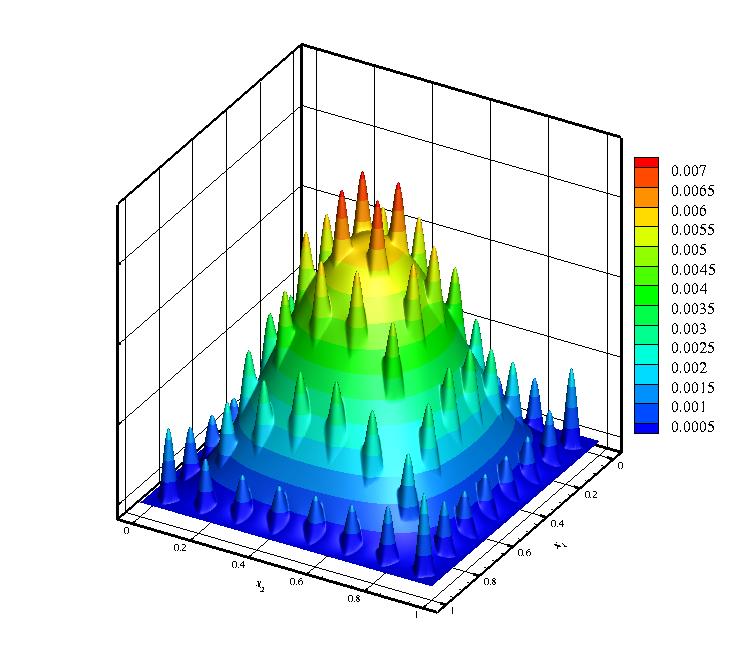}\\
  (d)
\end{minipage}
\begin{minipage}[c]{0.32\textwidth}
  \centering
  \includegraphics[width=45mm]{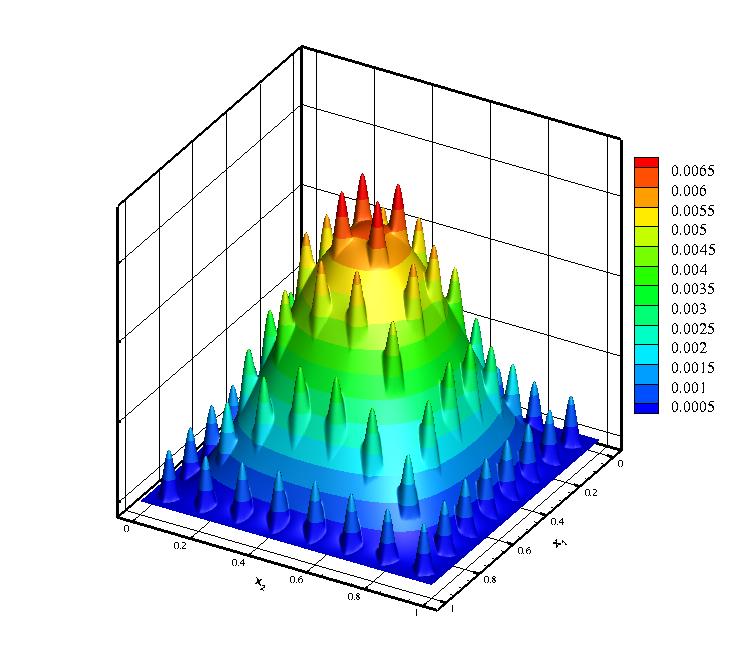}\\
  (e)
\end{minipage}
\caption{The transverse displacement of composite Kirchhoff plate with type-II PUC computed by Morley finite element: (a) $\omega^{(0)}$; (b) $\omega^{(2\epsilon)}$; (c) $\omega^{(3\epsilon)}$; (d) $\omega^{(4\epsilon)}$; (e) $\omega^\epsilon_{\text{D}}$.}\label{f8}
\end{figure}

\begin{figure}[!htb]
\centering
\begin{minipage}[c]{0.32\textwidth}
  \centering
  \includegraphics[width=45mm]{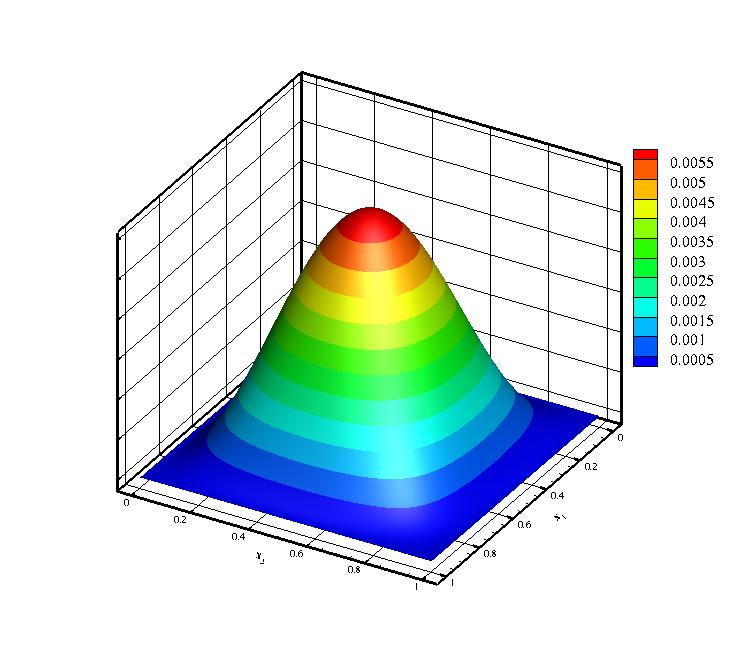}\\
  (a)
\end{minipage}
\begin{minipage}[c]{0.32\textwidth}
  \centering
  \includegraphics[width=45mm]{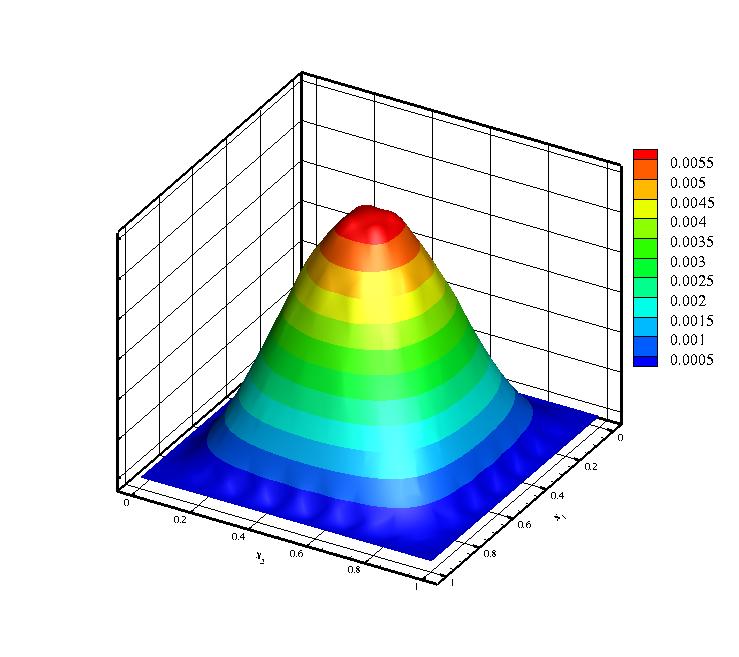}\\
  (b)
\end{minipage}
\begin{minipage}[c]{0.32\textwidth}
  \centering
  \includegraphics[width=45mm]{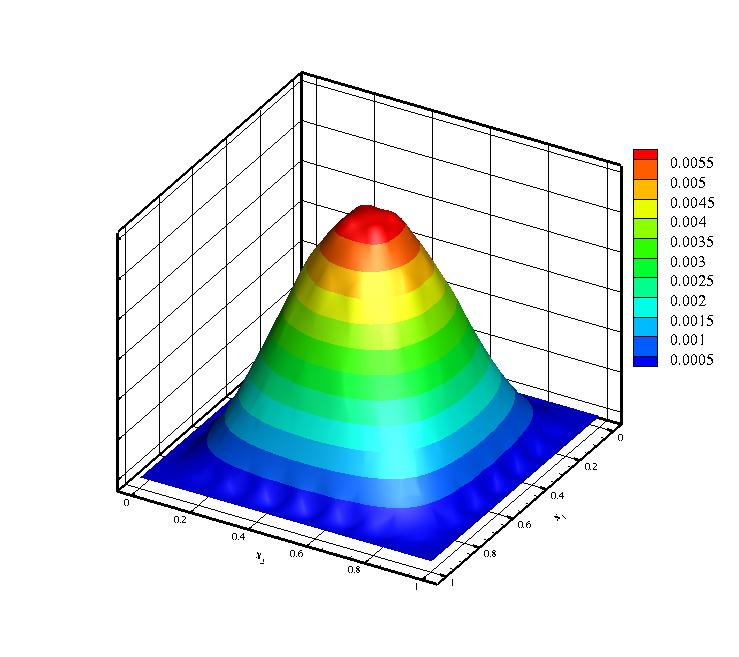}\\
  (c)
\end{minipage}
\begin{minipage}[c]{0.32\textwidth}
  \centering
  \includegraphics[width=45mm]{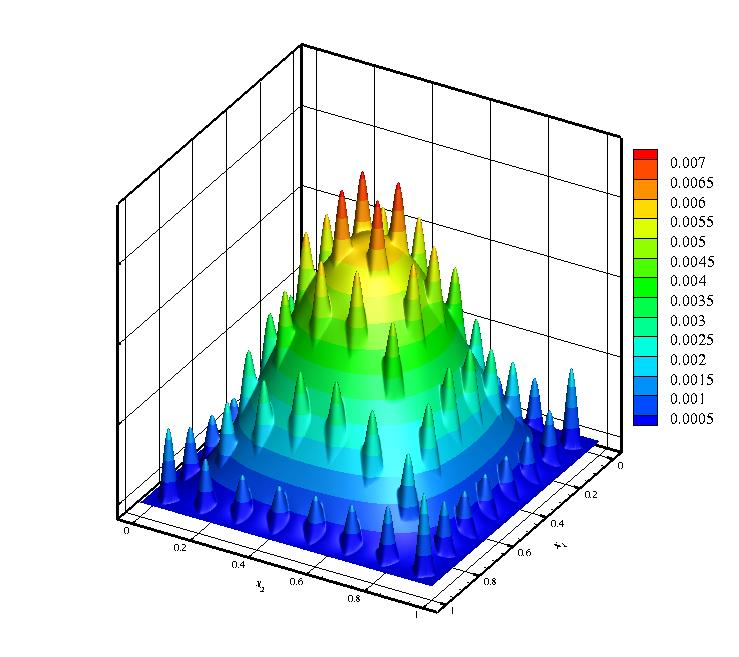}\\
  (d)
\end{minipage}
\begin{minipage}[c]{0.32\textwidth}
  \centering
  \includegraphics[width=45mm]{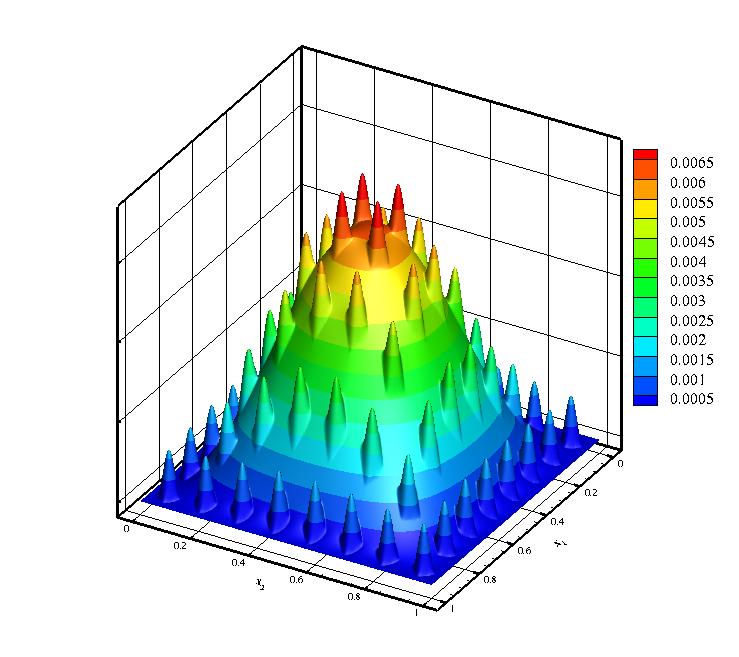}\\
  (e)
\end{minipage}
\caption{The transverse displacement of composite Kirchhoff plate with type-III PUC computed by Morley finite element: (a) $\omega^{(0)}$; (b) $\omega^{(2\epsilon)}$; (c) $\omega^{(3\epsilon)}$; (d) $\omega^{(4\epsilon)}$; (e) $\omega^\epsilon_{\text{D}}$.}\label{f8}
\end{figure}

\begin{figure}[!htb]
\centering
\begin{minipage}[c]{0.32\textwidth}
  \centering
  \includegraphics[width=45mm]{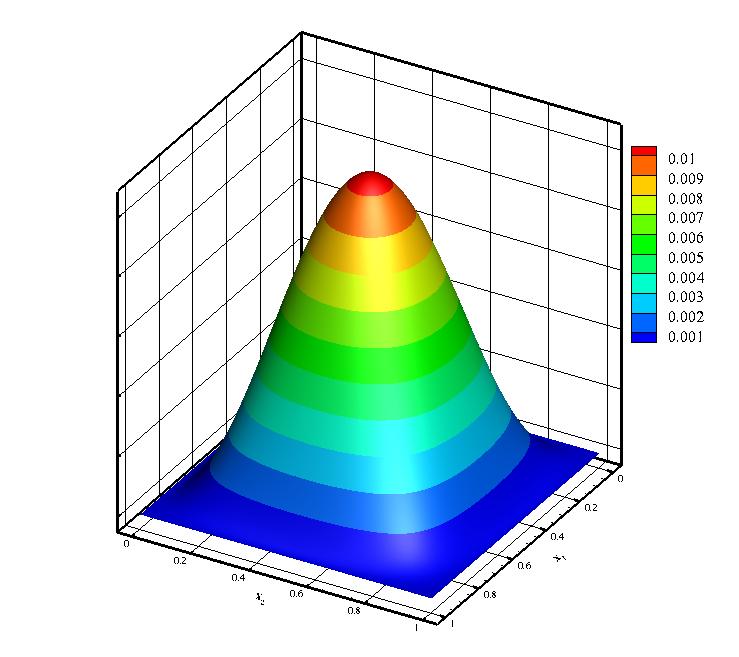}\\
  (a)
\end{minipage}
\begin{minipage}[c]{0.32\textwidth}
  \centering
  \includegraphics[width=45mm]{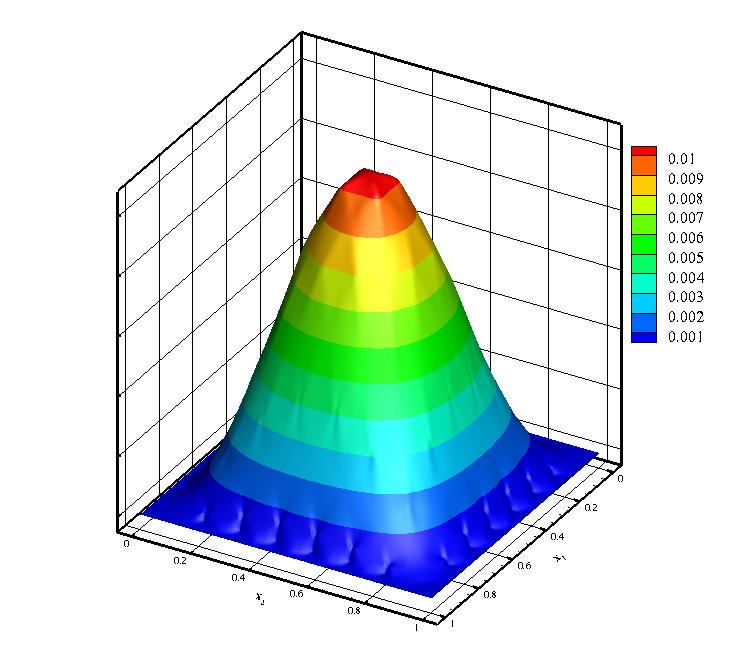}\\
  (b)
\end{minipage}
\begin{minipage}[c]{0.32\textwidth}
  \centering
  \includegraphics[width=45mm]{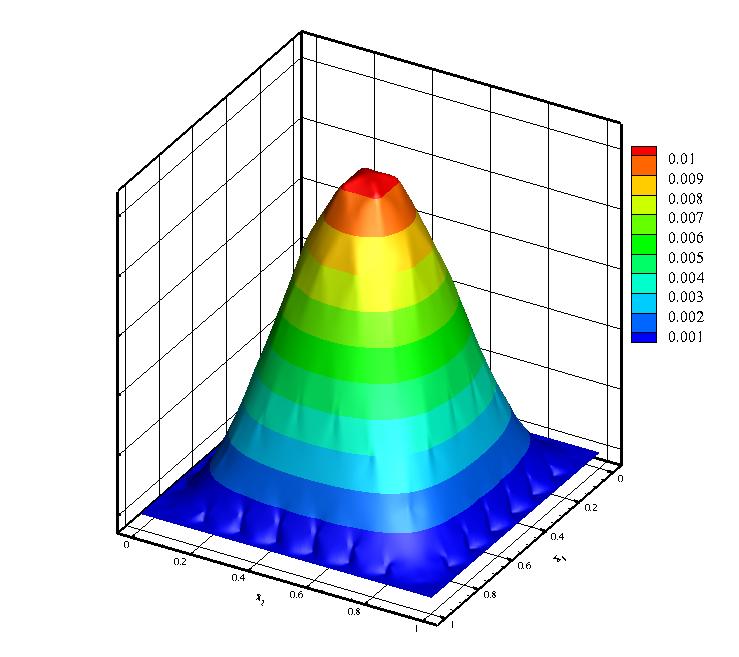}\\
  (c)
\end{minipage}
\begin{minipage}[c]{0.32\textwidth}
  \centering
  \includegraphics[width=45mm]{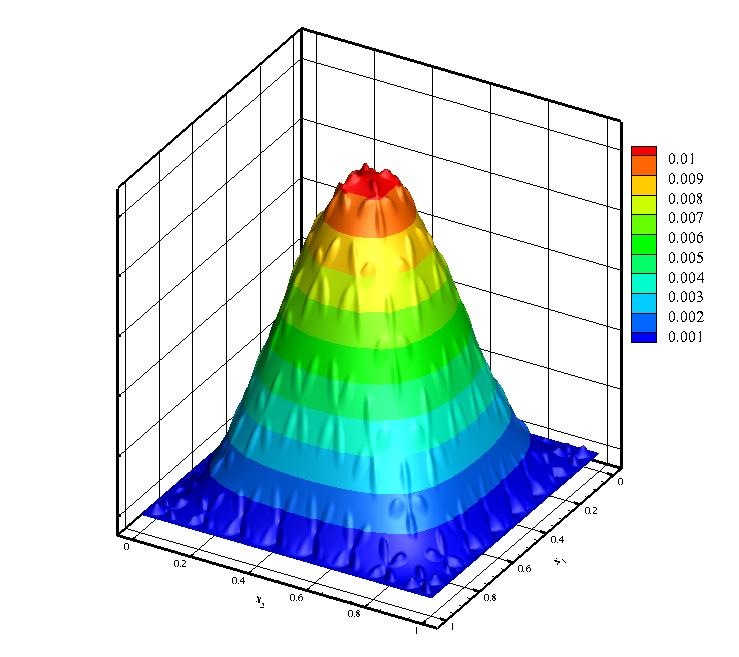}\\
  (d)
\end{minipage}
\begin{minipage}[c]{0.32\textwidth}
  \centering
  \includegraphics[width=45mm]{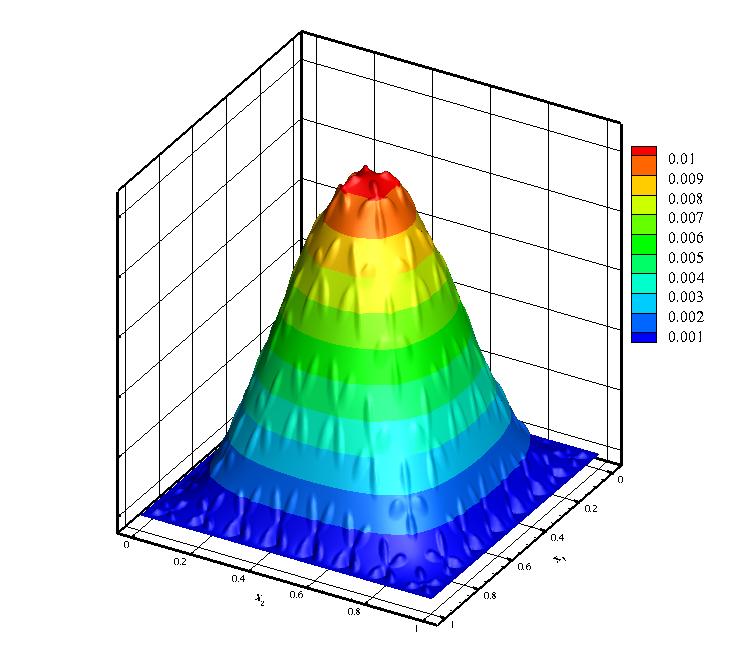}\\
  (e)
\end{minipage}
\caption{The transverse displacement of composite Kirchhoff plate with type-IV PUC computed by Morley finite element: (a) $\omega^{(0)}$; (b) $\omega^{(2\epsilon)}$; (c) $\omega^{(3\epsilon)}$; (d) $\omega^{(4\epsilon)}$; (e) $\omega^\epsilon_{\text{D}}$.}\label{f8}
\end{figure}

As shown in Table 3, it can be clearly found that the FOMS approach exhibits the optimal numerical accuracy for highly oscillating behaviors of four kinds of composite Kirchhoff plates, compared with macroscopic homogenized, second-order multi-scale and third-order multi-scale approaches. In addition, it can be obviously seen from Figs.~10-13 that the FOMS solutions for four kinds of composite Kirchhoff plates agree well with those of direct numerical simulation on highly refined meshes. In stark contrast, the macroscopic homogenized solution can merely reflect the macroscopic bending behavior of the composite plates, while the second- and third-order multi-scale solutions can capture the inadequate local information at the micro-scale. 

We should highlight that another potential advantage of the introduced FOMS methodology is its applicability to the effective computation of large-scale composite plates with a vast number of unit cells. In order to validate this advanced strength, the proposed FOMS methodology is utilized to compute the bending behaviors of the aforementioned four kinds of composite Kirchhoff plates with the characteristic periodic parameter $\epsilon=1/16$, as shown in Fig.~14.
\begin{figure}[!htb]
\centering
\begin{minipage}[c]{0.35\textwidth}
  \centering
  \includegraphics[width=1.0\linewidth,totalheight=1.7in]{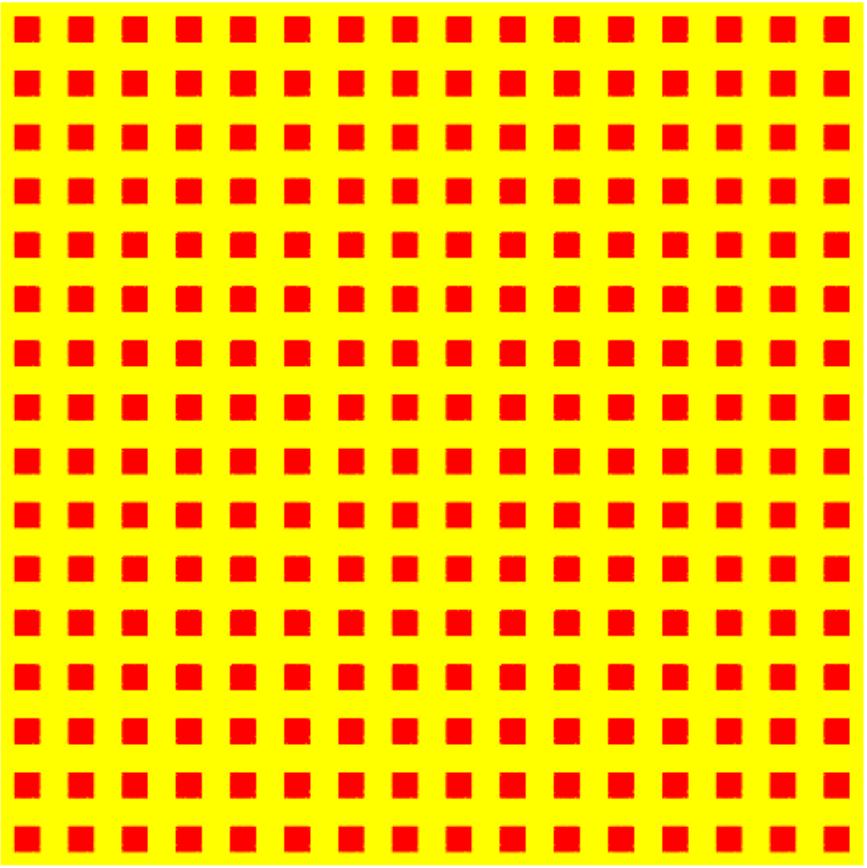} \\
  (a)
\end{minipage}
\begin{minipage}[c]{0.35\textwidth}
  \centering
  \includegraphics[width=1.0\linewidth,totalheight=1.7in]{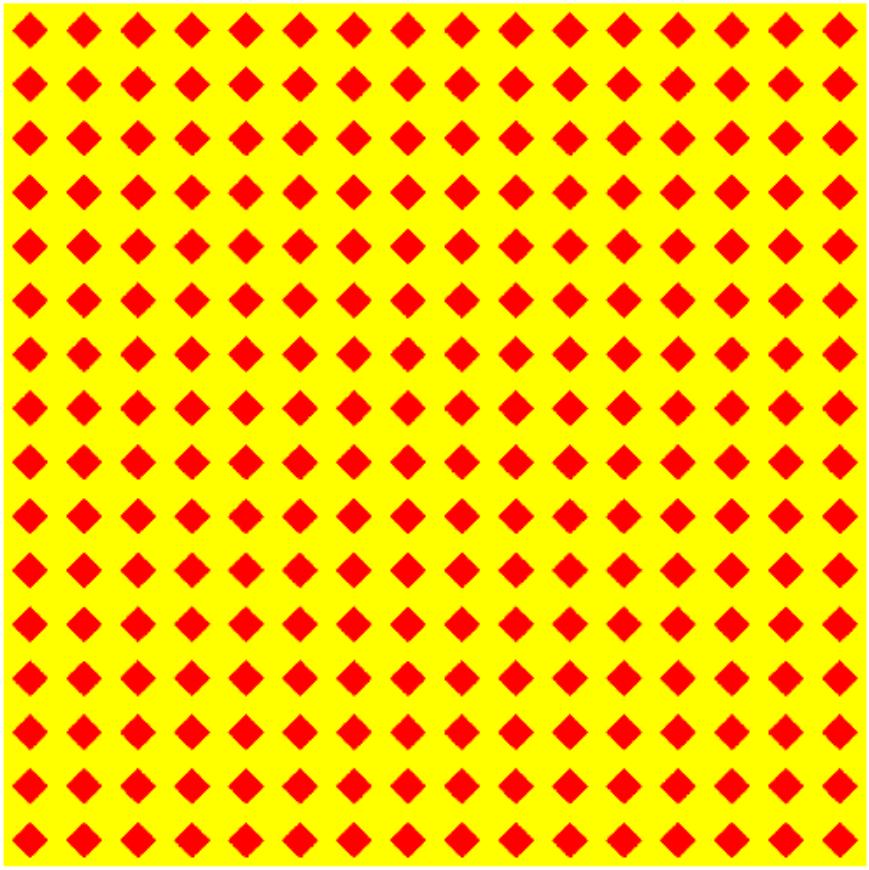} \\
  (b)
\end{minipage}
\begin{minipage}[c]{0.35\textwidth}
  \centering
  \includegraphics[width=1.0\linewidth,totalheight=1.7in]{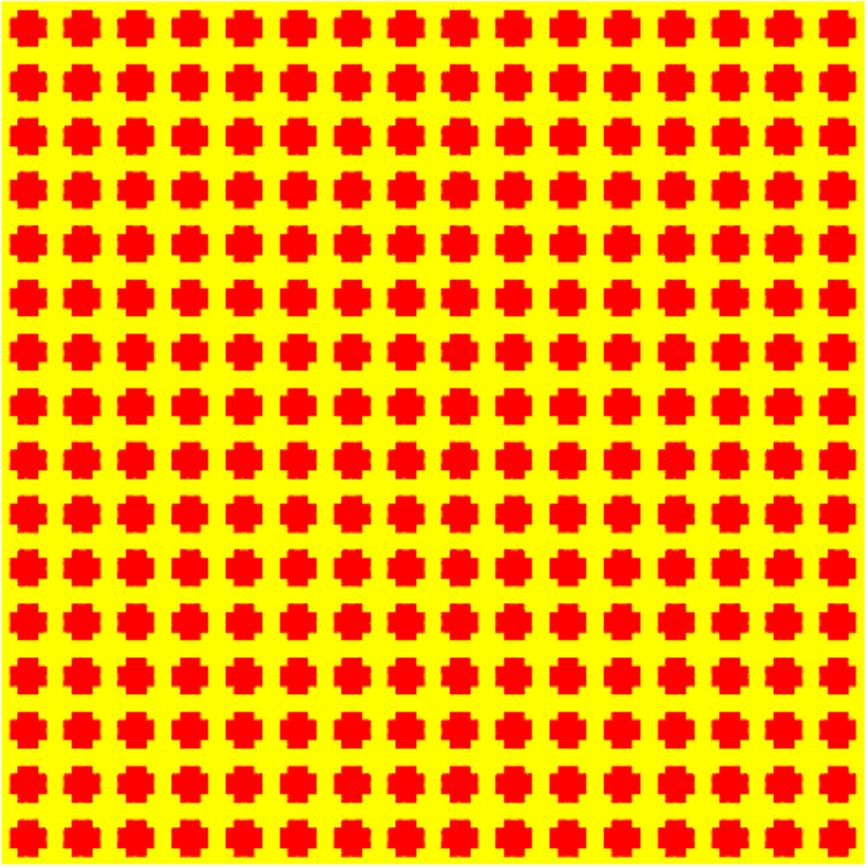} \\
  (c)
\end{minipage}
\begin{minipage}[c]{0.35\textwidth}
  \centering
  \includegraphics[width=1.0\linewidth,totalheight=1.7in]{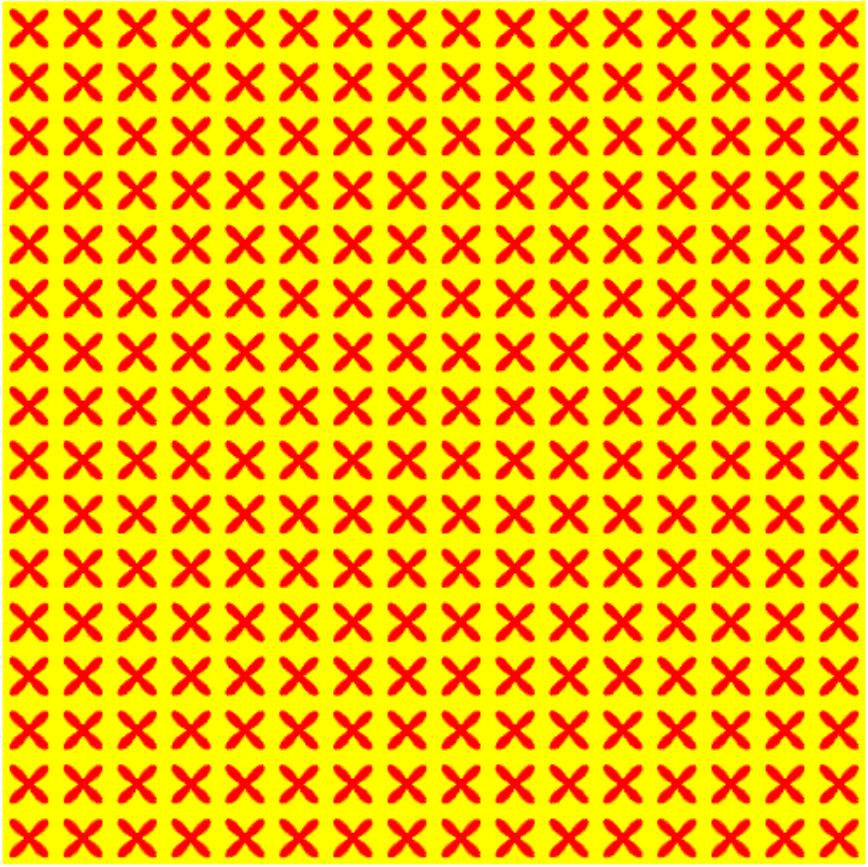} \\
  (d)
\end{minipage}
\caption{Four kinds of composite Kirchhoff plates with $\epsilon=1/16$: (a) type-I PUC $\mathbf{Y}$; (b) type-II PUC $\mathbf{Y}$; (c) type-III PUC $\mathbf{Y}$; (d) type-IV PUC $\mathbf{Y}$.}\label{f2}
\end{figure}

Moreover, the transverse displacement of four kinds of composite Kirchhoff plates with a large amount of PUCs is shown in Fig.~15. At this time, the Freefem++ software is unable to access the high-resolution reference solutions due to a highly refined mesh. However, the proposed FOMS method enables accurate and efficient computation of these fourth-order multi-scale problems with highly spatial heterogeneities.
\begin{figure}[!htb]
\centering
\begin{minipage}[c]{0.35\textwidth}
  \centering
  \includegraphics[width=1.0\linewidth,totalheight=1.7in]{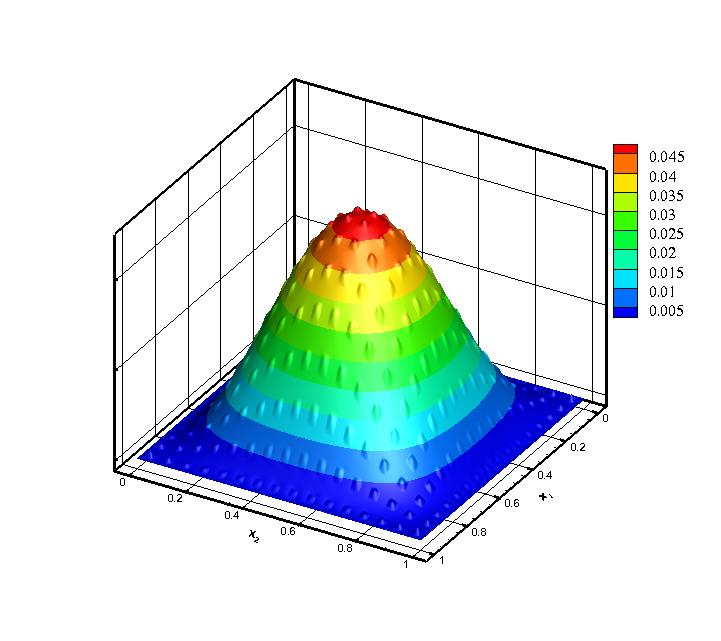} \\
  (a)
\end{minipage}
\begin{minipage}[c]{0.35\textwidth}
  \centering
  \includegraphics[width=1.0\linewidth,totalheight=1.7in]{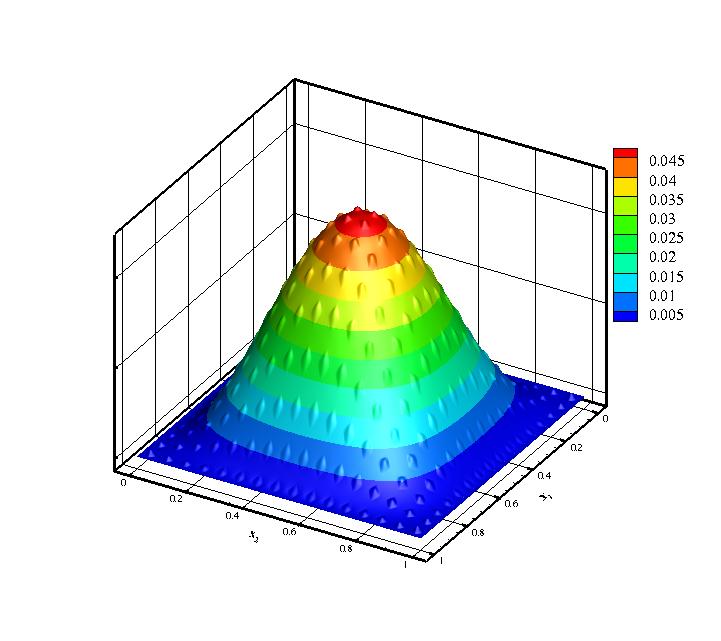} \\
  (b)
\end{minipage}
\begin{minipage}[c]{0.35\textwidth}
  \centering
  \includegraphics[width=1.0\linewidth,totalheight=1.85in]{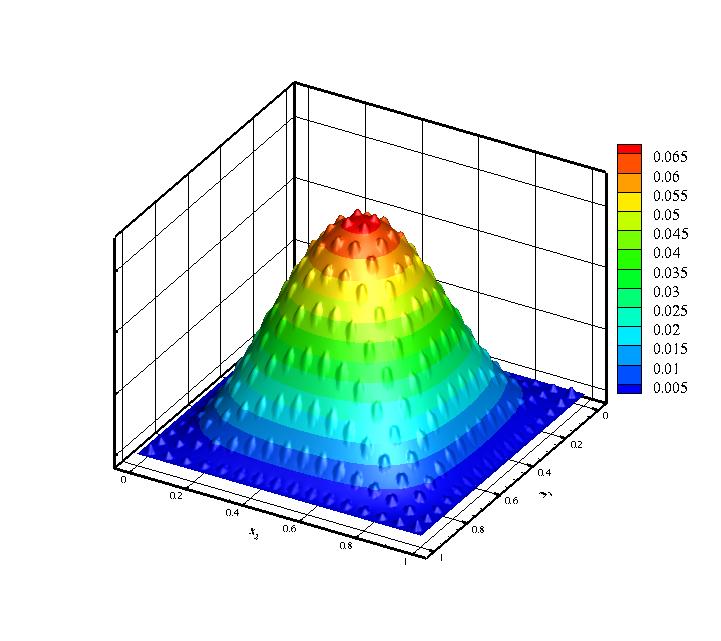} \\
  (c)
\end{minipage}
\begin{minipage}[c]{0.35\textwidth}
  \centering
  \includegraphics[width=1.0\linewidth,totalheight=1.7in]{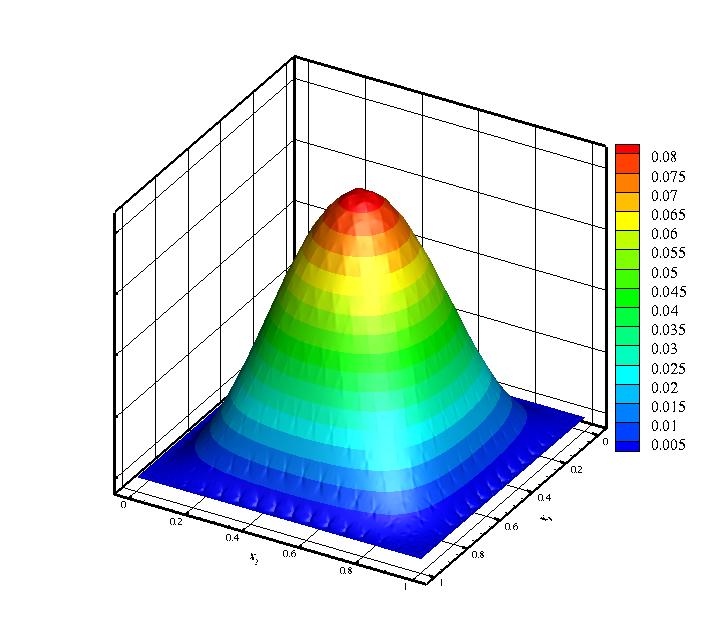} \\
  (d)
\end{minipage}
\caption{The transverse displacement of four kinds of composite Kirchhoff plates with $\epsilon=1/16$ computed by Morley finite element: (a) type-I PUC $\mathbf{Y}$; (b) type-II PUC $\mathbf{Y}$; (c) type-III PUC $\mathbf{Y}$; (d) type-IV PUC $\mathbf{Y}$.}\label{f2}
\end{figure}

In real-world engineering applications, it is impractical to obtain reference FEM solutions for large-scale composite plates due to prohibitive computational overhead. However, the FOMS method can effectively simulate large-scale composite plates with minimal computational resource consumption. Also, the proposed FOMS approach is robust for simulating the composite plates with massive PUCs without increasing the computational overhead, which is of great engineering values.
\section{\noindent \textbf{Conclusions and outlooks}}
Owing to the key properties of being lightweight, high-strength, and enabling large structural spans, composite plates have abroad application in the aerospace, construction, and machinery industries, etc. However, analytical approaches are prohibitively difficult to obtain an exact solution of fourth-order composite Kirchhoff plate models with rapidly oscillating and highly discontinuous coefficients. Moreover, direct simulation by classical numerical methods involves a substantial computational cost. Hence, the development of accurate and efficient multi-scale method is essential for effectively investigating the mechanical behavior of composite plate structures, offering significant scientific insights and practical applications.

This study proposes a novel fourth-order multi-scale method (fourth-order multi-scale computational model and numerical algorithm), which paves the way for simulating the composite Kirchhoff plates both accurately and efficiently. The innovations of this work are threefold: high-accuracy FOMS computational model with novel fourth-order correctors, local and global error analyses for validating the local balance preserving and proving the explicit convergence, efficient numerical algorithm with reduced computational resource. Furthermore, numerical experiments clearly demonstrated that the proposed FOMS methodology with local balance preserving can accurately capture the highly oscillating behavior of composite plates and provide high-fidelity numerical solutions for engineering applications. Also, the fourth-order multi-scale method owns excellent agreement between its theoretical framework and numerical results. Furthermore, numerical experiments also illustrate the proposed FOMS method can strongly economize the computer resource compared with finite element method. These computational advantages are of great application value for large-scale engineering simulation.

In the future, the introduced FOMS method will extend to simulate the dynamic vibration problem of composite Kirchhoff plates. Additionally, the FOMS approach will focus on the extension to composite plate considering three-scale, nonlinear and nonlocal mechanical behaviors.

\section*{Acknowledgments}
This research was supported by the National Natural Science Foundation of China (No.\hspace{1mm}12471387), Xidian University Specially Funded Project for Interdisciplinary Exploration (No.\hspace{1mm}TZJH2024008), Fundamental Research Funds for the Central Universities (No.\hspace{1mm}QTZX25082), and also supported by the High-Performance Computing Platform of Xidian University.

\bibliographystyle{model1-num-names}
\bibliography{paper}

\begin{thebibliography}{40}
\expandafter\ifx\csname natexlab\endcsname\relax\def\natexlab#1{#1}\fi
\providecommand{\url}[1]{\texttt{#1}}
\providecommand{\href}[2]{#2}
\providecommand{\path}[1]{#1}
\providecommand{\DOIprefix}{doi:}
\providecommand{\ArXivprefix}{arXiv:}
\providecommand{\URLprefix}{URL: }
\providecommand{\Pubmedprefix}{pmid:}
\providecommand{\doi}[1]{\href{http://dx.doi.org/#1}{\path{#1}}}
\providecommand{\Pubmed}[1]{\href{pmid:#1}{\path{#1}}}
\providecommand{\bibinfo}[2]{#2}
\ifx\xfnm\relax \def\xfnm[#1]{\unskip,\space#1}\fi
\bibitem[{Lewinski and Telega(2000)}]{R1}
\bibinfo{author}{T.~Lewinski}, \bibinfo{author}{J.~J. Telega},
  \bibinfo{title}{Plates, laminates and shells: asymptotic analysis and
  homogenization}, volume~\bibinfo{volume}{52}, \bibinfo{publisher}{World
  Scientific}, \bibinfo{year}{2000}.
\bibitem[{Cioranescu and Donato(1999)}]{R2}
\bibinfo{author}{D.~Cioranescu}, \bibinfo{author}{P.~Donato},
  \bibinfo{title}{An introduction to homogenization},
  \bibinfo{publisher}{Oxford university press}, \bibinfo{year}{1999}.
\bibitem[{Ole{\"\i}nik et~al.(1992)Ole{\"\i}nik, Shamaev, and Yosifian}]{R3}
\bibinfo{author}{O.~A. Ole{\"\i}nik}, \bibinfo{author}{A.~Shamaev},
  \bibinfo{author}{G.~Yosifian}, \bibinfo{title}{Mathematical problems in
  elasticity and homogenization}, volume~\bibinfo{volume}{26},
  \bibinfo{publisher}{Elsevier}, \bibinfo{year}{1992}.
\bibitem[{Hou et~al.(1999)Hou, Wu, and Cai}]{R4}
\bibinfo{author}{T.~Hou}, \bibinfo{author}{X.-H. Wu}, \bibinfo{author}{Z.~Cai},
\newblock \bibinfo{title}{Convergence of a multiscale finite element method for
  elliptic problems with rapidly oscillating coefficients},
\newblock \bibinfo{journal}{Mathematics of Computation} \bibinfo{volume}{68}
  (\bibinfo{year}{1999}) \bibinfo{pages}{913--943}.
\bibitem[{Efendiev and Hou(2009)}]{R5}
\bibinfo{author}{Y.~Efendiev}, \bibinfo{author}{T.~Y. Hou},
  \bibinfo{title}{Multiscale finite element methods: theory and applications},
  volume~\bibinfo{volume}{4}, \bibinfo{publisher}{Springer Science \& Business
  Media}, \bibinfo{year}{2009}.
\bibitem[{Abdulle et~al.(2012)Abdulle, Weinan, Engquist, and
  Vanden-Eijnden}]{R6}
\bibinfo{author}{A.~Abdulle}, \bibinfo{author}{E.~Weinan},
  \bibinfo{author}{B.~Engquist}, \bibinfo{author}{E.~Vanden-Eijnden},
\newblock \bibinfo{title}{The heterogeneous multiscale method},
\newblock \bibinfo{journal}{Acta Numerica} \bibinfo{volume}{21}
  (\bibinfo{year}{2012}) \bibinfo{pages}{1--87}.
\bibitem[{Hughes(1995)}]{R7}
\bibinfo{author}{T.~J. Hughes},
\newblock \bibinfo{title}{Multiscale phenomena: Green's functions, the
  dirichlet-to-neumann formulation, subgrid scale models, bubbles and the
  origins of stabilized methods},
\newblock \bibinfo{journal}{Computer Methods in Applied Mechanics and
  Engineering} \bibinfo{volume}{127} (\bibinfo{year}{1995})
  \bibinfo{pages}{387--401}.
\bibitem[{Xing et~al.(2010)Xing, Yang, and Wang}]{R8}
\bibinfo{author}{Y.~Xing}, \bibinfo{author}{Y.~Yang},
  \bibinfo{author}{X.~Wang},
\newblock \bibinfo{title}{A multiscale eigenelement method and its application
  to periodical composite structures},
\newblock \bibinfo{journal}{Composite Structures} \bibinfo{volume}{92}
  (\bibinfo{year}{2010}) \bibinfo{pages}{2265--2275}.
\bibitem[{Henning and M{\aa}lqvist(2014)}]{R9}
\bibinfo{author}{P.~Henning}, \bibinfo{author}{A.~M{\aa}lqvist},
\newblock \bibinfo{title}{Localized orthogonal decomposition techniques for
  boundary value problems},
\newblock \bibinfo{journal}{SIAM Journal on Scientific Computing}
  \bibinfo{volume}{36} (\bibinfo{year}{2014}) \bibinfo{pages}{A1609--A1634}.
\bibitem[{Efendiev et~al.(2013)Efendiev, Galvis, and Hou}]{R10}
\bibinfo{author}{Y.~Efendiev}, \bibinfo{author}{J.~Galvis},
  \bibinfo{author}{T.~Y. Hou},
\newblock \bibinfo{title}{Generalized multiscale finite element methods
  (gmsfem)},
\newblock \bibinfo{journal}{Journal of computational physics}
  \bibinfo{volume}{251} (\bibinfo{year}{2013}) \bibinfo{pages}{116--135}.
\bibitem[{Chung et~al.(2018)Chung, Efendiev, and Leung}]{R11}
\bibinfo{author}{E.~T. Chung}, \bibinfo{author}{Y.~Efendiev},
  \bibinfo{author}{W.~T. Leung},
\newblock \bibinfo{title}{Constraint energy minimizing generalized multiscale
  finite element method},
\newblock \bibinfo{journal}{Computer Methods in Applied Mechanics and
  Engineering} \bibinfo{volume}{339} (\bibinfo{year}{2018})
  \bibinfo{pages}{298--319}.
\bibitem[{Dong et~al.(2019)Dong, Cui, Nie, Yang, Ma, and Cheng}]{R12}
\bibinfo{author}{H.~Dong}, \bibinfo{author}{J.~Cui}, \bibinfo{author}{Y.~Nie},
  \bibinfo{author}{Z.~Yang}, \bibinfo{author}{Q.~Ma},
  \bibinfo{author}{X.~Cheng},
\newblock \bibinfo{title}{Multiscale computational method for transient heat
  conduction problems of periodic porous materials with diverse periodic
  configurations in different subdomains},
\newblock \bibinfo{journal}{Applied Numerical Mathematics}
  \bibinfo{volume}{136} (\bibinfo{year}{2019}) \bibinfo{pages}{215--234}.
\bibitem[{Dong et~al.(2021)Dong, Cui, Nie, Jin, Guan, and Yang}]{R13}
\bibinfo{author}{H.~Dong}, \bibinfo{author}{J.~Cui}, \bibinfo{author}{Y.~Nie},
  \bibinfo{author}{K.~Jin}, \bibinfo{author}{X.~Guan},
  \bibinfo{author}{Z.~Yang},
\newblock \bibinfo{title}{High-order three-scale computational method for
  elastic behavior analysis and strength prediction of axisymmetric composite
  structures with multiple spatial scales},
\newblock \bibinfo{journal}{Mathematics and Mechanics of Solids}
  \bibinfo{volume}{26} (\bibinfo{year}{2021}) \bibinfo{pages}{905--936}.
\bibitem[{Dong et~al.(2025)Dong, Guan, and Nie}]{R14}
\bibinfo{author}{H.~Dong}, \bibinfo{author}{X.~Guan}, \bibinfo{author}{Y.~Nie},
\newblock \bibinfo{title}{Multiscale method and convergence analysis for
  coupled nonlinear thermomechanical problems in heterogeneous shells},
\newblock \bibinfo{journal}{SIAM Journal on Scientific Computing}
  \bibinfo{volume}{47} (\bibinfo{year}{2025}) \bibinfo{pages}{B190--B219}.
\bibitem[{Cao et~al.(2002)Cao, Cui, and Zhu}]{R15}
\bibinfo{author}{L.-Q. Cao}, \bibinfo{author}{J.-Z. Cui},
  \bibinfo{author}{D.-C. Zhu},
\newblock \bibinfo{title}{Multiscale asymptotic analysis and numerical
  simulation for the second order helmholtz equations with rapidly oscillating
  coefficients over general convex domains},
\newblock \bibinfo{journal}{SIAM Journal on Numerical Analysis}
  \bibinfo{volume}{40} (\bibinfo{year}{2002}) \bibinfo{pages}{543--577}.
\bibitem[{Cao and Cui(2004)}]{R16}
\bibinfo{author}{L.-Q. Cao}, \bibinfo{author}{J.-Z. Cui},
\newblock \bibinfo{title}{Asymptotic expansions and numerical algorithms of
  eigenvalues and eigenfunctions of the dirichlet problem for second order
  elliptic equations in perforated domains},
\newblock \bibinfo{journal}{Numerische Mathematik} \bibinfo{volume}{96}
  (\bibinfo{year}{2004}) \bibinfo{pages}{525--581}.
\bibitem[{Pastukhova(2016)}]{R17}
\bibinfo{author}{S.~E. Pastukhova},
\newblock \bibinfo{title}{Estimates in homogenization of higher-order elliptic
  operators},
\newblock \bibinfo{journal}{Applicable Analysis} \bibinfo{volume}{95}
  (\bibinfo{year}{2016}) \bibinfo{pages}{1449--1466}.
\bibitem[{Pastukhova(2017)}]{R18}
\bibinfo{author}{S.~Pastukhova},
\newblock \bibinfo{title}{Operator error estimates for homogenization of fourth
  order elliptic equations},
\newblock \bibinfo{journal}{St. Petersburg Mathematical Journal}
  \bibinfo{volume}{28} (\bibinfo{year}{2017}) \bibinfo{pages}{273--289}.
\bibitem[{Suslina(2018{\natexlab{a}})}]{R19}
\bibinfo{author}{T.~Suslina},
\newblock \bibinfo{title}{Homogenization of the dirichlet problem for
  higher-order elliptic equations with periodic coefficients},
\newblock \bibinfo{journal}{St. Petersburg Mathematical Journal}
  \bibinfo{volume}{29} (\bibinfo{year}{2018}{\natexlab{a}})
  \bibinfo{pages}{325--362}.
\bibitem[{Suslina(2018{\natexlab{b}})}]{R20}
\bibinfo{author}{T.~Suslina},
\newblock \bibinfo{title}{Homogenization of the neumann problem for higher
  order elliptic equations with periodic coefficients},
\newblock \bibinfo{journal}{Complex Variables and Elliptic Equations}
  \bibinfo{volume}{63} (\bibinfo{year}{2018}{\natexlab{b}})
  \bibinfo{pages}{1185--1215}.
\bibitem[{Niu et~al.(2018)Niu, Shen, and Xu}]{R21}
\bibinfo{author}{W.~Niu}, \bibinfo{author}{Z.~Shen}, \bibinfo{author}{Y.~Xu},
\newblock \bibinfo{title}{Convergence rates and interior estimates in
  homogenization of higher order elliptic systems},
\newblock \bibinfo{journal}{Journal of Functional Analysis}
  \bibinfo{volume}{274} (\bibinfo{year}{2018}) \bibinfo{pages}{2356--2398}.
\bibitem[{Niu and Xu(2019)}]{R22}
\bibinfo{author}{W.~Niu}, \bibinfo{author}{Y.~Xu},
\newblock \bibinfo{title}{Uniform boundary estimates in homogenization of
  higher-order elliptic systems},
\newblock \bibinfo{journal}{Annali di Matematica Pura ed Applicata (1923-)}
  \bibinfo{volume}{198} (\bibinfo{year}{2019}) \bibinfo{pages}{97--128}.
\bibitem[{Xu and Niu(2021)}]{R23}
\bibinfo{author}{Y.~Xu}, \bibinfo{author}{W.~Niu},
\newblock \bibinfo{title}{Convergence rates in almost-periodic homogenization
  of higher-order elliptic systems},
\newblock \bibinfo{journal}{Asymptotic Analysis} \bibinfo{volume}{123}
  (\bibinfo{year}{2021}) \bibinfo{pages}{95--137}.
\bibitem[{Huang et~al.(2022)Huang, Xing, and Gao}]{R24}
\bibinfo{author}{Z.~Huang}, \bibinfo{author}{Y.~Xing},
  \bibinfo{author}{Y.~Gao},
\newblock \bibinfo{title}{Two-scale asymptotic homogenization method for
  composite kirchhoff plates with in-plane periodicity},
\newblock \bibinfo{journal}{Aerospace} \bibinfo{volume}{9}
  (\bibinfo{year}{2022}) \bibinfo{pages}{751}.
\bibitem[{Huang et~al.(2025)Huang, Zeng, Wang, and Liu}]{R25}
\bibinfo{author}{Z.~Huang}, \bibinfo{author}{X.~Zeng},
  \bibinfo{author}{C.~Wang}, \bibinfo{author}{C.~Liu},
\newblock \bibinfo{title}{Two-scale convergence analysis and numerical
  simulation for periodic kirchhoff plates},
\newblock \bibinfo{journal}{Heliyon} \bibinfo{volume}{11}
  (\bibinfo{year}{2025}).
\bibitem[{Kolpakov(1999)}]{R26}
\bibinfo{author}{A.~Kolpakov},
\newblock \bibinfo{title}{Variational principles for stiffnesses of a
  non-homogeneous plate},
\newblock \bibinfo{journal}{Journal of the Mechanics and Physics of Solids}
  \bibinfo{volume}{47} (\bibinfo{year}{1999}) \bibinfo{pages}{2075--2092}.
\bibitem[{Lok and Cheng(2000)}]{R27}
\bibinfo{author}{T.-S. Lok}, \bibinfo{author}{Q.-H. Cheng},
\newblock \bibinfo{title}{Elastic stiffness properties and behavior of
  truss-core sandwich panel},
\newblock \bibinfo{journal}{Journal of Structural Engineering}
  \bibinfo{volume}{126} (\bibinfo{year}{2000}) \bibinfo{pages}{552--559}.
\bibitem[{Buannic et~al.(2003)Buannic, Cartraud, and Quesnel}]{R28}
\bibinfo{author}{N.~Buannic}, \bibinfo{author}{P.~Cartraud},
  \bibinfo{author}{T.~Quesnel},
\newblock \bibinfo{title}{Homogenization of corrugated core sandwich panels},
\newblock \bibinfo{journal}{Composite structures} \bibinfo{volume}{59}
  (\bibinfo{year}{2003}) \bibinfo{pages}{299--312}.
\bibitem[{Cai et~al.(2014)Cai, Xu, and Cheng}]{R29}
\bibinfo{author}{Y.~Cai}, \bibinfo{author}{L.~Xu}, \bibinfo{author}{G.~Cheng},
\newblock \bibinfo{title}{Novel numerical implementation of asymptotic
  homogenization method for periodic plate structures},
\newblock \bibinfo{journal}{International Journal of Solids and Structures}
  \bibinfo{volume}{51} (\bibinfo{year}{2014}) \bibinfo{pages}{284--292}.
\bibitem[{Dong et~al.(2016)Dong, Li, Wang, and Wu}]{R30}
\bibinfo{author}{B.~Dong}, \bibinfo{author}{C.~Li}, \bibinfo{author}{D.~Wang},
  \bibinfo{author}{C.-T. Wu},
\newblock \bibinfo{title}{Consistent multiscale analysis of heterogeneous thin
  plates with smoothed quadratic hermite triangular elements},
\newblock \bibinfo{journal}{International Journal of Mechanics and Materials in
  Design} \bibinfo{volume}{12} (\bibinfo{year}{2016})
  \bibinfo{pages}{539--562}.
\bibitem[{Ren et~al.(2017)Ren, Cong, Wang, and Guo}]{R31}
\bibinfo{author}{M.~Ren}, \bibinfo{author}{J.~Cong}, \bibinfo{author}{B.~Wang},
  \bibinfo{author}{X.~Guo},
\newblock \bibinfo{title}{Extended multiscale finite element method for
  small-deflection analysis of thin composite plates with aperiodic
  microstructure characteristics},
\newblock \bibinfo{journal}{Composite Structures} \bibinfo{volume}{160}
  (\bibinfo{year}{2017}) \bibinfo{pages}{422--434}.
\bibitem[{Li et~al.(2025)Li, Sharif~Khodaei, and Aliabadi}]{R32}
\bibinfo{author}{H.~Li}, \bibinfo{author}{Z.~Sharif~Khodaei},
  \bibinfo{author}{M.~Aliabadi},
\newblock \bibinfo{title}{Fft-based homogenisation for efficient concurrent
  multiscale modelling of thin plate structures},
\newblock \bibinfo{journal}{Computational Mechanics} \bibinfo{volume}{75}
  (\bibinfo{year}{2025}) \bibinfo{pages}{637--653}.
\bibitem[{Li et~al.(2005)Li, Soric, Jarak, and Atluri}]{R37}
\bibinfo{author}{Q.~Li}, \bibinfo{author}{J.~Soric},
  \bibinfo{author}{T.~Jarak}, \bibinfo{author}{S.~N. Atluri},
\newblock \bibinfo{title}{A locking-free meshless local petrov--galerkin
  formulation for thick and thin plates},
\newblock \bibinfo{journal}{Journal of Computational Physics}
  \bibinfo{volume}{208} (\bibinfo{year}{2005}) \bibinfo{pages}{116--133}.
\bibitem[{Huang et~al.(2011)Huang, Huang, and Xu}]{R40}
\bibinfo{author}{J.~Huang}, \bibinfo{author}{X.~Huang},
  \bibinfo{author}{Y.~Xu},
\newblock \bibinfo{title}{Convergence of an adaptive mixed finite element
  method for kirchhoff plate bending problems},
\newblock \bibinfo{journal}{SIAM journal on numerical analysis}
  \bibinfo{volume}{49} (\bibinfo{year}{2011}) \bibinfo{pages}{574--607}.
\bibitem[{Wang and Xu(2004)}]{R33}
\bibinfo{author}{L.~Wang}, \bibinfo{author}{X.~Xu}, \bibinfo{title}{The
  Mathematical Basis of Finite Element Method}, \bibinfo{publisher}{Science
  Press}, \bibinfo{year}{2004}.
\bibitem[{Feng and Shi(2010)}]{R34}
\bibinfo{author}{K.~Feng}, \bibinfo{author}{Z.~Shi}, \bibinfo{title}{The
  Mathematical Theory of Elastic Structures}, \bibinfo{publisher}{Science
  Press}, \bibinfo{year}{2010}.
\bibitem[{Cao(2006)}]{R35}
\bibinfo{author}{L.-Q. Cao},
\newblock \bibinfo{title}{Multiscale asymptotic expansion and finite element
  methods for the mixed boundary value problems of second order elliptic
  equation in perforated domains},
\newblock \bibinfo{journal}{Numerische Mathematik} \bibinfo{volume}{103}
  (\bibinfo{year}{2006}) \bibinfo{pages}{11--45}.
\bibitem[{Shen et~al.(2023)Shen, Wu, and Zhu}]{R39}
\bibinfo{author}{X.~Shen}, \bibinfo{author}{R.~Wu}, \bibinfo{author}{S.~Zhu},
\newblock \bibinfo{title}{Numerical method for two-dimensional linearly elastic
  clamped plate model},
\newblock \bibinfo{journal}{International Journal of Computer Mathematics}
  \bibinfo{volume}{100} (\bibinfo{year}{2023}) \bibinfo{pages}{1779--1793}.
\bibitem[{Cui(2001)}]{R36}
\bibinfo{author}{J.~Cui},
\newblock \bibinfo{title}{Multiscale computational method for unified design of
  structure, components and their materials},
\newblock \bibinfo{journal}{Proceedings on Computational Mechanics in Science
  and Engineering, CCCM-2001, Guangzhou}  (\bibinfo{year}{2001})
  \bibinfo{pages}{5--8}.
\bibitem[{Hecht(2012)}]{R38}
\bibinfo{author}{F.~Hecht},
\newblock \bibinfo{title}{New development in freefem++},
\newblock \bibinfo{journal}{Journal of numerical mathematics}
  \bibinfo{volume}{20} (\bibinfo{year}{2012}) \bibinfo{pages}{1--14}.

\end{thebibliography}







\end{document}